\documentclass{amsart} 
\usepackage{float}
\usepackage{amsmath}
\usepackage{epsfig}
\usepackage{siunitx} 
\usepackage{color,graphics}
\usepackage{booktabs}
\usepackage{caption}
\usepackage{subcaption}
\usepackage{amsthm}
\usepackage{amsmath}
\usepackage{amssymb}
\usepackage{lipsum}
\usepackage{pdflscape,longtable,array}
\usepackage{afterpage}
\usepackage{subcaption}
\usepackage{pdflscape}
\usepackage{fancyhdr} 
\usepackage{xcolor}
\usepackage{csquotes}
\usepackage{algorithm}
\usepackage{algpseudocode}

\usepackage{tabularx} 
\usepackage{lineno}
\usepackage{diagbox} 
\numberwithin{equation}{section}
\DeclareMathOperator{\sech}{sech}
\usepackage{xcolor}
\newcolumntype{^}{>{\global\let\currentrowstyle\relax}}
\newcolumntype{_}{>{\currentrowstyle}}
\usepackage{multirow}
\usepackage{epstopdf}
\usepackage[all]{xy}
\usepackage{fixmath, bm, amssymb, amsfonts, latexsym}
\newtheorem{theorem}{Theorem}[section]

\newtheorem{example}[theorem]{Example}

\usepackage[square,numbers]{natbib}
\begin{document}
\author{Lavanya V Salian  $^\dagger$, Samala Rathan$^{\dagger}$, and Debojyoti Ghosh$^{\S}$}
\title{A novel central compact finite-difference scheme for third derivatives with high spectral resolution}
\thanks{
$^\dagger$Department of Humanities and Sciences, Indian Institute of Petroleum and Energy-Visakhapatnam,
India-530003 ({lavanya\_vs@iipe.ac.in, rathans.math@iipe.ac.in (Corresponding author: Samala Rathan)}),
\newline
$^{\S}$ Center for Applied Scientific Computing, Lawrence Livermore National Laboratory, Livermore, CA 94550, United States  ({ghosh5@llnl.gov})}

\date{\today}

\maketitle

\begin{abstract}
In this paper, we introduce a novel category of central compact schemes inspired by existing cell-node and cell-centered compact finite difference schemes, that offer a superior spectral resolution for solving the dispersive wave equation. In our approach, we leverage both the function values at the cell nodes and cell centers to calculate third-order spatial derivatives at the cell nodes. To compute spatial derivatives at the cell centers, we employ a technique that involves half-shifting the indices within the formula initially designed for the cell-nodes. In contrast to the conventional compact interpolation scheme, our proposed method effectively sidesteps the introduction of transfer errors. We employ the Taylor-series expansion-based method to calculate the finite difference coefficients. By conducting systematic Fourier analysis and numerical tests, we note that the methods exhibit exceptional characteristics such as high order, superior resolution, and low dissipation. Computational findings further illustrate the effectiveness of high-order compact schemes, particularly in addressing problems with a third derivative term.
\end{abstract}



\pagestyle{myheadings} \thispagestyle{plain} \markboth{Lavanya V Salian, Samala Rathan, Debojyoti Ghosh}{Central compact schemes for third derivatives}

\section{Introduction}
\label{sec:1}
In this paper, we propose  a new central compact scheme with enhanced spectral resolution to approximate the third-order derivative term that evolves the nonlinear, possibly prototypical, dispersion equation of the form
\begin{equation} 
\begin{split}
u_t+g(u)_x + f(u)_{xxx}&=0, \quad (x,t) \in \Omega \times (0,T],\\
u(x,0) &= u_0 (x),
\end{split}
\label{eqn:A}
\end{equation}
where $g(u)$ and $f(u)$ are arbitrary nonlinear functions that may be smooth. The Korteweg–de Vries (KdV) equation  is nonlinear dispersive partial differential equation (PDE) that represents several physical phenomena. Examples include the shallow water waves and plasma ion acoustic waves \cite{guo2018study} and long-distance propagation in optical fibers \cite{Iordache2010study}. The KdV equation governs weakly--nonlinear long waves and is a fundamental nonlinear evolution equation of the form
\begin{equation}
\label{KdV}
u_t+a u u_x + u_{xxx} =0, \quad (x,t) \in \Omega \times (0,T],\\
\end{equation}
 where $u$ is a function of $(x,t)$, and $u_{xxx}$  is the linear dispersion term. Solitary waves are localized traveling waves characterized by a consistent speed and shape that gradually diminish to zero amplitude at a significant distance \cite{wazwaz2010partial}. The pioneering investigations by Zabusky and Kruskal \cite{zabusky1965interaction} reveal interesting characteristics of solitary waves and their involvement in nonlinear interactions. These lead to the emergence of waves that preserve both their original shape and amplitude. Such unique solitary waves are commonly referred to as \textit{solitons}. The KdV equation elegantly captures the dynamics of waves in nonlinear, dispersive media with applications in multiple fields. These include semiconductor device simulations \cite{gardner1994quantum}, aeroacoustics~\cite{hu1996low, hixon2000compact}, electromagnetic simulations~\cite{shang1999high}, and tectonic dynamics~\cite{lakshmanan2011solitons}. Challenges arise in \enquote{convection--dominated} scenarios categorized, where third derivative terms may have small or even zero coefficients. Stable, efficient, and high-order numerical methods for such applications remain a persistent challenge. Notably, in applications like aeroacoustics and electromagnetism, researchers have employed high-order compact finite difference schemes to minimize error accumulation for accurate simulations of wave propagation over extended periods~\cite{Lele, JV, hixon2000compact}.
\par
The significance of the dispersive KdV equation and its applications has led to a variety of analytical and numerical methods. These approaches include the differential quadrature method~\cite{karunakar2019differential}, the inverse scattering technique~\cite{ablowitz1991solitons}, and the variational iteration method (VIM) employed by Wazwaz~\cite{wazwaz2007variational} to address problems related to the Burgers, cubic Boussinesq, KdV, and K$\left(2,2\right)$ equations. Additionally, the local discontinuous Galerkin method~\cite{LSY}, the adaptive mesh refinement (AMR)--based line method~\cite{SWSZ}, and high-order compact schemes with high-order low-pass filters~\cite{JV} have been successfully applied. Solutions to the dispersive wave equation~(\ref{eqn:A}) share similarities with those of hyperbolic conservation laws, such as the potential existence of sharp fronts and the finite wave propagation speeds, Ahmat and Qiu~\cite{AQ} applied a fifth--order weighted essentially non--oscillatory (WENO) scheme utilizing polynomial bases. An alternative approach by Lavanya and Rathan~\cite{LS} involves exponential basis WENO reconstruction, incorporating a tension parameter that can be adjusted according to the given data.
\par
High-order finite difference (FD) methods can be divided into explicit and compact or Pad\'{e}-type schemes. Explicit schemes compute numerical derivatives directly at grid points using stencils whose size typically increases with the order of accuracy. On the other hand, compact schemes use smaller stencils but solve linear systems of equations to determine derivatives along grid lines. Although compact schemes are more accurate than explicit schemes of the same order, they come with the added complexity inverting a scalar tridiagonal or pentadiagonal matrix. Lele's seminal work~\cite{Lele} influenced the field of compact schemes, focusing on derivatives, interpolation, and filtering. Systematic Fourier analysis reveals the spectral--like resolution of Lele's compact schemes, particularly for short waves, surpassing the resolution of cell node compact schemes. Staggered compact schemes developed by Nagarajan et al. \cite{nagarajan2003robust} and Boersma \cite{boersma2005staggered} exhibit robustness in numerical tests. Despite the perceived advantage of a lower aliasing error in staggered compact schemes, challenges arise due to the inclusion of cell-centered values, necessitating interpolation from grid nodes (cell boundaries) as articulated by Lele~\cite{Lele}. To mitigate the transfer error introduced by the interpolation, Zhang et al.~\cite{zhang2008development} devised a shock--capturing weighted compact scheme using a weighted interpolation method based on the WENO idea~\cite{liu1994weighted, jiang1996efficient}. Many enhanced versions of the WENO schemes \cite{balsara2000monotonicity, henricketal2005, martin2006bandwidth, borgesetal2008, yamaleevcarpenter2009_1, zhu2010trigonometric, ha2013modified, ha2016sixth, fuetal2016, rathan2018modified, rathan2020simple, rathan2017improved, rathan2020l1, abedian2022high, rathan2023sixth} and nonlinear compact schemes~\cite{cockburn1994nonlinearly, weinan1996essentially, adams1996high, deng1997compact, gaitonde1997optimized, yee1997explicit, ekaterinaris1999implicit, deng2000developing, lee2002new, pirozzoli2002conservative, sengupta2003analysis, ghosh2012compact, pengshen2015, pengshen2017, fidalgoetal2018, subramaniametal2019, chenetal2021, hiejima2022} have been reported in the literature. Liu et al.'s~\cite{LZZS} approach for first derivatives combines cell-centered and cell node values on the right-hand side of the compact scheme, enhancing accuracy, order, and wave resolution properties without additional computation cost, albeit with increased memory requirements. Leveraging prior work, Wang et al.~\cite{WLWXC} innovatively broadened the application of weighted summation to incorporate second-order spatial derivatives within the acoustic wave equation. Their approach ingeniously equates these summations, calculated both on cell nodes and centers, enabling Taylor series expansion and optimization-based methods. 
\par 
In this paper, a novel central compact scheme with an enhanced spectral resolution for the dispersive wave equation is developed by extending the concepts from existing cell-node \cite{Lele} and cell-centered compact finite difference schemes. In this approach, we utilize function values at both the cell-nodes and cell-centers to calculate third-order spatial derivatives at the cell-nodes.  Furthermore, the cell-centered values are treated as independent variables and are evolved with the cell-node values. Spatial derivatives at the cell-centers are determined by ``half-shifting'' the formula initially devised for the cell-nodes. While this method increases the memory required, there is no corresponding increase in computational cost. This is because the compact interpolation used for calculating values on the half-grid is substituted with a compact formula for computing spatial derivatives (and the updating the residual) at these half-grid points; both incur comparable computational costs. In comparison to the traditional compact interpolation method, this proposed approach effectively eliminates the introduction of transfer errors. The finite difference coefficients are computed using either truncation error-based or least-squares-based methods optimized for superior spectral properties. The third derivative central compact scheme (TDCCS) is shown to have superior spectral properties as compared to the third derivative cell node compact scheme (TDCNCS) and third derivative cell centered compact scheme (TDCCCS)  of the same order (or even higher). Numerical experiments are conducted and compared with the TDCNCS. 
\par
The structure of this paper is outlined as follows:  Section~\ref{sec:cncccs} initiates with a review of third derivative cell-node and cell-centered compact finite difference (FD) schemes.   In Section~\ref{sec:new_ccs}, we present the derivation of the new scheme, providing detailed insights into the approaches for determining the FD coefficients. In Section~\ref{sec:Fourier}, a Fourier analysis is conducted to systematically evaluate the wave resolution of the proposed schemes. Section~\ref{sec:Filter} briefly introduces the high-order central explicit filtering scheme employed to mitigate numerical oscillations.  Section~\ref{sec:time_integration} presents time integration and linear stability analysis. In Section~\ref{sec:num_ex}, numerical examples are presented to validate the advantages of the proposed method. Lastly, Section~\ref{sec:con} provides concluding remarks.
\section{Cell-node and cell-centered compact schemes}
\label{sec:cncccs}
We begin by examining Lele's cell-node compact scheme, initially designed to achieve accuracy up to sixth-order in approximating the third derivative. In this section, we expand upon Lele's work in two aspects. First, we enhance the cell-node compact scheme, achieving an improved tenth-order accuracy. Secondly, we extend it to the cell-centered compact scheme with an accuracy of up to tenth-order.
\par
\begin{figure}[htbp!]
\begin{center}
\minipage{0.9\textwidth}
  \includegraphics[width=\linewidth]{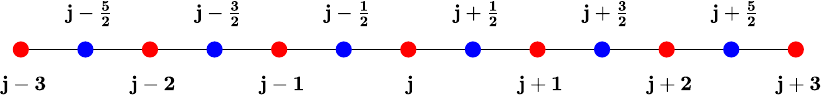}
\endminipage\hfill
\caption{\scriptsize{ The stencil of cell-center and cell-node compact schemes. The cell nodes and cell-centers are denoted by the red circles and blue circles, respectively.}}\label{Fig:F_1}
\end{center}
\end{figure}

We consider numerical approximations to the 1D prototypical dispersion equations of the form
\begin{equation}
\label{eqn:1a}
\frac{\partial u}{\partial t} + \frac{\partial g(u)}{\partial x} + \frac{\partial^3 f(u)}{\partial x^3} =0.
\end{equation}
The framework for describing a semidiscrete finite difference is given by
\begin{equation}
\label{eqn:1b}
\frac{d u_j}{d t} = - g'_j -f'''_j.
\end{equation}
Here, $g'_j$ and $f'''_j$ represent approximations to the spatial derivatives $\frac{\partial g(u)}{\partial x}$ and $\frac{\partial^3 f(u)}{\partial x^3}$ at the grid node $x_j$ respectively. The computational domain is discretized uniformly into $N$ points: $x_1, x_2, ..., x_{j-1}, x_j, x_{j+1}, ..., x_N$. The mesh size is denoted as $h=\Delta x = x_{j+1} - x_j$. Figure \ref{Fig:F_1} illustrates the stencil of cell nodes and the cell centers. For the computation of the first derivatives $g'_j$, the cell node compact scheme introduced by Lele \cite{Lele} and the central compact scheme proposed by Liu et al. \cite{LZZS} are referenced. 
The linear cell-node compact scheme with up to sixth-order accuracy \cite{Lele} is given by
\begin{equation}
\label{eqn:2a}
\begin{split}
\beta f'''_{j-2}+\alpha f'''_{j-1}+f'''_{j}+\alpha f'''_{j+1}+\beta f'''_{j+2} &= a \frac{f_{j+2}-2f_{j+1}+2f_{j-1} - f_{j-2}}{2h^3}+b \frac{f_{j+3}-3f_{j+1}+3f_{j-1} - f_{j-3}}{8h^3}.\\ 
\end{split}
\end{equation}
\par We have extended the expression (\ref{eqn:2a}) to achieve up to tenth-order accuracy, denoted as the third derivative cell-node compact scheme (TDCNCS). The general form is given by
\begin{equation}
\label{eqn:2}
\begin{split}
\beta f'''_{j-2}+\alpha f'''_{j-1}+f'''_{j}+\alpha f'''_{j+1}+\beta f'''_{j+2} &= a \frac{f_{j+2}-2f_{j+1}+2f_{j-1} - f_{j-2}}{2h^3}+b \frac{f_{j+3}-3f_{j+1}+3f_{j-1} - f_{j-3}}{8h^3}\\ 
&+c \frac{f_{j+4}-4f_{j+1}+4f_{j-1} - f_{j-4}}{20h^3}.
\end{split}
\end{equation}
The third derivative cell-centered compact scheme (TDCCCS) is given by
\begin{equation}
\label{eqn:3}
\begin{split}
\beta f'''_{j-2}+\alpha f'''_{j-1}+f'''_{j}+\alpha f'''_{j+1}+\beta f'''_{j+2}&= a \frac{f_{j+\frac{3}{2}}-3f_{j+\frac{1}{2}}+3f_{j-\frac{1}{2}} -f_{j-\frac{3}{2}}}{h^3}+b \frac{f_{j+\frac{5}{2}}-5f_{j+\frac{1}{2}}+5f_{j-\frac{1}{2}} -f_{j-\frac{5}{2}}}{5h^3}\\
&+c \frac{f_{j+\frac{7}{2}}-7f_{j+\frac{1}{2}}+7f_{j-\frac{1}{2}} -f_{j-\frac{7}{2}}}{14h^3}.
\end{split}
\end{equation}
\begin{table}[!ht]
\centering
\captionof{table}{The coefficients of TDCNCS schemes.}\label{Table:T_1a}
\setlength{\tabcolsep}{0pt}
\begin{tabular*}{\textwidth}{@{\extracolsep{\fill}} l *{9}{c} }
\toprule
\textbf{Scheme} &
\textbf{a} &
\textbf{b} &
\textbf{c} &
\textbf{$\alpha$} &
\textbf{$\beta$} &
\textbf{Order}\\
\midrule
TDCNCS-E2 & 1 & 0 & 0 & 0 & 0 & 2\\ [0.15cm]
TDCNCS-E4 & 2 & -1 & 0 & 0 & 0 & 4\\ [0.15cm]
TDCNCS-E6 & $\frac{169}{60}$ & $-\frac{12}{5}$ & $\frac{7}{12}$ & 0 & 0 & 6\\ [0.15cm]
TDCNCS-T4 & 2 & 0 & 0 & $\frac{1}{2}$ & 0 & 4\\ [0.15cm]
TDCNCS-T6 & 2 & $-\frac{1}{8}$ & 0 & $\frac{7}{16}$ & 0 & 6\\ [0.15cm]
TDCNCS-T8 & $\frac{2367}{1180}$ & $-\frac{167}{1180}$ & $\frac{1}{236}$ & $-\frac{205}{472}$ & 0 & 8\\ [0.15cm]
TDCNCS-P6 & $\frac{40}{21}$ & 0 & 0 & $\frac{4}{9}$ & $\frac{1}{126}$ & 6\\ [0.15cm]
TDCNCS-P8 & $\frac{160}{83}$ & $-\frac{5}{166}$ & 0 & $\frac{147}{332}$ & $-\frac{1}{166}$ & 8\\ [0.15cm]
TDCNCS-P10 & $\frac{18221}{5478}$ & $-\frac{1846}{913}$ & $\frac{5}{66}$ & $\frac{799}{2739}$ & $-\frac{557}{5478}$ & 10\\ [0.15cm]
\bottomrule
\end{tabular*}
\end{table}
\begin{table}[!ht]
\centering
\captionof{table}{The coefficients of TDCCCS schemes.}\label{Table:T_1b}
\setlength{\tabcolsep}{0pt}
\begin{tabular*}{\textwidth}{@{\extracolsep{\fill}} l *{9}{c} }
\toprule
\textbf{Scheme} &
\textbf{a} &
\textbf{b} &
\textbf{c} &
\textbf{$\alpha$} &
\textbf{$\beta$} &
\textbf{Order}\\
\midrule
TDCCCS-E2 & 1 & 0 & 0 & 0 & 0 & 2\\ [0.15cm]
TDCCCS-E4 & $\frac{13}{8}$ & $-\frac{5}{8}$ & 0 & 0 & 0 & 4\\ [0.15cm]
TDCCCS-E6 & $\frac{1299}{640}$ & $-\frac{499}{384}$ & $\frac{259}{960}$ & 0 & 0 & 6\\ [0.15cm]
TDCCCS-T4 & $\frac{4}{3}$ & 0 & 0 & $\frac{1}{6}$ & 0 & 4\\ [0.15cm]
TDCCCS-T6 & $\frac{205}{166}$ & $\frac{35}{166}$ & 0 & $\frac{37}{166}$ & 0 & 6\\ [0.15cm]
TDCCCS-T8 & $\frac{1058279}{975200}$ & $\frac{96627}{195040}$ & $-\frac{24787}{487600}$ & $\frac{3229}{12190}$ & 0 & 8\\ [0.15cm]
TDCCCS-P6 & $\frac{320}{233}$ & 0 & 0 & $\frac{134}{699}$ & $-\frac{7}{1398}$ & 6\\ [0.15cm]
TDCCCS-P8 & $\frac{49720}{79903}$ & $\frac{91400}{79903}$ & 0 & $\frac{28838}{79903}$ & $\frac{3541}{159806}$ & 8\\ [0.15cm]
TDCCCS-P10 & $\frac{55463611}{150617762}$ & $\frac{677644345}{451853286}$ & $\frac{6301771}{225926643}$ & $\frac{93443398}{225926643}$ & $\frac{15505921}{451853286}$ & 10\\ [0.15cm]
\bottomrule
\end{tabular*}
\end{table}
\par
The left-hand sides of both equations (\ref{eqn:2}) and (\ref{eqn:3}) involve spatial derivatives $f'''_j$ computed at the grid nodes. While the right-hand side of equation (\ref{eqn:2}) relies solely on function values $f'''_j$ at the grid node $x_j$, equation (\ref{eqn:3}) incorporates function values $f_{j+\frac{1}{2}}$ at the center $x_{j+\frac{1}{2}}=\frac{1}{2}(x_j + x_{j+1})$ within each interval $x \in [x_j , x_{j+1}]$. The accuracy of these schemes is contingent upon specific choices for the coefficients $(\alpha, \beta, a, b, c)$. By matching terms in the Taylor series expansion around the point $x_j$, we can derive conditions for achieving different orders of accuracy. Tables (\ref{Table:T_1a}) and (\ref{Table:T_1b}) explicitly lists the coefficients for the TDCNCS and TDCCCS schemes, respectively. By restricting the parameter $\alpha = \beta = 0$, we obtain a family of explicit schemes. Further, upon setting $\alpha \neq 0$ and $\beta = 0$ yields tridiagonal schemes. The combination of $\alpha \neq 0$ and $\beta \neq 0$ produces pentadiagonal schemes. These three distinct categories are referred to as TDCNCS-E, TDCNCS-T, and TDCNCS-P, respectively. Their formal order of accuracy is appended to their names for convenient identification. The truncation error for eighth-order accuracy is expressed as $Q f_j^{(11)}(x) h^8 + \mathcal{O}(h^{10})$, where $Q$ is $3.12192 \times 10^{-5}$ for equation (\ref{eqn:2}) and $6.57252 \times 10^{-5}$ for equation (\ref{eqn:3}).
\par Fourier analysis serves as an effective means for evaluating the accuracy and resolution characteristics of a FD schemes. The modified wavenumbers associated with the third-order spatial derivative for equation (\ref{eqn:2}) and (\ref{eqn:3}) are
\begin{equation}
    \omega'''_{\text{TDCNCS}} = \frac{a[2 \sin(\omega) - \sin(2\omega)]+\frac{b}{4}[3 \sin(\omega) - \sin(3\omega)]+\frac{c}{10}[4 \sin(\omega) - \sin(4\omega)]}{1+2\alpha \cos(\omega)+2\beta \cos(2\omega)},
\end{equation}

\begin{equation}
    \omega'''_{\text{TDCCCS}} = \frac{2a[3 \sin(\frac{\omega}{2}) - \sin(\frac{3\omega}{2})]+\frac{2b}{5}[5 \sin(\frac{\omega}{2}) - \sin(\frac{5\omega}{2})]+\frac{c}{7}[7 \sin(\frac{\omega}{2}) - \sin(\frac{7\omega}{2})]}{1+2\alpha \cos(\omega)+2\beta \cos(2\omega)}.
\end{equation}
Here, the parameter $\omega = kh$ represents the scaled wavenumber, and $\omega'''$ signifies the scaled modified wavenumber.  Exact difference approximation occurs when the scaled wavenumber $\omega$ and the scaled modified wavenumber  $\omega'''$  coincide, as represented by the equation $\omega = \omega'''$.
\begin{figure}[htbp!]
  \begin{minipage}[b]{0.3\linewidth}
    \centering
    \includegraphics[trim=0cm 0cm 0cm 0cm, clip=true,width=\linewidth]{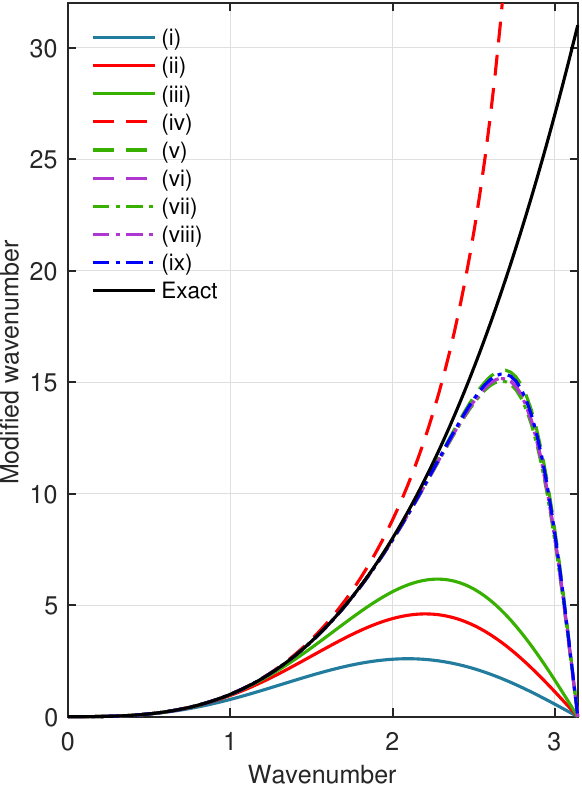}
    \subcaption*{\normalsize{\centering (a) TDCNCS}}
  \end{minipage}\hfill
    \begin{minipage}[b]{0.3\linewidth}
    \centering
    \includegraphics[trim=0cm 0cm 0cm 0cm, clip=true,width=\linewidth]{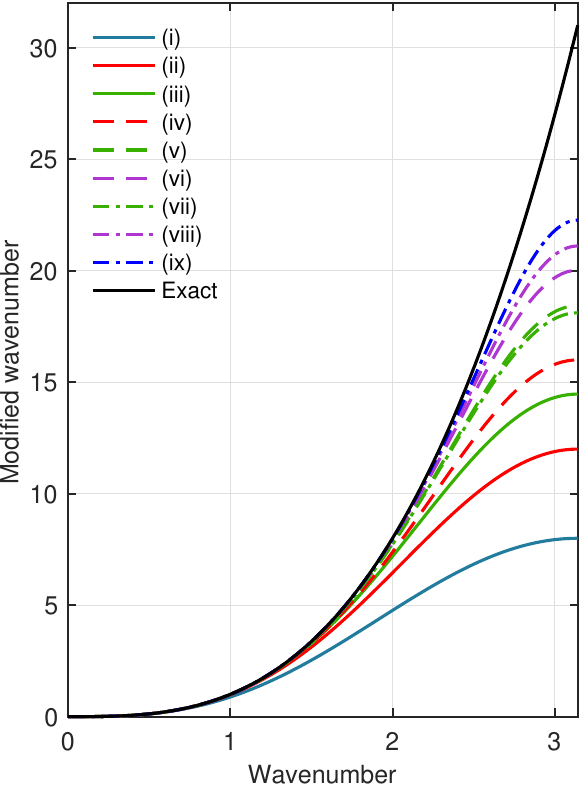}
    \subcaption*{\normalsize{\centering (b) TDCCCS}}
  \end{minipage}\hfill
  \begin{minipage}[b]{0.3\linewidth}
    \centering
   \includegraphics[trim=0cm 0cm 0cm 0cm, clip=true,width=\linewidth]{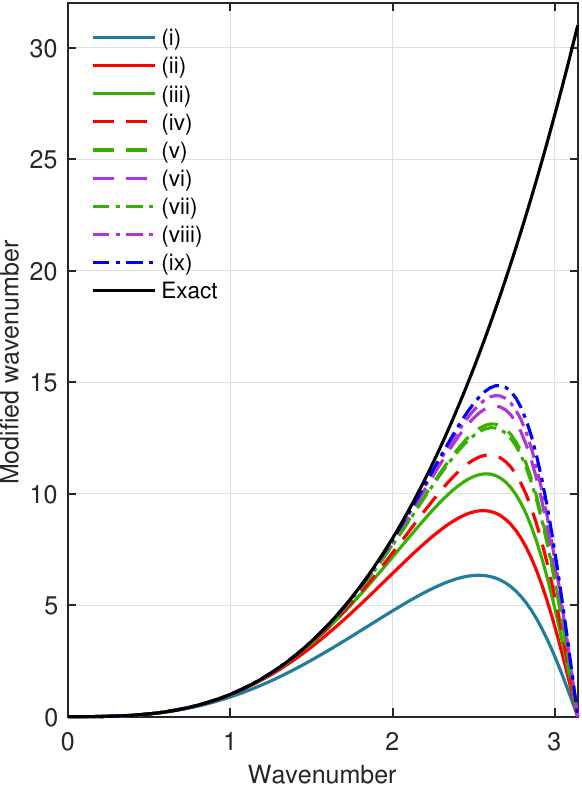}
       \subcaption*{\normalsize{\centering (c) TDCCCS-CI}}
  \end{minipage} 
  \caption{ Plot of modified wavenumber versus wavenumber for (a) TDCNCS, (b) TDCCCS and (c) TDCCCS-CI : (i) E2; (ii) E4; (iii) E6; (iv) T4; (v) T6; (vi) T8; (vii) P6; (viii) P8; (ix) P10. The tenth-order pentadiagonal scheme is adopted for the compact interpolation in TDCCCS-CI.}\label{Fig:F_2}
\end{figure}
\par Figure \ref{Fig:F_2} illustrates the plot of modified wavenumber versus wavenumber, representing the resolution properties for various schemes: (a) TDCNCS, (b) TDCCCS. The schemes include (i) second-order explicit scheme (E2), (ii) fourth-order explicit scheme (E4), (iii) sixth-order explicit scheme (E6), (iv) fourth-order tridiagonal scheme (T4), (v) sixth-order tridiagonal scheme (T6), (vi) eighth-order tridiagonal scheme (T8), (vii) sixth-order pentadiagonal scheme (P6), (viii) eighth-order pentadiagonal scheme (P8), and (ix) tenth-order pentadiagonal scheme (P10). Notably, for a same operator length, TDCCCS exhibits significantly enhanced resolution compared to TDCNCS. However, the implementation of TDCCCS requires the estimation of function values at the cell centers. The high-order compact interpolation (CI) scheme represents the most straightforward method to perform this estimation,
\begin{equation}
\label{eqn:4a}
\beta \hat{f}_{j-\frac{3}{2}}+\alpha \hat{f}_{j-\frac{1}{2}}+\hat{f}_{j+\frac{1}{2}}+\alpha \hat{f}_{j+\frac{3}{2}}+\beta \hat{f}_{j+\frac{5}{2}}= a\frac{(f_{j+1}+f_j)}{2}+b\frac{(f_{j+2}+f_{j-1})}{2}+c\frac{(f_{j+3}+f_{j-2})}{2},
\end{equation}
where $\hat{f}_{j+\frac{1}{2}}$ represent the interpolated values at the midpoints $(j+\frac{1}{2})$. The transfer function corresponding to (\ref{eqn:4a}) is
\begin{equation}
\label{eqn:4b}
T_{CI}(\omega)=  \frac{a \cos(\frac{\omega}{2})+b \cos(\frac{3\omega}{2})+c \cos(\frac{5\omega}{2})}{1+2\alpha \cos(\omega)+2\beta \cos(2\omega)}.
\end{equation}
The coefficients in (\ref{eqn:4a}) can be determined by matching the truncated expansion (TE) coefficients corresponding to different orders of accuracy, which are listed in Table (\ref{Table:T_1c}). We use the notation TDCCCS-CI to represent TDCCCS combined with the CI scheme. To achieve greater precision in the interpolated function values at the midpoints, a larger interpolation stencil is required, leading to an increased computational cost. In this study, we employ the tenth-order penta-diagonal CI scheme to calculate function values at the cell-centers within the framework of TDCCCS-CI. Figure \ref{Fig:F_2}(c) illustrates the dispersion relations of TDCCCS-CI for various accuracy levels. When compared to Figure \ref{Fig:F_2}(b), it becomes evident that the use of CI results in a noticeable reduction in resolution for the TDCCCS, as it introduces transfer errors.
\begin{table}[!ht]
\centering
\captionof{table}{The coefficients of the transfer function.}\label{Table:T_1c}
\setlength{\tabcolsep}{0pt}
\begin{tabular*}{\textwidth}{@{\extracolsep{\fill}} l *{9}{c} }
\toprule
\textbf{Scheme} &
\textbf{a} &
\textbf{b} &
\textbf{c} &
\textbf{$\alpha$} &
\textbf{$\beta$} &
\textbf{Order}\\
\midrule
CI-E2 & 1 & 0 & 0 & 0 & 0 & 2\\ [0.15cm]
CI-E4 & $\frac{9}{8}$ & -$\frac{1}{8}$ & 0 & 0 & 0 & 4\\ [0.15cm]
CI-E6 & $\frac{75}{64}$ & $-\frac{25}{128}$ & $\frac{3}{128}$ & 0 & 0 & 6\\ [0.15cm]
CI-T4 & $\frac{4}{3}$ & 0 & 0 & $\frac{1}{6}$ & 0 & 4\\ [0.15cm]
CI-T6 & $\frac{3}{2}$ & $\frac{1}{10}$ & 0 & $\frac{3}{10}$ & 0 & 6\\ [0.15cm]
CI-T8 & $\frac{25}{16}$ & $\frac{5}{32}$ & $-\frac{1}{224}$ & $\frac{5}{14}$ & 0 & 8\\ [0.15cm]
CI-P6 & $\frac{64}{45}$ & 0 & 0 & $\frac{2}{9}$ & $-\frac{1}{90}$ & 6\\ [0.15cm]
CI-P8 & $\frac{8}{5}$ & $\frac{8}{35}$ & 0 & $\frac{2}{5}$ & $\frac{1}{70}$ & 8\\ [0.15cm]
CI-P10 & $\frac{5}{3}$ & $\frac{5}{14}$ & $\frac{1}{126}$ & $\frac{10}{21}$ & $\frac{5}{126}$ & 10\\ [0.15cm]
\bottomrule
\end{tabular*}
\end{table}
\section{A new class of central compact schemes}\label{sec:new_ccs}
This section introduces the methodology for designing third derivative central compact schemes (TDCCS). The stencil in the cell-centered compact schemes, as defined by equation (\ref{eqn:3}), includes both grid points and half-grid points, denoted as $\{ j-\frac{5}{2},j-2,j-\frac{3}{2},j-1,j-\frac{1}{2}, j, j+\frac{1}{2},j+1,j+\frac{3}{2},j+2, j+\frac{5}{2}\}$. However, only the values corresponding to the cell centers, specifically $\{ j-\frac{5}{2},j-\frac{3}{2},j-\frac{1}{2}, j+\frac{1}{2},j+\frac{3}{2},j+\frac{5}{2}\}$, are utilized in computing derivatives at the cell nodes $\{ j-2, j-1, j, j+1, j+2 \}$. Employing both the values at the cell nodes $\{ j-2, j-1, j, j+1, j+2 \}$ and the cell centers $\{ j-\frac{5}{2},j-\frac{3}{2},j-\frac{1}{2}, j+\frac{1}{2},j+\frac{3}{2},j+\frac{5}{2}\}$ could potentially result in a compact scheme with increased order accuracy and improved resolution~\cite{LZZS,WLWXC}. Inspired by this concept, we propose a new category of third derivative central compact schemes (TDCCS) represented by the following formula:
\begin{equation}
\label{eqn:5}
\begin{split}
\beta f'''_{j-2}+\alpha f'''_{j-1}+f'''_{j}+\alpha f'''_{j+1}+\beta f'''_{j+2} &=  a \frac{4f_{j+1}-8f_{j+\frac{1}{2}}+8f_{j-\frac{1}{2}} -4f_{j-1}}{h^3}\\&+b  \frac{8f_{j+\frac{3}{2}}-12f_{j+1}+12f_{j-1} -8f_{j-\frac{3}{2}}}{5h^3}\\ &+ c  \frac{8f_{j+\frac{5}{2}}-20f_{j+1}+20f_{j-1} -8f_{j-\frac{5}{2}}}{35h^3}.\\
\end{split}
\end{equation}
Note that in equation (\ref{eqn:5}), we must initially compute the function values at the cell-centers. These values at the cell-centers can be determined using the high-order CI method, as described in equation (\ref{eqn:4a}). Nevertheless, as previously discussed in the preceding section, employing the high-order CI approach may introduce transfer errors that undermine the precision of TDCCS. To ensure the accuracy of the cell-center values, the values at the cell centers are stored as independent computational variables, and the identical scheme is employed to compute updating values on cell nodes. This involves a straightforward approach by shifting the indices in equation (\ref{eqn:5}) by $\frac{1}{2}$. 
\begin{equation}
\label{eqn:6}
\begin{split}
\beta f'''_{j-\frac{5}{2}}+\alpha f'''_{j-\frac{3}{2}}+f'''_{j-\frac{1}{2}}+\alpha f'''_{j+\frac{1}{2}}+\beta f'''_{j+\frac{3}{2}} &=  a \frac{4f_{j+\frac{1}{2}}-8f_{j}+8f_{j-1} -4f_{j-\frac{3}{2}}}{h^3}\\&+b \frac{8f_{j+1}-12f_{j+\frac{1}{2}}+12f_{j-\frac{3}{2}} -8f_{j-2}}{5h^3}\\ &+ c  \frac{8f_{j+2}-20f_{j+\frac{1}{2}}+20f_{j-\frac{3}{2}} -8f_{j-3}}{35h^3}.\\
\end{split}
\end{equation}
It is important to observe that this modification results in a higher memory demand for storing function values at cell centers.  However, there is no corresponding increase in computational cost, as the compact interpolation (\ref{eqn:4a}) is substituted with the comparable-cost compact updating (\ref{eqn:6}). Both equations (\ref{eqn:5}) and (\ref{eqn:6}), for the same accuracy order, utilize an identical set of coefficients. These coefficients, denoted by $\alpha, \beta, a, b $ and $c$  can be determined through two approaches: either by matching the TE coefficients for different accuracy levels, or by optimizing a misfit function. Equation (\ref{eqn:5}) was chosen through analysis and comparison with alternative combinations outlined in the Appendix~\ref{Appendix:A}. Among the four combinations considered, section ~\ref{sec:Fourier} demonstrates that equation (\ref{eqn:5}) offers superior spectral resolutions.
\subsection{Determining FD Coefficients Based on TE:}
To derive the relationships among the coefficients $a, b, c, \alpha,$ and $\beta$ in equation (\ref{eqn:5}), we match the Taylor series coefficients of different orders. Solving the resulting set of linear equations yields schemes ranging from second to tenth orders. The relationships for different orders are presented as follows:\\
 Second order:
\begin{equation}
\label{eqn:5a}
    1+2 \alpha+ 2 \beta = a+b+c
\end{equation}
Fourth order:
\begin{equation}
\label{eqn:5b}
   \alpha + 2^2 \beta=\frac{a}{16} +\frac{13b}{80}+\frac{29 c}{80}
\end{equation}
Sixth order:
\begin{equation}
\label{eqn:5c}
    \alpha + 2^4 \beta=\frac{13a}{160} +\frac{93b}{160}+\frac{2451 c}{1120}
\end{equation}
 Eighth order:
\begin{equation}
\label{eqn:5d}
    \alpha + 2^6 \beta=\frac{205a}{2688} +\frac{4069b}{2688}+\frac{30025 c}{2688}
\end{equation}
Tenth order:
\begin{equation}
\label{eqn:5e}
   \alpha + 2^8 \beta=\frac{671a}{7680} +\frac{36991b}{7680}+\frac{534991 c}{7680}
\end{equation}
By solving equations (\ref{eqn:5a})–(\ref{eqn:5e}), we can determine the coefficients for TDCCS. When the schemes are constrained to ~$\alpha=\beta=0$, they yield an explicit family of TDCCS. Conversely, if the schemes are limited to $\alpha \ne 0$, various tridiagonal TDCCS are derived. Moreover, when both $\alpha \ne 0$ and $\beta \ne 0$, the result is a family of pentadiagonal TDCCS. We label these three distinct types of schemes as TDCCS-E, TDCCS-T, and TDCCS-P, respectively. For ease of identification and unambiguous referencing, the formal order of accuracy of each scheme type is appended to its respective acronym. Table (\ref{Table:T_2}) presents the coefficients for TDCCS. The CCS-T6 and CCS-T8 schemes strike a remarkable balance between resolution, accuracy, and efficiency. They achieve sixth- and eighth-order accuracy, respectively, while maintaining a tridiagonal matrix structure, leading to significant computational savings compared to pentadiagonal schemes with equivalent accuracy.  For the later numerical analysis, we are considering eighth-order TDCCS. The truncation error for eighth-order accuracy TDCCS (\ref{eqn:5}) is given by $2.1882 \times 10^{-6} f_j^{(11)}(x) h^8 + \mathcal{O}(h^{10}) $. Note that the magnitude of the leading error term in the TDCCS scheme is an order of magnitude lower than the corresponding TDCNCS scheme of the same order. We denote TDCCS with the truncated expansion coefficients as TDCCS-TE.
\begin{table}[!ht]
\centering
\captionof{table}{The coefficients of TDCCS schemes.}\label{Table:T_2}
\setlength{\tabcolsep}{0pt}
\begin{tabular*}{\textwidth}{@{\extracolsep{\fill}} l *{9}{c} }
\toprule
\textbf{Scheme} &
\textbf{a} &
\textbf{b} &
\textbf{c} &
\textbf{$\alpha$} &
\textbf{$\beta$} &
\textbf{Order}\\
\midrule
TDCCS-E4 & $\frac{13}{8}$ & $-\frac{5}{8}$ & 0 & 0 & 0 & 4\\ [0.15cm]
TDCCS-E6 & $\frac{361}{192}$ & $-\frac{129}{128}$ & $\frac{49}{384}$ & 0 & 0 & 6\\ [0.15cm]
TDCCS-T4 & $\frac{8}{7}$ & 0 & 0 & $\frac{1}{14}$ & 0 & 4\\ [0.15cm]
TDCCS-T6 & 5 & -5 & 0 & $-\frac{1}{2}$ & 0 & 6\\ [0.15cm]
TDCCS-T8 & $\frac{58021}{14120}$ & $-\frac{109007}{28240}$ & $\frac{1029}{28240}$ & $-\frac{1261}{3530}$ & 0 & 8\\ [0.15cm]
TDCCS-P6 & $\frac{320}{273}$ & 0 & 0 & $\frac{74}{819}$ & $-\frac{1}{234}$ & 6\\ [0.15cm]
TDCCS-P8 & $\frac{19640}{4621}$ & $-\frac{353000}{87799}$ & 0 & $-\frac{33746}{87799}$ & $-\frac{147}{175598}$ & 8\\ [0.15cm]
TDCCS-P10 & $\frac{74390155}{19635801}$ & $-\frac{45752035}{13090534}$ & $\frac{4684435}{39271602}$ & $-\frac{5803114}{19635801}$ & $\frac{74747}{39271602}$ & 10 \\
\bottomrule
\end{tabular*}
\end{table}
\subsection{Determining FD Coefficients Based on the optimization:}
Optimization techniques play a crucial role in estimating FD coefficients, effectively mitigating numerical dispersion. Studies have shown that sacrificing the formal order of accuracy through optimized compact schemes can yield surprisingly better wave propagation performance \cite{zhou2011prefactored}, allowing the scheme to operate effectively over a wider band of wavenumbers. $L^2-$norm-based objective functions are commonly employed for assessing differences, owing to their compatibility with least square (LS) methods \cite{kim1996optimized, zhou2011prefactored}. We formulate the misfit function by minimizing the weighted deviation between the scaled true wavenumber $\omega^3$ and the scaled modified wavenumber $\omega'''$ across a specified wavenumber range given in Eqn. (\ref{eqn:7}). Subsequently, we employ the LS approach to solve the objective function, ensuring efficient optimization.
\begin{equation}
\label{eqn:LS}
    \mathrm{E}(a,b,c,\alpha,\beta) = \int_{0}^{r\pi}\bigl(\omega'''-\omega^3\bigr)^2 W(\omega)d\omega,
\end{equation}
where $\mathrm{E}$ is the integral error over the effective wavenumber, $W(\omega)$ is the weighting function, and $r$ is a scalar to control the effective wavenumber range $0 < r \leq 1$ because the value range of $\omega$ is $[0, \pi]$. The choice of the weighting function is crucial to ensure the analytical integrability \cite{kim1996optimized} of the equation (\ref{eqn:LS}) while simultaneously enhancing or diminishing specific wavenumber ranges. The employed weighting function $W(\omega)$ is $[1+2\alpha \cos(\omega)+2\beta \cos(2\omega)]^2$. In this work, we opt for $r=1$ to establish the finite difference coefficients.
\begin{equation}
    \min_{a,b,c,\alpha,\beta}\mathrm{E}= \int_{0}^{r\pi}\biggl\{\frac{2a[8 \sin(\frac{\omega}{2}) - 4 \sin(\omega)]+\frac{2b}{5}[12 \sin(\omega) - 8 \sin(\frac{3\omega}{2})]+\frac{2c}{35}[20 \sin(\omega) - 8 \sin(\frac{5\omega}{2})]}{1+2\alpha \cos(\omega)+2\beta \cos(2\omega)}-\omega^3\biggr\}^2 W(\omega)d\omega.
\end{equation}
\par The LS method solves the optimization problem by setting the derivatives of the misfit function for the unknown parameters to zero and then solving the resulting linear algebraic system:
\begin{itemize}
    \item Fourth-order tri-diagonal TDCCS 
    \begin{equation*}
	   \begin{cases}
			\frac{\partial \mathrm{E}}{\partial a}=0,\\
			1+2 \alpha+ 2 \beta = a+b+c, b=0, c=0, \beta=0. 
	   \end{cases}
    \end{equation*}
    \item Sixth-order tri-diagonal TDCCS 
    \begin{equation*}
	   \begin{cases}
			\frac{\partial \mathrm{E}}{\partial a}=0,\frac{\partial E}{\partial b}=0, \\
			1+2 \alpha+ 2 \beta = a+b+c, c=0, \beta=0. 
	   \end{cases}
    \end{equation*}
    \item Eighth-order tri-diagonal TDCCS 
    \begin{equation*}
	   \begin{cases}
			\frac{\partial \mathrm{E}}{\partial a}=0,\frac{\partial E}{\partial b}=0, \frac{\partial E}{\partial c}=0, \\
			1+2 \alpha+ 2 \beta = a+b+c, \beta=0. 
	   \end{cases}
    \end{equation*}
    \item Tenth-order penta-diagonal TDCCS 
    \begin{equation*}
	   \begin{cases}
			\frac{\partial \mathrm{E}}{\partial a}=0,\frac{\partial E}{\partial b}=0, \frac{\partial E}{\partial c}=0, \\
			1+2 \alpha+ 2 \beta = a+b+c, \\
			 \alpha + 2^2 \beta=\frac{a}{16} +\frac{13b}{80}+\frac{29 c}{80}.
	   \end{cases}
    \end{equation*}
\end{itemize}
Note that we denote TDCCS with the least squares-based coefficients as TDCCS-LS.

\section{Fourier analysis of the errors}
\label{sec:Fourier}
The main incentive behind developing TDCCS schemes is to precisely resolve small scales in multiscale physical problems. Consequently, Fourier analysis is conducted on these optimal schemes to evaluate their spectral characteristics. In this section, we investigate the dispersion and dissipation properties of TDCCS through Fourier analysis. TDCCS, being a central difference, eliminates numerically dissipative errors. The Fourier transformation, a frequently employed tool in finite difference scheme analysis, is applied to equation (2.5), and by utilizing Euler's formula, the modified wavenumber $\omega'''$ of TDCCS can be derived. It is:
\begin{equation}\label{eqn:7}
    \omega'''_{\text{TDCCS}} = \frac{2a[8 \sin(\frac{\omega}{2}) - 4 \sin(\omega)]+\frac{2b}{5}[12 \sin(\omega) - 8 \sin(\frac{3\omega}{2})]+\frac{2c}{35}[20 \sin(\omega) - 8 \sin(\frac{5\omega}{2})]}{1+2\alpha \cos(\omega)+2\beta \cos(2\omega)}.
\end{equation}
Here, $\omega $ represents the scaled wavenumber with $\omega = kh$,  and $\omega'''$ represents the scaled modified wavenumber.
\begin{figure}[htbp!]   
  \begin{minipage}[b]{0.3\linewidth}
    \centering
    \includegraphics[trim=0cm 0cm 0cm 0cm, clip=true,width=\linewidth]{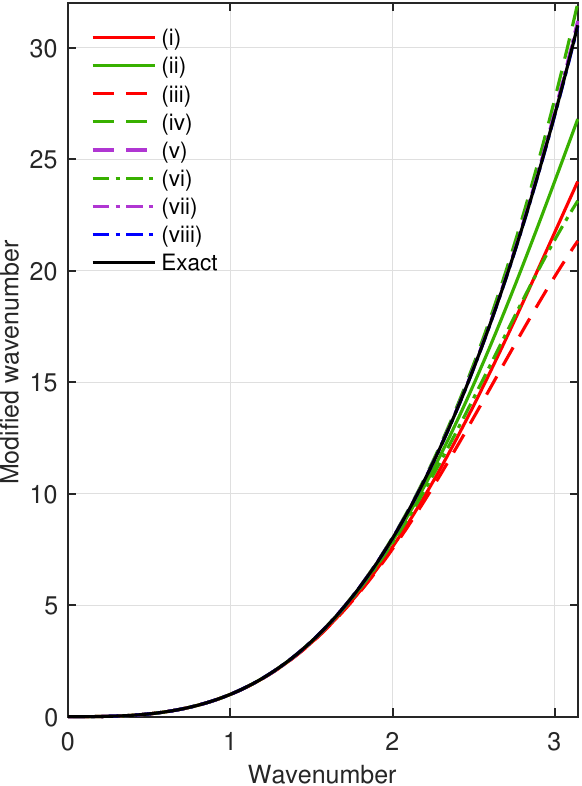}
    \subcaption*{\normalsize{\centering (a) TDCCS-TE}}
  \end{minipage}\hfill
    \begin{minipage}[b]{0.3\linewidth}
    \centering
    \includegraphics[trim=0cm 0cm 0cm 0cm, clip=true,width=\linewidth]{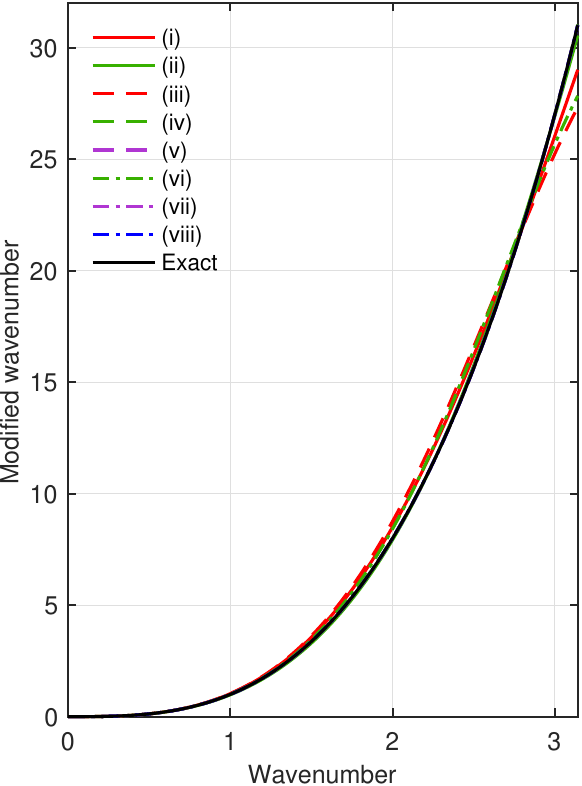}
    \subcaption*{\normalsize{\centering (b) TDCCS-LS}}
  \end{minipage}\hfill
      \begin{minipage}[b]{0.3\linewidth}
    \centering
    \includegraphics[trim=0cm 0cm 0cm 0cm, clip=true,width=\linewidth]{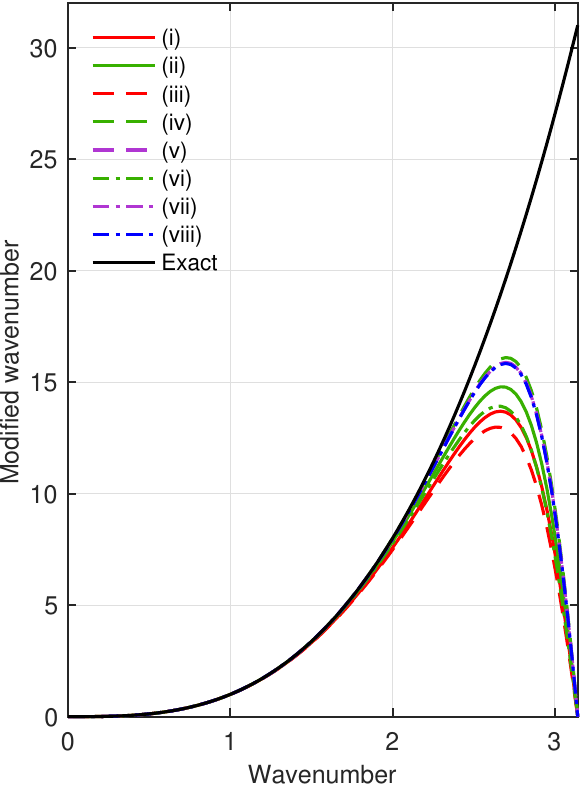}
    \subcaption*{\normalsize{\centering (c) TDCCS-CI}}
  \end{minipage}
  \caption{ Plot of modified wavenumber versus wavenumber for (a) TDCCS-TE, (b) TDCCS-LS and (c) TDCCS-CI : (i) E4; (ii) E6; (iii) T4; (iv) T6; (v) T8; (vi) P6; (vii) P8; (viii) P10.}\label{Fig:F_3}
\end{figure}
\par Figure \ref{Fig:F_3} presents plots comparing modified wavenumber versus wavenumber for third derivative approximations using different methods: (a) the Taylor expansion-based method (TDCCS-TE), (b) the least square optimization-based method (TDCCS-LS), and (c) the tenth-order pentadiagonal compact interpolation (TDCCS-CI). The schemes encompass various orders and types, including (i) fourth-order explicit scheme (E4), (ii) sixth-order explicit scheme (E6), (iii) fourth-order tridiagonal scheme (T4), (iv) sixth-order tridiagonal scheme (T6), (v) eighth-order tridiagonal scheme (T8), (vi) sixth-order pentadiagonal scheme (P6), (vii) eighth-order pentadiagonal scheme (P8), and (viii) tenth-order pentadiagonal scheme (P10). The resolutions achieved by TDCCS-TE and TDCCS-LS surpass those of TDCCCS and TDCNCS. The difference between the modified wavenumber of TDCCS-TE and TDCCS-LS with the exact wavenumber is minimal. In particular, for pentadiagonal schemes, the wave number of TDCCS cannot be distinguished from the exact wave number on the graph. The explicit TDCCS-TE and TDCCS-LS even have higher resolution than the tenth-order pentadiagonal TDCNCS, TDCCCS and TDCCCS-CI. Hence, these schemes exhibit spectral-like resolution. To compare the difference between TDCCS-TE and TDCCS-LS, we define the relative modified wavenumber factor as
\begin{equation*}
    R_{\omega'''}=\frac{\omega'''}{\omega^3} = \frac{2a[8 \sin(\frac{\omega}{2}) - 4 \sin(\omega)]+\frac{2b}{5}[12 \sin(\omega) - 8 \sin(\frac{3\omega}{2})]+\frac{2c}{35}[20 \sin(\omega) - 8 \sin(\frac{5\omega}{2})]}{[1+2\alpha \cos(\omega)+2\beta \cos(2\omega)]\omega^3}.
\end{equation*}

\begin{figure}[htbp!]
  \begin{minipage}[b]{0.45\linewidth}
    \centering
    \includegraphics[trim=0cm 0cm 0cm 0cm, clip=true,width=\linewidth]{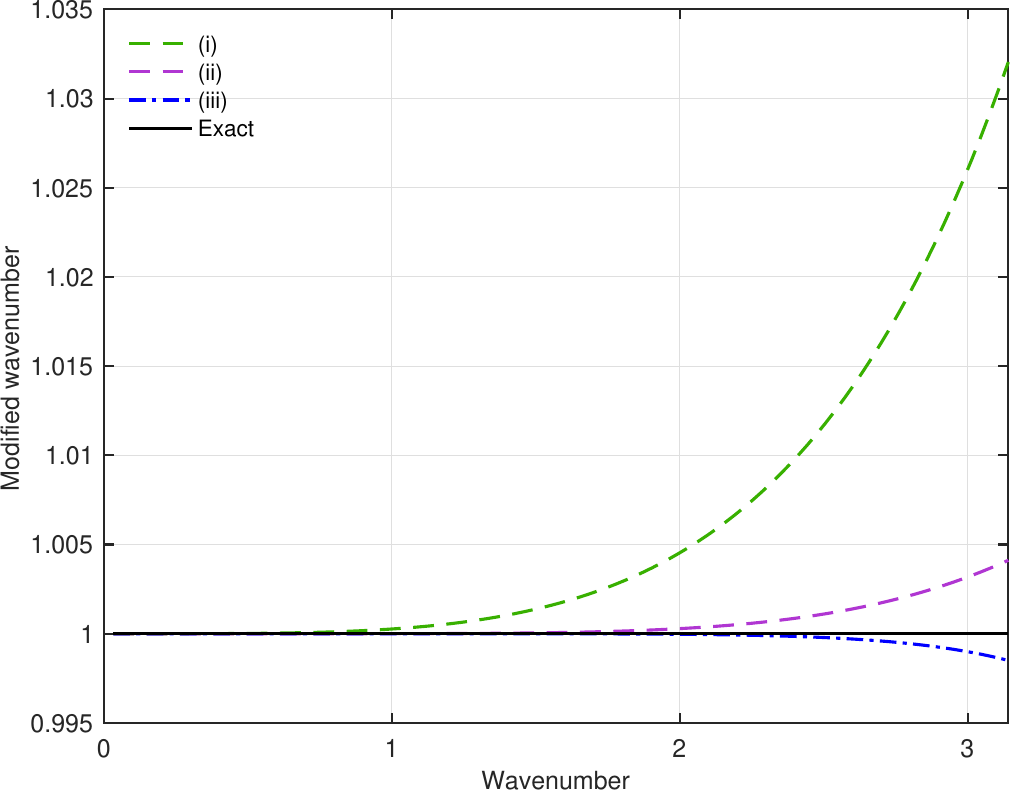}
    \subcaption*{\normalsize{\centering (a) TDCCS-TE}}
  \end{minipage}\hfill
    \begin{minipage}[b]{0.45\linewidth}
    \centering
    \includegraphics[trim=0cm 0cm 0cm 0cm, clip=true,width=\linewidth]{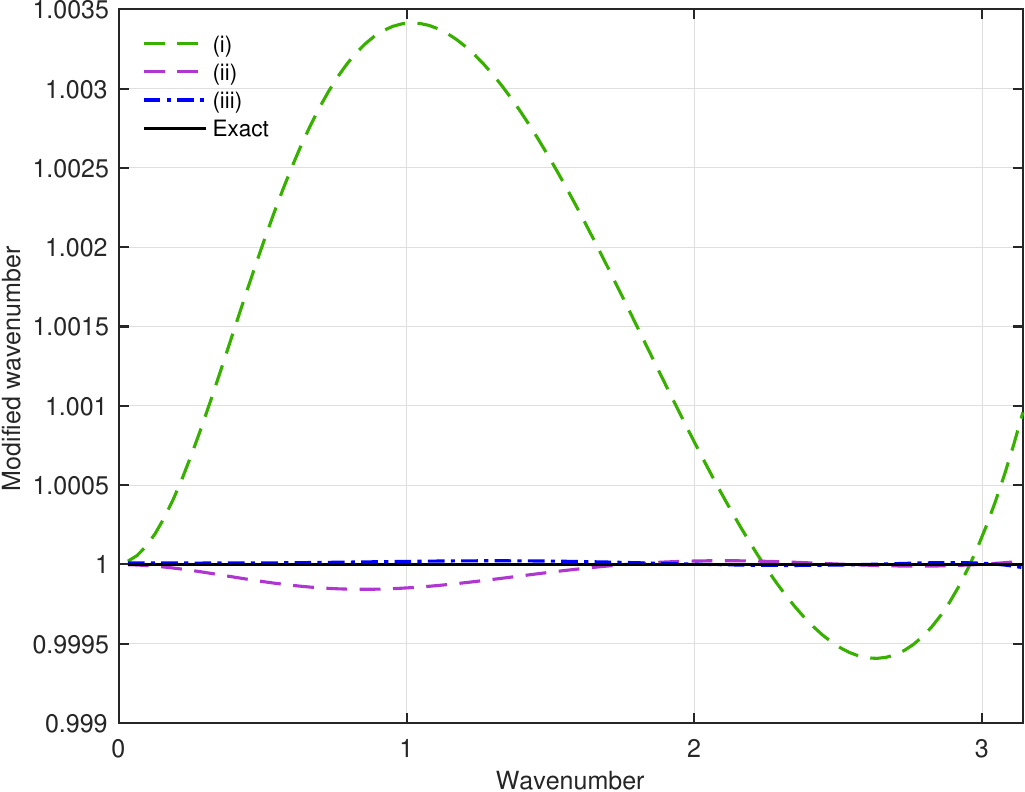}
    \subcaption*{\normalsize{\centering (b) TDCCS-LS}}
  \end{minipage}\hfill
  \caption{ Variations of the relative modified wavenumber factor $R_{\omega'''}$ with wavenumber for TDCCS with coefficients determined by (a) the TE-based method
and (b) the least-square optimization-based method : (i) T6; (ii) T8; (iii) P10.}\label{Fig:F_4}
\end{figure}
 In Figure \ref{Fig:F_4}, the changes in the relative modified wavenumber factor, denoted as $R_{\omega'''}$, versus wavenumbers for TDCCS-TE and TDCCS-LS are shown. In comparison to TDCCS-TE, TDCCS-LS notably expands the range of wavenumbers while maintaining the same order of truncation error. The eighth-order tridiagonal TDCCS-LS even exhibits comparable accuracy to the tenth-order pentadiagonal TDCCS-TE.
 \par
 The bandwidth resolving efficiency~\cite{Lele} is a quantitative indicator of spectral resolution. The resolving efficiency of the FD scheme, with a specified error tolerance, is defined as 
 \begin{equation*}
    e = \frac{\omega_f}{\pi},
\end{equation*}
where $\omega_f$ is the shortest well-resolved wave component satisfying 
\begin{equation*}
   \biggl|\frac{\omega'''(\omega)-\omega^3}{\omega^3}\biggr| \leq \epsilon_t,
\end{equation*}
and $\epsilon_t$ represents the error tolerance threshold. This threshold remains constant when comparing various compact FD schemes. Tables (\ref{Table:T_2a}) and (\ref{Table:T_2b}) present the values of the bandwidth resolving efficiency ($e$) for various compact finite difference (FD) schemes, with error tolerances set at $\epsilon_t$ = 0.001 and $\epsilon_t$ = 0.0001, respectively. Tables (\ref{Table:T_2a}) and (\ref{Table:T_2b}) highlight that TDCCS exhibits the highest resolving efficiency, effectively capturing a broader range of wavenumbers under the same error tolerance. Specifically, for $\epsilon_t$ = 0.001, the resolving efficiency of an eighth-order TDCNCS is 0.5018, while that of TDCCS-TE and TDCCS-LS is 0.7828 and 0.9998, respectively. Similarly, for $\epsilon_t$ = 0.0001, the resolving efficiency of an eighth-order TDCNCS is 0.3855, whereas that of TDCCS-TE and TDCCS-LS is 0.5376 and 0.9998, respectively.  Additionally, when compared to other potential combinations such as TDCCS-TE-1, TDCCS-TE-2, and TDCCS-TE-3, TDCCS-TE exhibits favorable resolving efficiency. Therefore, we adopt the TDCCS-TE scheme in this paper.
\begin{table}[!ht]
\centering
\captionof{table}{The shortest well-resolved wave $\omega_f$ and resolving efficiency $e$ of different schemes with tolerance error
$\epsilon_t$ = 0.001.}\label{Table:T_2a}
\setlength{\tabcolsep}{0pt}
\begin{tabular*}{\textwidth}{@{\extracolsep{\fill}} l *{8}{c} }
\toprule
\textbf{Schemes/order} &
\multicolumn{2}{c}{\textbf{4-T}} &
\multicolumn{2}{c}{\textbf{6-T}} &
\multicolumn{2}{c}{\textbf{8-T}} &
\multicolumn{2}{c}{\textbf{10-P}}\\
\cmidrule{2-3} \cmidrule{4-5} \cmidrule{6-7} \cmidrule{8-9}
& $\omega_f$ & $e$
& $\omega_f$ & $e$
& $\omega_f$ & $e$
& $\omega_f$ & $e$\\
\midrule
TDCNCS & 0.693 & 0.2205 & 1.735 & 0.5523 & 1.576 & 0.5018 & 1.635 & 0.5205 \\ 
TDCCCS & 0.716 & 0.2278 & 1.164 & 0.3705 & 1.468 & 0.4672 & 1.845 & 0.5874 \\
TDCCCS-CI & 0.715 & 0.2277 & 1.162 & 0.3699 & 1.445 & 0.4600 & 1.682 & 0.5354 \\
\midrule
TDCCS-CI & 0.722 & 0.2297 & 2.037 & 0.6483 & 1.769 & 0.5631 & 1.748 & 0.5565 \\
TDCCS-TE & 0.722 & 0.2297 & 1.386 & 0.4411 & 2.459 & 0.7828 & 2.998 & 0.9542 \\
TDCCS-LS & 2.795 & 0.8898 & 3.141 & 0.9998 & 3.141 & 0.9998 & 3.141 & 0.9998 \\
\midrule
TDCCS-CI-1 & 0.722 & 0.2297 & 2.037 & 0.6483 & 2.037 & 0.6483 & 2.037 & 0.6483 \\
TDCCS-TE-1 & 0.722 & 0.2297 & 1.386 & 0.4411 & 2.412 & 0.7679 & 2.928 & 0.9321 \\ 
TDCCS-LS-1 & 2.795 & 0.8898 & 3.141 & 0.9998 & 3.141 & 0.9998 & 3.141 & 0.9998 \\ 
\midrule
TDCCS-TE-2 & 0.722 & 0.2297 & 1.244 & 0.3959 & 1.612 & 0.5131 & 2.044 & 0.6506 \\
\midrule
TDCCS-TE-3 & 0.722 & 0.2297 & 1.244 & 0.3959 & 1.618 & 0.5152 & 1.919 & 0.6110 \\
\bottomrule
\end{tabular*}
\end{table}
\begin{table}[!ht]
\centering
\captionof{table}{The shortest well-resolved wave $\omega_f$ and resolving efficiency $e$ of different schemes with tolerance error $\epsilon_t$ = 0.0001.}\label{Table:T_2b}
\setlength{\tabcolsep}{0pt}
\begin{tabular*}{\textwidth}{@{\extracolsep{\fill}} l *{8}{c} }
\toprule
\textbf{Schemes/order} &
\multicolumn{2}{c}{\textbf{4-T}} &
\multicolumn{2}{c}{\textbf{6-T}} &
\multicolumn{2}{c}{\textbf{8-T}} &
\multicolumn{2}{c}{\textbf{10-P}}\\
\cmidrule{2-3} \cmidrule{4-5} \cmidrule{6-7} \cmidrule{8-9}
& $\omega_f$ & $e$
& $\omega_f$ & $e$
& $\omega_f$ & $e$
& $\omega_f$ & $e$\\
\midrule
TDCNCS & 0.392 & 0.1248 & 1.524 & 0.4850 & 1.211 & 0.3855 & 1.293 & 0.4114 \\
TDCCCS & 0.405 & 0.1290 & 0.800 & 0.2545 & 1.105 & 0.3518 & 1.493 & 0.4753 \\
TDCCCS-CI & 0.405 & 0.1290 & 0.799 & 0.2544 & 1.097 & 0.3490 & 1.362 & 0.4336 \\
\midrule
TDCCS-CI & 0.407 & 0.1294 & 1.994 & 0.6347 & 1.468 & 0.4672 & 1.420 & 0.4521 \\
TDCCS-TE & 0.407 & 0.1294 & 0.785 & 0.2497 & 1.689 & 0.5376 & 2.288 & 0.7284 \\
TDCCS-LS & 2.792 & 0.8888 & 2.982 & 0.9493 & 3.141 & 0.9998 & 3.044 & 0.9690\\
\midrule
TDCCS-CI-1 & 0.407 & 0.1294 & 1.994 & 0.6347 & 1.994 & 0.6347 & 1.994 & 0.6347 \\
TDCCS-TE-1 & 0.407 & 0.1294 & 0.785 & 0.2497 & 1.659 & 0.5280 & 2.228 & 0.7093 \\ 
TDCCS-LS-1 & 2.792 & 0.8888 & 2.982 & 0.9493 & 3.141 & 0.9998 & 3.141 & 0.9998 \\ 
\midrule
TDCCS-TE-2 & 0.407 & 0.1294 & 0.854 & 0.2718 & 1.231 & 0.3917 & 1.585 & 0.5046 \\
\midrule
TDCCS-TE-3 & 0.407 & 0.1294 & 0.854 & 0.2718 & 1.234 & 0.3928 & 1.515 & 0.4823 \\
\bottomrule
\end{tabular*}
\end{table}


\section{Low-pass spatial filter}
\label{sec:Filter}
With high-order finite difference schemes, it is necessary to artificially mitigate all spurious waves while preserving the accuracy and resolution of the computed solution. This can be performed by regularizing the computed solution via a high-order low-pass spatial filter (LPSF)~\cite{Lele}.  A tridiagonal filter of high order can be expressed as
\begin{equation}
\label{eqn:Filter}
\alpha_F \hat{f}_{j-1} + \hat{f}_{j} + \alpha_F \hat{f}_{j+1} = \sum_{n=0}^{N} \frac{a_n}{2} (f_{j+n} + f_{j-n})
\end{equation}
where $f_j$ denotes the given value at point $j$, and $\hat{f}_j$ denotes the value after filtering. The most natural formulation of the problem involves the transfer function associated with the equation (\ref{eqn:Filter}).
\begin{equation}
T(\omega) = \frac{\sum_{n=0}^{N} a_n \cos(n\omega)}{1+2\alpha_F \cos(\omega)} 
\end{equation}
To determine the unknown coefficients, we require the exclusion of the highest-frequency mode by imposing the condition $T(\pi) = 0$. For adaptability, we maintain $\alpha_F$ as an unrestricted parameter. Subsequently, the remaining $N$ equations can be deduced by equating the Taylor series coefficients of the left and right sides. Through this process, equation (\ref{eqn:Filter}) yields a $2N$-th order formula within a $2N +1$ point stencil. It is important to note that $T(\omega)$ is real, signifying that the filter alters only the amplitude of each wave component without impacting the phase. Here $\alpha_F$ is a free parameter that satisfies $|\alpha_F|< 0.5$. The coefficients of a family of 8th-order (F8), 10th-order (F10), and 12th-order filter (F12) are given in the Table \ref{Table:Filter}.
\begin{table}[!ht]
\centering
\captionof{table}{Coefficients for the filter formula.}
\setlength{\tabcolsep}{0pt}
\begin{tabular*}{\textwidth}{@{\extracolsep{\fill}} l *{9}{c} }
\toprule
\textbf{Scheme} &
\textbf{$a_0$} &
\textbf{$a_1$} &
\textbf{$a_2$} &
\textbf{$a_3$} &
\textbf{$a_4$} &
\textbf{$a_5$} &
\textbf{$a_6$} &
\textbf{Order}\\
\midrule
F-8 & $\frac{(93+70 \alpha_F)}{128}$ & $\frac{(7+18 \alpha_F)}{16}$ & $\frac{(-7+14 \alpha_F)}{32}$ & $\frac{(1-2 \alpha_F)}{16}$ & $\frac{(-1+2 \alpha_F)}{128}$ & 0 & 0 & 8\\ [0.20cm]
F-10 & $\frac{(193+126 \alpha_F)}{256}$ & $\frac{(105+302 \alpha_F)}{256}$ & $\frac{15(-1+2 \alpha_F)}{64}$ & $\frac{45(1-2 \alpha_F)}{512}$ & $\frac{5(-1+2 \alpha_F)}{256}$ & $\frac{(1-2 \alpha_F)}{512}$ & 0 & 10\\ [0.20cm]
F-12 & $\frac{(793 + 462 \alpha_F) }{1024}$ & $\frac{ (99 + 314 \alpha_F) }{256}$ & $\frac{495 (-1 +2 \alpha_F)}{2048}$ & $\frac{55 (1 -2 \alpha_F)}{512}$ & $\frac{33 (-1 +2 \alpha_F)}{1024}$ & $\frac{3 (1 -2 \alpha_F)}{512}$ & $\frac{ (-1 +2 \alpha_F)}{2048}$ & 12\\ [0.15cm]
\bottomrule
\end{tabular*}\label{Table:Filter}
\end{table}


\section{Time Integration}
\label{sec:time_integration}
\begin{figure}[htbp!]
  \begin{minipage}[b]{0.45\linewidth}
    \centering
    \includegraphics[trim=0cm 0cm 0cm 0cm, clip=true,width=\linewidth]{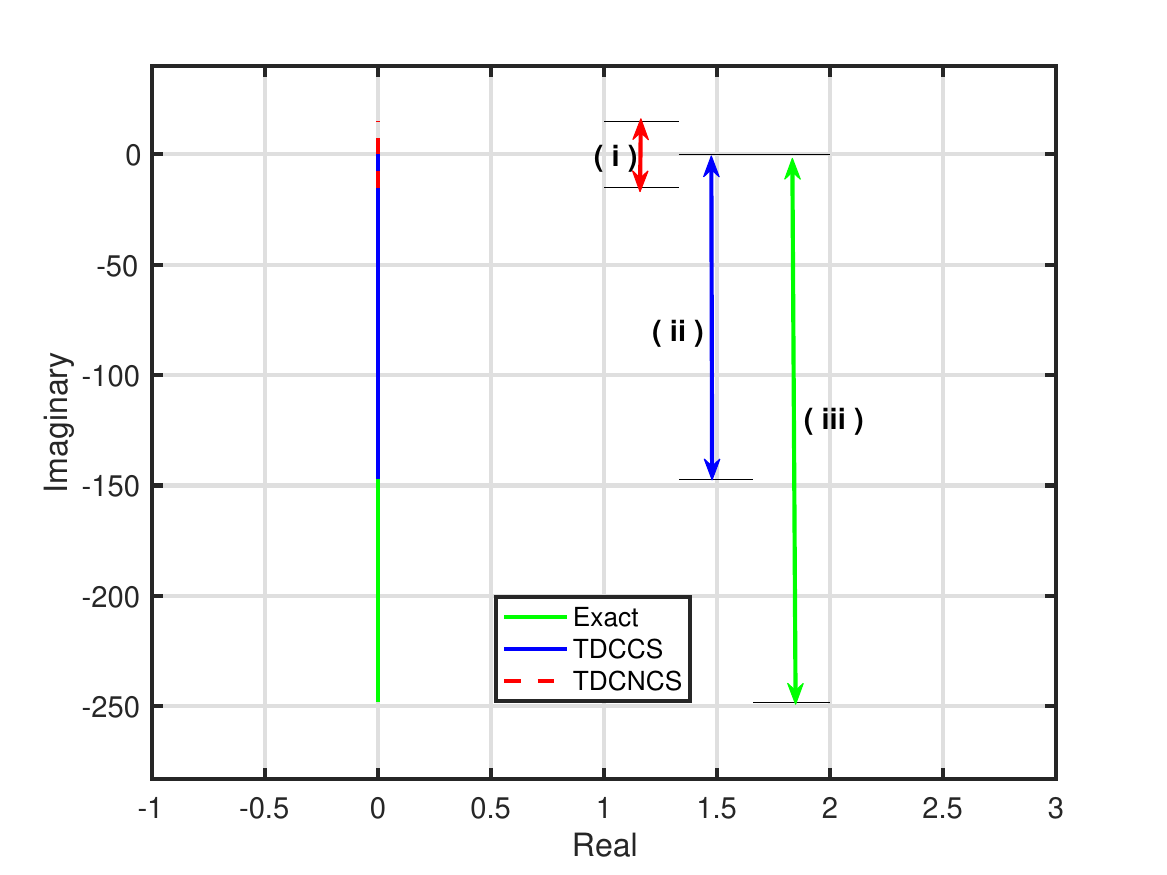}
    \subcaption*{\normalsize{\centering (a)}}
  \end{minipage}\hfill
    \begin{minipage}[b]{0.45\linewidth}
    \centering
    \includegraphics[trim=0cm 0cm 0cm 0cm, clip=true,width=\linewidth]{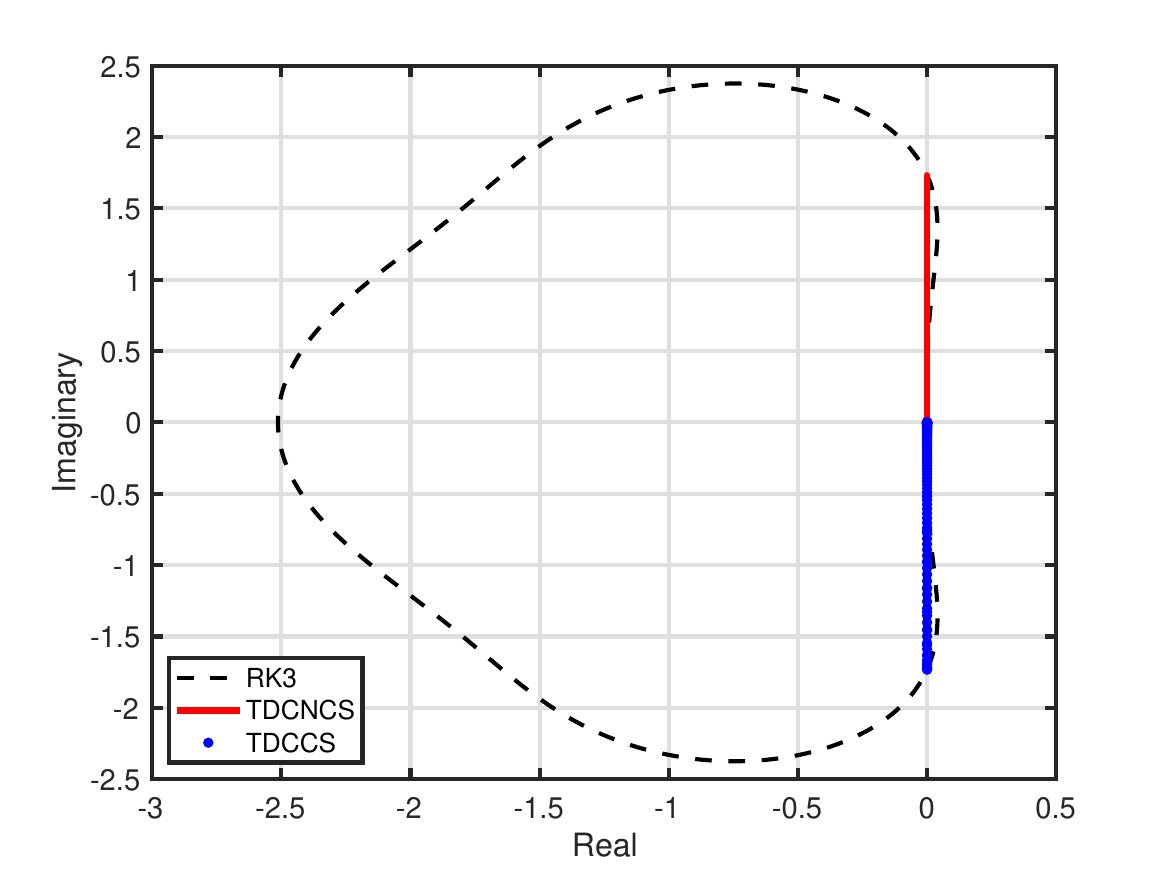}
    \subcaption*{\normalsize{\centering (b)}}
  \end{minipage}\hfill
  \caption{ (a) Eigenvalues of the (i) Exact, (ii) TDCCS and (iii) TDCNCS. (b) Stability regions of three-stage Runge–Kutta time-integration scheme and eigenvalues of the TDCNCS and TDCCS scaled by their maximum stable CFL.}\label{Fig:Stabilty}
\end{figure}
Discretizing the spatial derivatives using the compact scheme described in equations~(\ref{eqn:5})--(\ref{eqn:6}) yields an ordinary differential equation (ODE), given as
\begin{equation}
    \frac{du}{dt} =S(u),
\end{equation}
where $S(u)$ is the spatially-discretized approximation of the right-hand-side. The ODE is evolved in time using the third-order TVD Runge-Kutta (TVDRK3) method~\cite{gottlieb1998total}:
\begin{equation}
\label{eqn:TVD}
\begin{split}
u^{(1)} &= u^n + \Delta t S(u^n),\\
u^{(2)} &= \frac{3}{4} u^n + \frac{1}{4} u^{(1)} + \frac{1}{4} \Delta t S(u^{(1)}),\\
u^{n+1} &= \frac{1}{3} u^n + \frac{2}{3} u^{(2)} + \frac{2}{3} \Delta t S(u^{(2)}).\\
\end{split}
\end{equation}
Figure~\ref{Fig:Stabilty}(a) shows the eigenvalues of the exact, TDCNCS, and TDCCS schemes, while Fig.~\ref{Fig:Stabilty}(b) shows the stability region boundary of the TVDRK3 method. Linear stability requires a time step $\Delta t$ such that the eigenvalues of the spatial discretization scheme scaled by $\Delta t/\Delta x^3$ fall within the stability region of the time integration method. Since the eigenvalues of the TDCNCS and TDCCS scheme are purely imaginary, we consider the value at which the stability region boundary of TVDRK3 intercepts the imaginary axis, i.e., $\pm 1.732$. The maximum magnitude of the eigenvalues for the TDCNCS and TDCCS schemes are $15.157$ and $147.168$, respectively. This results in the following linear stability bound on the time step:
\begin{align}
    \frac{\Delta t}{\Delta x^3} \leq 0.11\ \left({\rm TDCNCS}\right),\ \ 0.012\ \left({\rm TDCCS}\right).
\end{align}
Figure~\ref{Fig:Stabilty}(b) shows the eigenvalues of the TDCNCS and TDCCS schemes scaled by these limits, and they lie within the stability region of TVDRK3. The TDCCS scheme has a very restrictive time step bound compared with the TDCNCS scheme, and this is a significant drawback. In future work, we propose to develop optimized high-order explicit Runge-Kutta methods for the TDCCS schemes~\cite{kubatkoetal2014} to allow larger time steps.

\section{Numerical results}
\label{sec:num_ex}
This section presents numerous numerical examples to showcase the accuracy and effectiveness of the method. The numerical illustrations aim to emphasize the high-order accuracy of the approach in addressing one-dimensional problems, both linear and non-linear. Additionally, the method is evaluated for its performance in handling convection-dominated cases, especially instances where the coefficients of the third derivative terms are small. For numerical experiments involving both convection and dispersion terms, spatial discretization is performed using two distinct numerical schemes: the eighth-order cell-node compact scheme (CNCS)~\cite{Lele} (first derivative) with the eighth-order TDCNCS~(\ref{eqn:2}) (third derivative), and the eighth-order central compact scheme (CCS)~\cite{LZZS} (first derivative) with the eighth-order TDCCS~(\ref{eqn:5}) (third derivative).We employed eighth-order TDCCS-TE instead of TDCCS-LS in the spatial discretization, as TDCCS-LS fails to preserve order of accuracy. The following error norms are used to compute the accuracy of schemes.
  \begin{equation*}
L^{\infty} = \max_{0\leq i \leq N} |u_e - u_{a}|,\,\,\,\,\,\,\,\,\,\, L^{1} = \dfrac{1}{N+1} \sum_{i=0}^{N} |u_e - u_{a}|, \,\,\,\,\,\,\,\,\,\, L^{2} = \biggl(\dfrac{1}{N+1} \sum_{i=0}^{N} |u_e - u_{a}|^2\biggr)^{1/2},
 \end{equation*}
where $u_e$ and $u_{a}$ denote the exact and approximate solutions of the PDE. To ensure a stable numerical solution for linear cases, a CFL value of 0.01 has been selected in the calculation of time step size $\Delta t$, given by $\text{CFL} \times \Delta x^3$. This choice falls within the stable region for both the TDCNCS and TDCCS schemes, ensuring dependable and physically meaningful results.

\begin{example}\label{example:1}
\normalfont We solve the following linear one-dimensional KdV equation
\begin{equation}
\left\{
\begin{aligned}
u_t+c^{-2}u_{xxx} &= 0, \quad (x,t) \in[0,2\pi] \times [0,1], \\
u_0(x)&= \sin(cx), \quad x \in[0,2\pi],
\end{aligned}
\right.  
 \end{equation}
 with periodic boundary conditions for $c=1,8$. The exact solution is represented by a left-moving wave given as $u(x,t)=\sin(c(x+t))$ for $x \in [0, 2\pi]$ and $t \in [0,1]$. Initially, we solve~(\ref{example:1}) with $c=8$, employing eighth-order TDCNCS and TDCCS schemes with a step size of $\Delta t = 0.01 \Delta x^3$. Figure~\ref{Figure:E1a} shows the resulting solutions and pointwise errors for $N=40$ at times $t=0,0.3, 0.7$, and $1$. Additionally, Table (\ref{Table:E1a}) provides a summary of the $L^{\infty}$-error, $L^{1}$-error, $L^{2}$-error, and the order of convergence for both considered schemes at the final time $t = 1$.  Notably, the absolute error is an order of magnitude lower for the TDCCS scheme compared to the TDCNCS. Specifically, for $N=40$, the $L^{\infty}$-error for TDCNCS is $1.0796\times 10^{-3}$, whereas for TDCCS, it is $1.1749\times 10^{-4}$. For a larger value of $N=160$, the $L^{\infty}$-error for TDCNCS is $1.2708\times 10^{-8}$, whereas for TDCCS, it is $5.1353\times 10^{-9}$. However, it is worth noting that TDCNCS maintains eighth-order accuracy, while TDCCS achieves close to eighth-order for a large value of $N$.
 \begin{figure}[htbp!]
  \begin{minipage}[b]{0.45\linewidth}
    \centering
    \includegraphics[trim=0cm 0cm 0cm 0cm, clip=true,width=\linewidth]{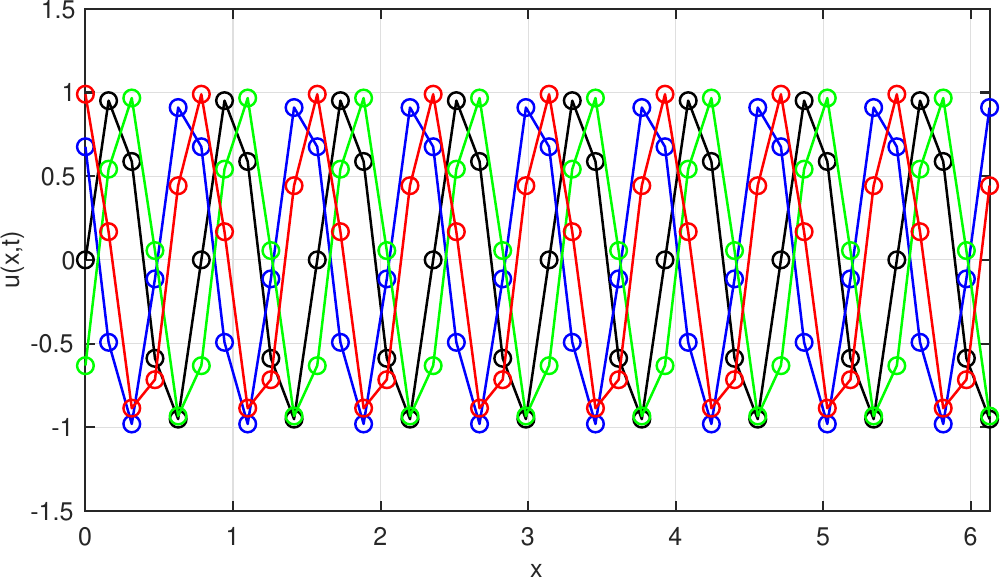}
  \end{minipage}\hfill
    \begin{minipage}[b]{0.45\linewidth}
    \centering
    \includegraphics[trim=0cm 0cm 0cm 0cm, clip=true,width=\linewidth]{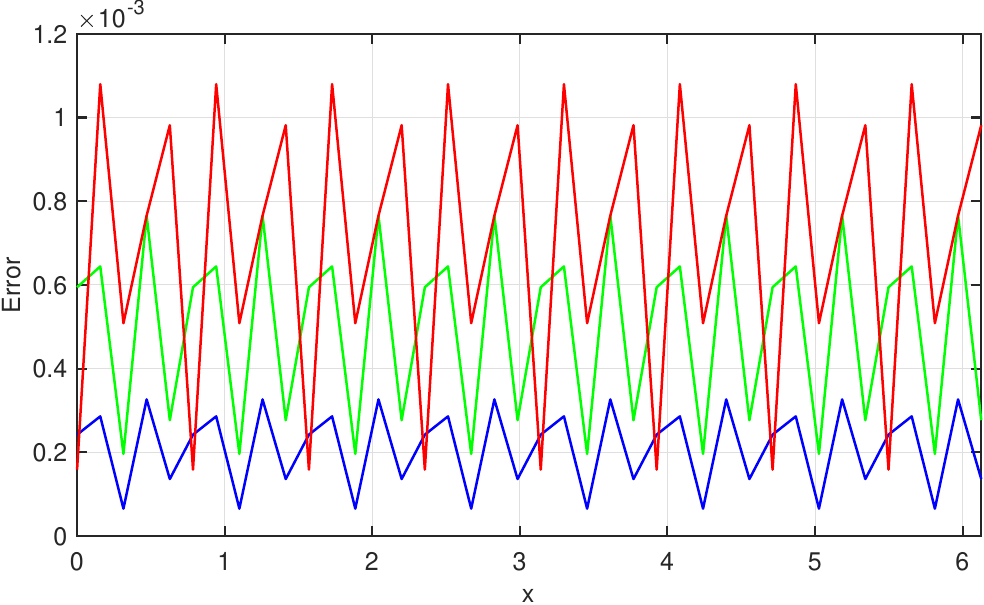}
  \end{minipage}\hfill
  \begin{minipage}[b]{0.45\linewidth}
    \centering
    \includegraphics[trim=0cm 0cm 0cm 0cm, clip=true,width=\linewidth]{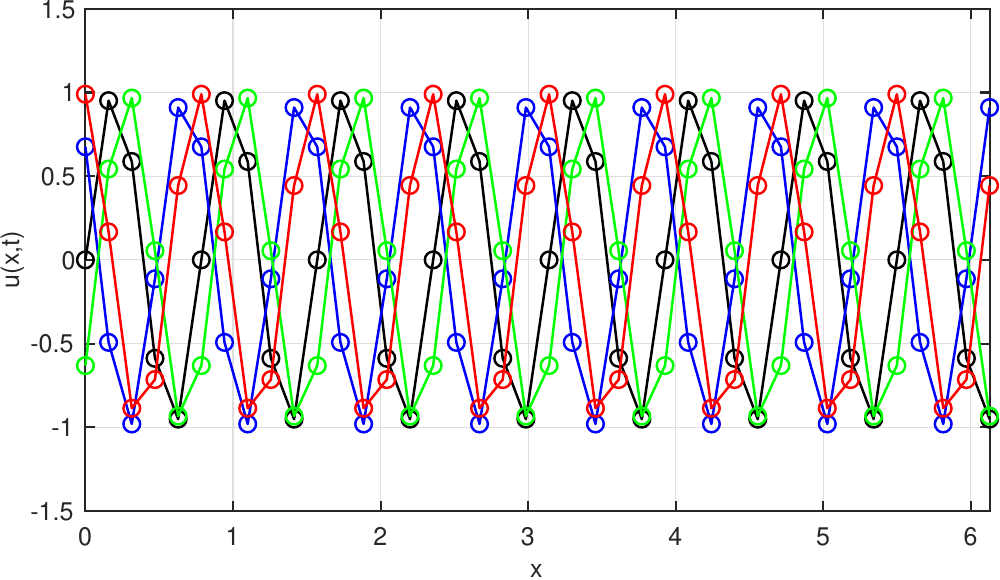}
  \end{minipage}\hfill
  \begin{minipage}[b]{0.45\linewidth}
    \centering
    \includegraphics[trim=0cm 0cm 0cm 0cm, clip=true,width=\linewidth]{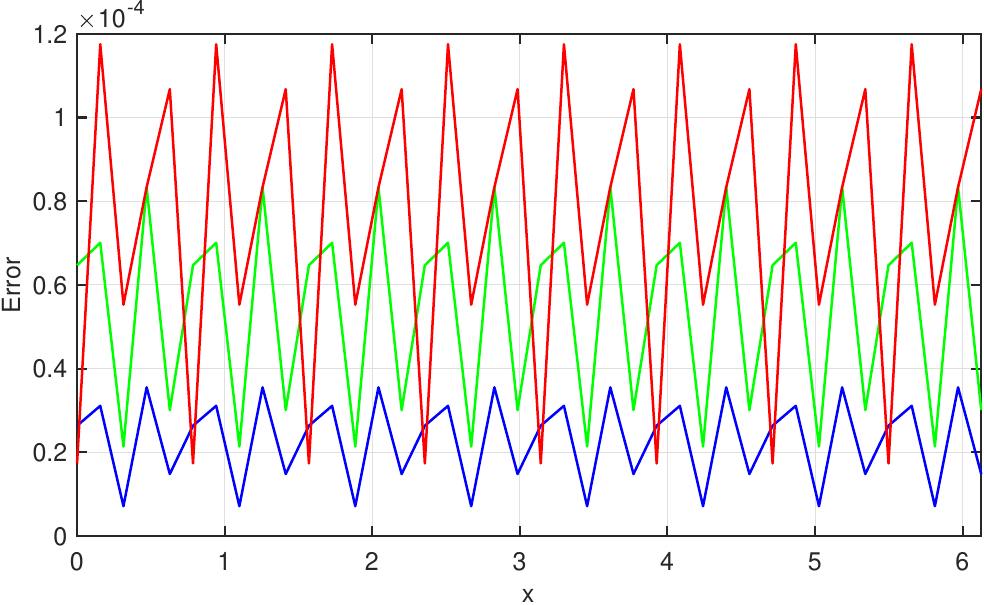}
  \end{minipage}\hfill
  \caption{Solutions and errors for Example \ref{example:1} with c=8. The first and second rows are obtained by TDCNCS-T8 and TDCCS-T8 respectively. Left column shows the exact (solid line) and calculated  ($\bigcirc$) left moving wave at $t =$ 0 (black), $t =$ 0.3 (blue), $t =$ 0.7 (green), $t =$ 1 (red). Right column shows the corresponding errors increasing with time.}\label{Figure:E1a}
\end{figure}
\begin{figure}[htbp!]
  \begin{minipage}[b]{0.45\linewidth}
    \centering
    \includegraphics[trim=0cm 0cm 0cm 0cm, clip=true,width=\linewidth]{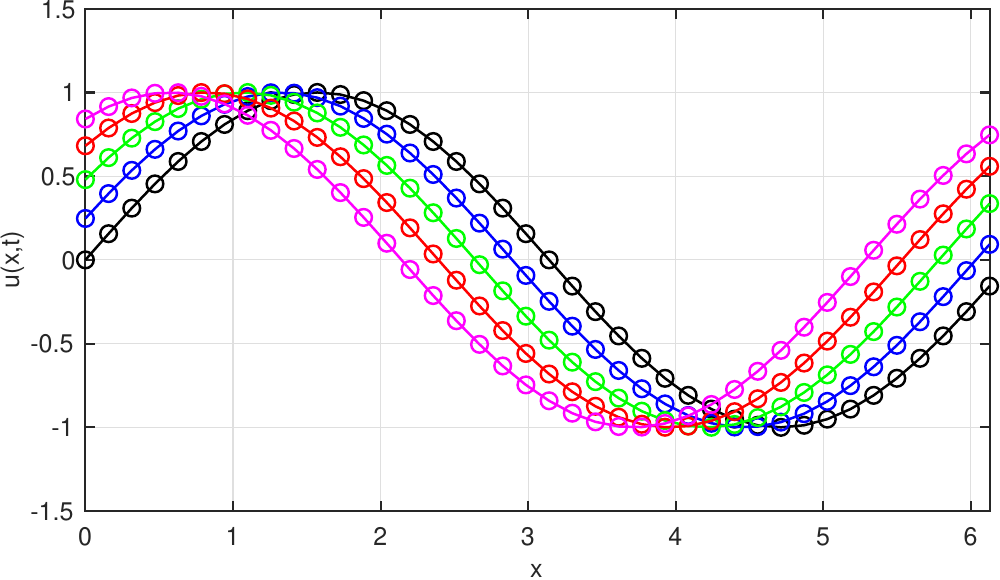}
  \end{minipage}\hfill
    \begin{minipage}[b]{0.45\linewidth}
    \centering
    \includegraphics[trim=0cm 0cm 0cm 0cm, clip=true,width=\linewidth]{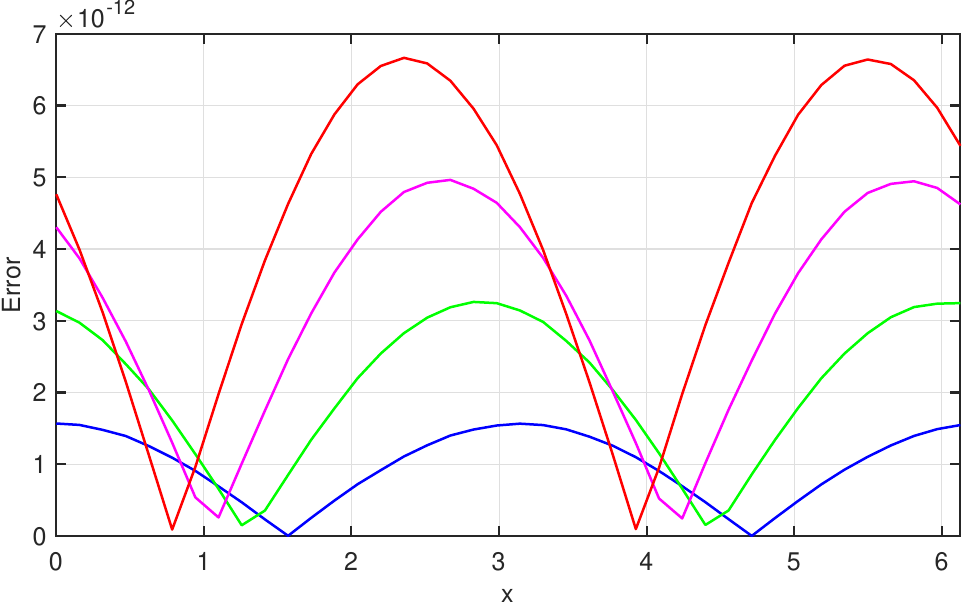}
  \end{minipage}\hfill
  \begin{minipage}[b]{0.45\linewidth}
    \centering
    \includegraphics[trim=0cm 0cm 0cm 0cm, clip=true,width=\linewidth]{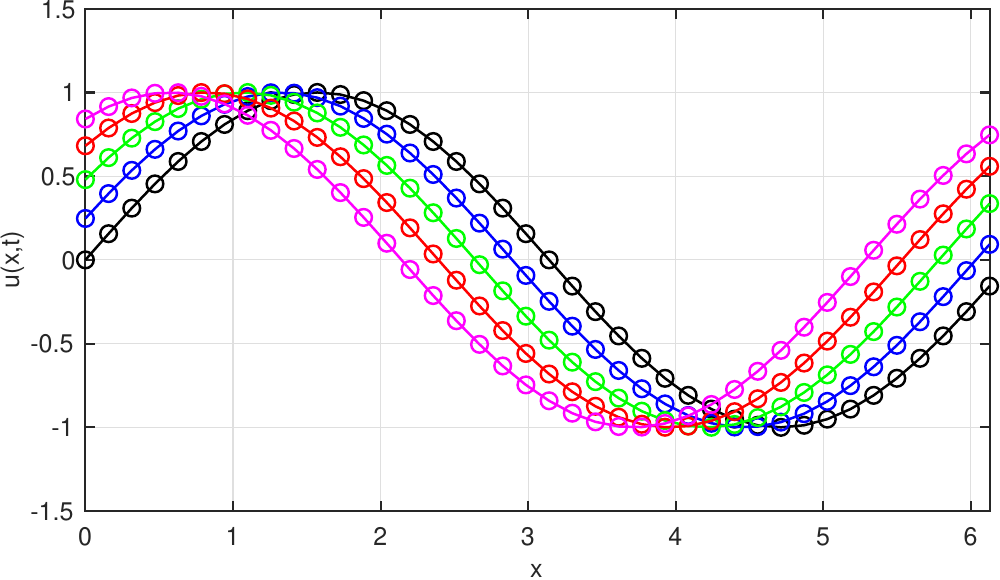}
  \end{minipage}\hfill
  \begin{minipage}[b]{0.45\linewidth}
    \centering
    \includegraphics[trim=0cm 0cm 0cm 0cm, clip=true,width=\linewidth]{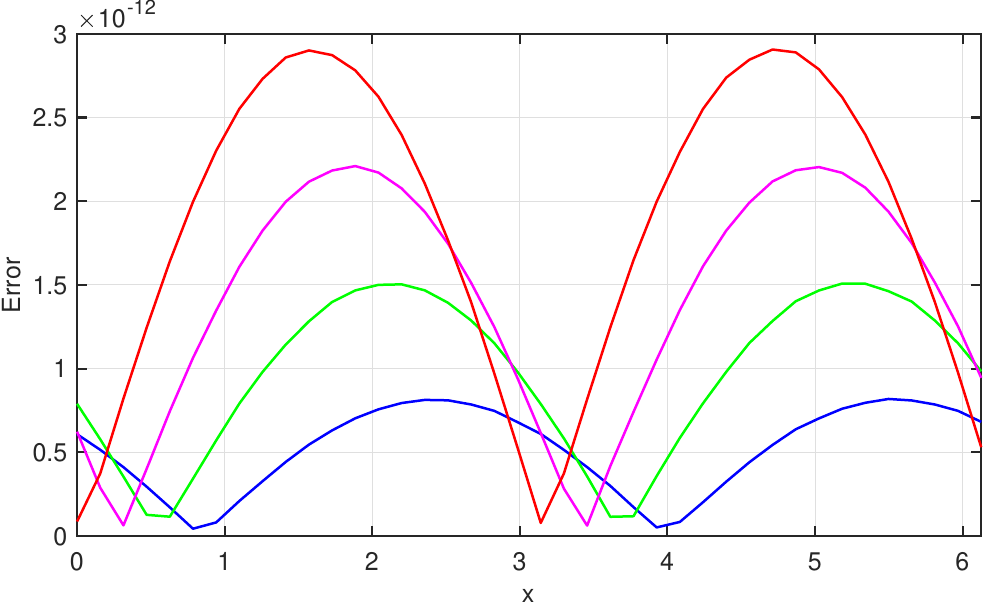}
  \end{minipage}\hfill
  \caption{ Solutions and errors for Example \ref{example:1} with c=1. The first and second rows are obtained by TDCNCS-T8 and TDCCS-T8 respectively. Left column shows the exact (solid line) and calculated ($\bigcirc$) left moving wave at $t =$ 0 (black), $t =$ 0.25 (blue), $t =$ 0.5 (green), $t =$ 0.75 (magenta), $t =$ 1 (red). The right column shows the corresponding errors increasing with time.}\label{Figure:E1b}
\end{figure}
\begin{table}[htbp!]
\centering
\captionof{table}{Errors and orders of convergence with eight-order methods for Example \ref{example:1} with $c=8$ at $t = 1$.}
\begin{tabular*}{\textwidth}{@{\extracolsep{\fill}} l *{8}{c} }
\midrule
\textbf{Scheme} & \textbf{N} & $\boldsymbol{L^{\infty}}$\textbf{-error} & \textbf{Rate} & $\boldsymbol{L^{1}}$\textbf{-error} & \textbf{Rate} & $\boldsymbol{L^{2}}$\textbf{-error} & \textbf{Rate}\\
\toprule
\textbf{TDCNCS} & 20 & 7.9125e-01 & - & 5.1211e-01 & - & 5.4678e-01 & - \\ 
\textbf{(c=8)} & 40 & 1.0796e-03 & 9.5175 & 6.9871e-04 & 9.5175 & 7.6486e-04 & 9.4816 \\ 
& 60 & 3.6487e-05 & 8.3543 & 2.3271e-05 & 8.3904 & 2.5610e-05 & 8.3773 \\ 
& 80 & 3.4195e-06 & 8.2294 & 2.2132e-06 & 8.1784 & 2.4374e-06 & 8.1759 \\ 
& 100 & 5.6767e-07 & 8.0473 & 3.6163e-07 & 8.1184 & 3.9977e-07 & 8.1014 \\ 
& 120 & 1.3038e-07 & 8.0686 & 8.3154e-08 & 8.0623 & 9.1887e-08 & 8.0645 \\ 
& 140 & 3.7691e-08 & 8.0508 & 2.4003e-08 & 8.0604 & 2.6558e-08 & 8.0520 \\ 
& 160 & 1.2708e-08 & 8.1418 & 8.2008e-09 & 8.0427 & 9.0447e-09 & 8.0667 \\ 
\midrule
\textbf{TDCCS}  & 20 & 8.9768e-03 & - & 5.8099e-03 & - & 6.2784e-03 & - \\ 
\textbf{(c=8)} & 40 & 1.1749e-04 & 6.2556 & 7.6040e-05 & 6.2556 & 8.3230e-05 & 6.2371 \\ 
& 60 & 7.5798e-06 & 6.7598 & 4.8343e-06 & 6.7960 & 5.3201e-06 & 6.7826 \\ 
& 80 & 9.5509e-07 & 7.2004 & 6.1815e-07 & 7.1494 & 6.8077e-07 & 7.1469 \\ 
& 100 & 1.8581e-07 & 7.3364 & 1.1837e-07 & 7.4074 & 1.3085e-07 & 7.3906 \\ 
& 120 & 4.6838e-08 & 7.5582 & 2.9873e-08 & 7.5518 & 3.3012e-08 & 7.5536 \\ 
& 140 & 1.4426e-08 & 7.6397 & 9.1866e-09 & 7.6497 & 1.0165e-08 & 7.6411 \\ 
& 160 & 5.1353e-09 & 7.7352 & 3.3195e-09 & 7.6232 & 3.6606e-09 & 7.6489 \\ 
\bottomrule
\end{tabular*}\label{Table:E1a}
\end{table} 
\begin{table}[htbp!]
\centering
\captionof{table}{Errors and orders of convergence with eight-order methods for Example \ref{example:1} with $c=1$ at $t = 1$.}
\begin{tabular*}{\textwidth}{@{\extracolsep{\fill}} l *{8}{c} }
\midrule
\textbf{Scheme} & \textbf{N} & $\boldsymbol{L^{\infty}}$\textbf{-error} & \textbf{Rate} & $\boldsymbol{L^{1}}$\textbf{-error} & \textbf{Rate} & $\boldsymbol{L^{2}}$\textbf{-error} & \textbf{Rate}\\
\toprule
\textbf{TDCNCS} & 10 & 4.1920e-07 & - & 2.7131e-07 & - & 2.9230e-07 & - \\ 
\textbf{(c=1)} & 20 & 1.6089e-09 & 8.0254 & 1.0249e-09 & 8.0483 & 1.1120e-09 & 8.0382 \\ 
& 30 & 6.2433e-11 & 8.0135 & 3.9821e-11 & 8.0104 & 4.3504e-11 & 7.9935 \\ 
& 40 & 6.6641e-12 & 7.7772 & 4.2295e-12 & 7.7944 & 4.6437e-12 & 7.7771 \\ 
\midrule
\textbf{TDCCS}  & 10 & 1.1729e-07 & - & 7.5910e-08 & - & 8.1641e-08 & - \\ 
\textbf{(c=1)} & 20 & 6.4028e-10 & 7.5171 & 4.0809e-10 & 7.5393 & 4.4262e-10 & 7.5271 \\ 
& 30 & 2.6549e-11 & 7.8500 & 1.6931e-11 & 7.8485 & 1.8471e-11 & 7.8342 \\ 
& 40 & 2.9058e-12 & 7.6900 & 1.8481e-12 & 7.6995 & 2.0281e-12 & 7.6791 \\ 
\bottomrule
\end{tabular*}\label{Table:E1b}
\end{table}
\par We next examine equation (\ref{example:1}) with $c = 1$, employing eighth-order TDCNCS and TDCCS compact schemes. The solutions and pointwise errors for $N=40$ at various times $t = 0, 0.25, 0.5, 0.75$, and $1$ are depicted in Figure (\ref{Figure:E1b}). Table (\ref{Table:E1b}) presents a summary of the $L^{\infty}$-error, $L^{1}$-error, $L^{2}$-error, and the order of convergence for both considered schemes at the final time $t = 1$.  For $N=20$, the $L^{\infty}$-error for TDCNCS is $1.6089\times 10^{-9}$, while for TDCCS, it is $6.4028\times 10^{-10}$. For $N=40$, the $L^{\infty}$-error for TDCNCS is $6.6641\times 10^{-12}$, whereas for TDCCS, it is $2.9058\times 10^{-12}$. However, as $N$ increases, the error converges to machine epsilon, with no decrease in errors observed. No filters are used in all computations for this example.
\end{example} 
\newpage
\begin{example}\label{example:2}
\normalfont To evaluate the method's accuracy in addressing non-linear problems, we compute the classical soliton solution of the KdV equation.  
\begin{equation}
\left\{
\begin{aligned}
u_t - 3(u^2)_x+ u_{xxx}&=0, \quad x \in[-10,12], \quad t \ge 0, \\
u_0(x)&= -2\sech^2(x), \quad x \in[-10,12].
\end{aligned}
\right.  
 \end{equation}
\begin{table}[htbp!]
\centering
\captionof{table}{Errors and orders of convergence with eight-order methods for Example \ref{example:2} at $t = 0.5$.}
\setlength{\tabcolsep}{0pt}
\begin{tabular*}{\textwidth}{@{\extracolsep{\fill}} l *{8}{c} }
\midrule
\textbf{Scheme} & \textbf{N} & $\boldsymbol{L^{\infty}}$\textbf{-error} & \textbf{Rate} & $\boldsymbol{L^{1}}$\textbf{-error} & \textbf{Rate} & $\boldsymbol{L^{2}}$\textbf{-error} & \textbf{Rate}\\
\toprule
\textbf{TDCNCS} & 20 & 5.4854e-01 & - & 1.4320e-01 & - & 2.1954e-01 & - \\
& 40 & 1.2988e-02 & 5.4003 & 3.6053e-03 & 5.3118 & 4.3622e-03 & 5.6533  \\
& 60 & 3.2825e-04 & 9.0711 & 5.7634e-05 & 10.2008 & 8.2101e-05 & 9.7981  \\
& 80 & 3.3159e-05 & 7.9687 & 5.6523e-06 & 8.0716 & 7.9914e-06 & 8.0978 \\
& 100 & 5.6867e-06 & 7.9017 & 9.0292e-07 & 8.2197 & 1.2949e-06 & 8.1558 \\
& 120 & 1.3256e-06 & 7.9874 & 2.0602e-07 & 8.1046 & 2.9639e-07 & 8.0875 \\
& 140 & 3.7695e-07 & 8.1575 & 5.8893e-08 & 8.1236 & 8.4988e-08 & 8.1036 \\
& 160 & 1.2706e-07 & 8.1438 & 2.0016e-08 & 8.0819 & 2.8872e-08 & 8.0851 \\
\midrule
\textbf{TDCCS} & 20 & 2.0778e-02 & - & 4.0951e-03 & - & 6.6827e-03 & - \\
& 40 & 2.6255e-04 & 6.3064 & 4.6850e-05 & 6.4497 & 7.2033e-05 & 6.5356 \\
& 60 & 1.7256e-05 & 6.7140 & 2.9090e-06 & 6.8542 & 4.2206e-06 & 6.9973 \\
& 80 & 2.3533e-06 & 6.9255 & 3.8793e-07 & 7.0034 & 5.7101e-07 & 6.9533 \\
& 100 & 4.8859e-07 & 7.0450 & 8.2305e-08 & 6.9479 & 1.1756e-07 & 7.0829 \\
& 120 & 1.3222e-07 & 7.1692 & 2.2022e-08 & 7.2311 & 3.1553e-08 & 7.2138 \\
& 140 & 4.2335e-08 & 7.3876 & 7.1750e-09 & 7.2750 & 1.0300e-08 & 7.2626  \\
& 160 & 1.7606e-08 & 6.5705 & 3.0410e-09 & 6.4285 & 4.2286e-09 & 6.6671 \\
\bottomrule
\end{tabular*}\label{Table:E2}
\end{table}
 
The exact solution is given by $u(x,t)=-2\sech^2(x-4t)$ for $x \in [-10,12]$ and $t \in [0,0.5]$. Employing periodic boundary conditions, we solve the problem (\ref{example:2}) using eighth-order TDCNCS and TDCCS schemes with a step size of $\Delta t = 0.01 \Delta x^3$. The resulting solutions and errors for $N=80$ at different time points $t = 0, 0.25, 0.5$ are visualized in Figure (\ref{Figure:E2a}), while Table (\ref{Table:E2}) offers a detailed comparison of their $L^{\infty}$-, $L^{1}$- errors and $L^{2}$-errors and their respective convergence rates for various grid sizes ($N$ ranging from 20 to 160). Tabulations confirm TDCCS's advantage in delivering better accuracy. It consistently outperformed TDCNCS, delivering errors at least an order of magnitude lower across all tested grid sizes. This remarkable precision, however, comes at a computational cost. TDCCS demands significantly more resources, with its runtime exceeding TDCNCS by at least a factor of two. While TDCCS exhibits a convergence rate of around the seventh order for larger $N$ values, TDCNCS maintains a consistent rate close to the theoretical eighth order. Using TDCNCS with $N=160$ yields an error of $1.2706\times 10^{-8}$, which closely corresponds to the error obtained by TDCCS with $N=100$ given equivalent computational time. It is important to note that further increasing grid points beyond $N = 160$ doesn't yield any improvement in accuracy due to limitations in machine precision. In summary, TDCCS offers superior accuracy but with substantially higher computational demands.
\begin{figure}[htbp!]
  \begin{minipage}[b]{0.45\linewidth}
    \centering
    \includegraphics[trim=0cm 0cm 0cm 0cm, clip=true,width=\linewidth]{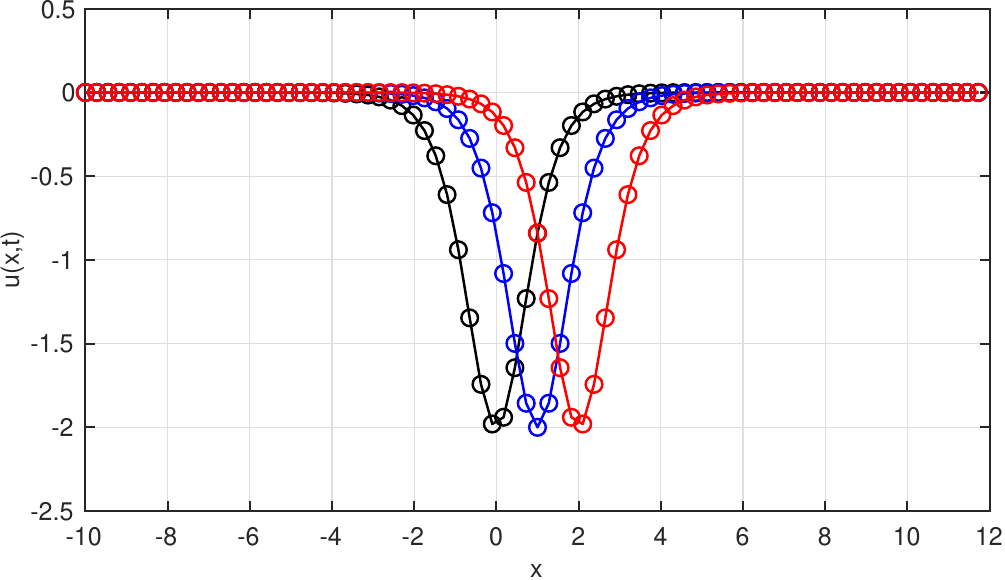}
  \end{minipage}\hfill
    \begin{minipage}[b]{0.45\linewidth}
    \centering
    \includegraphics[trim=0cm 0cm 0cm 0cm, clip=true,width=\linewidth]{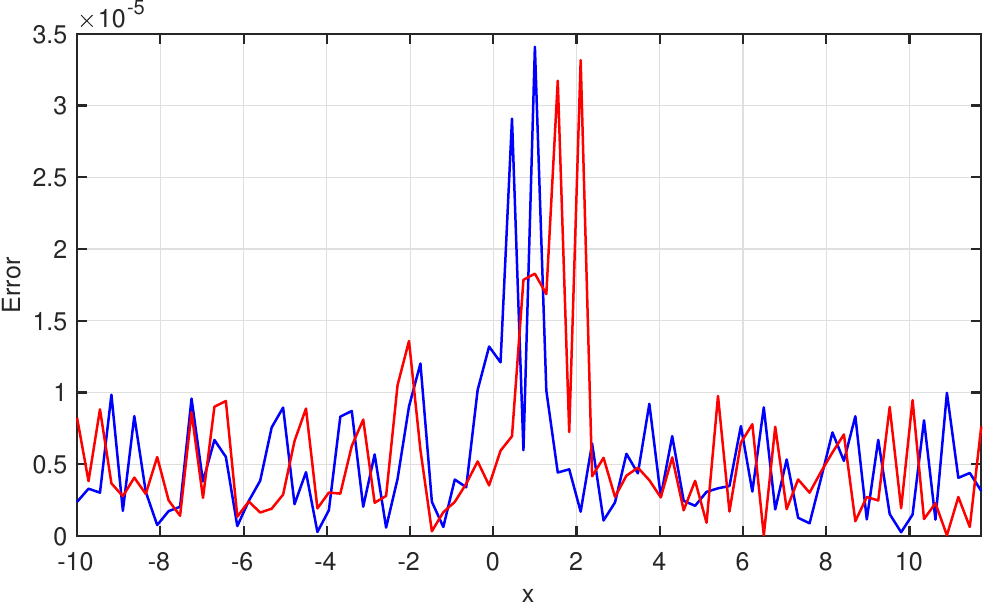}
  \end{minipage}\hfill
  \begin{minipage}[b]{0.45\linewidth}
    \centering
    \includegraphics[trim=0cm 0cm 0cm 0cm, clip=true,width=\linewidth]{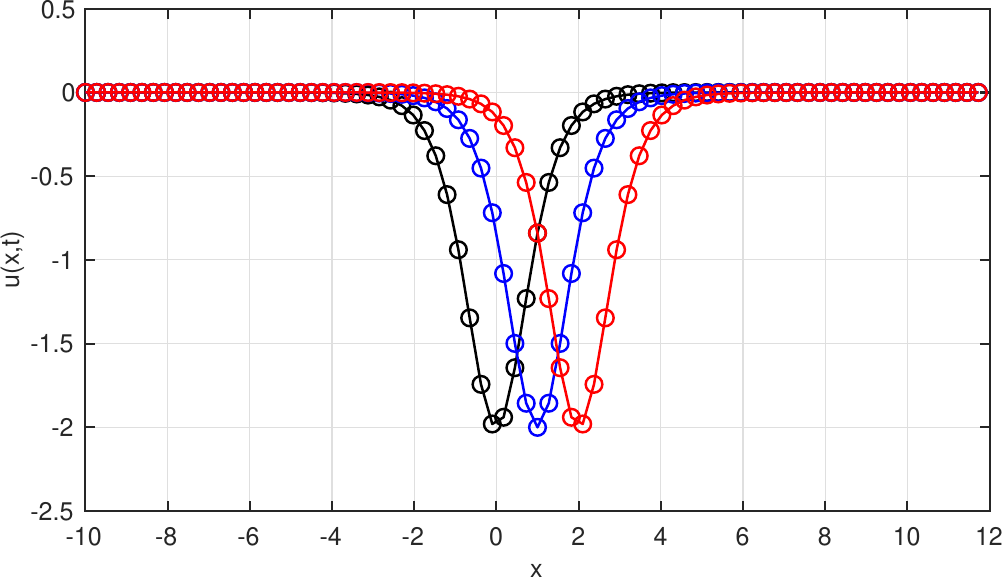}
  \end{minipage}\hfill
  \begin{minipage}[b]{0.45\linewidth}
    \centering
    \includegraphics[trim=0cm 0cm 0cm 0cm, clip=true,width=\linewidth]{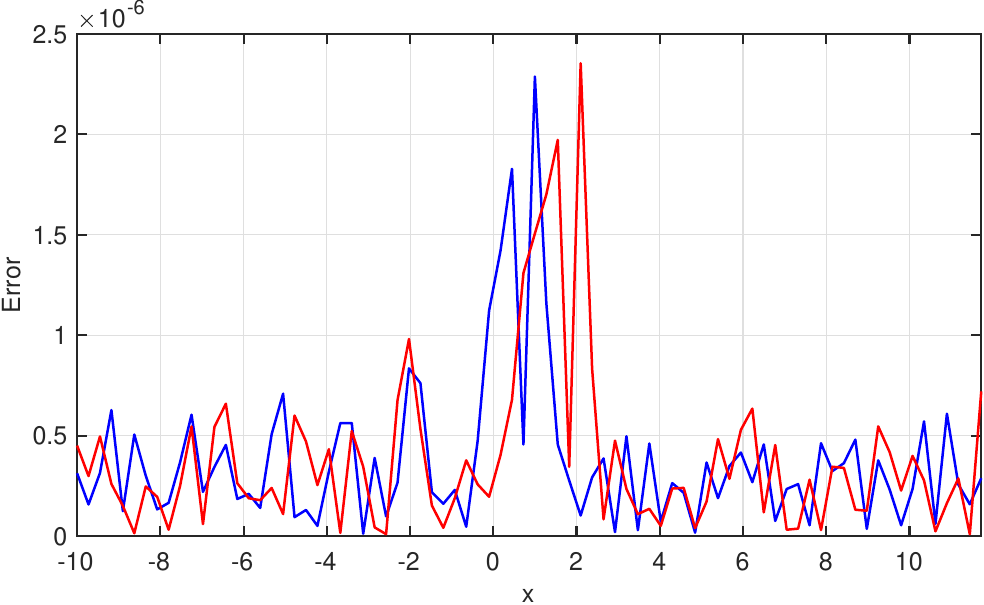}
  \end{minipage}\hfill
  \caption{Solutions and errors for Example \ref{example:2} at $t =0.25$ and $0.5$. The first and second rows are obtained by TDCNCS-T8 and TDCCS-T8, respectively. Left column shows the exact (solid line) and numerical values ($\bigcirc$) at $t =$0 (black), $t =$ 0.25 (blue), $t =$ 0.5 (red). Right column shows the corresponding errors.}\label{Figure:E2a}
\end{figure}
\subsection*{Performance analysis of filter} This study investigates the influence of LPSF filters on the precision of numerical solutions, aiming to mitigate spurious wave amplitudes while preserving phase accuracy. Specifically, we examine the effects of 8th, 10th, and 12th order filters on solutions with $N=40$ grid points, $t=0.1$ temporal duration (time steps=60), with both TDCNCS and TDCCS schemes depicted in Figure (\ref{Figure:E2b}). The filter's impact is evaluated every 2 and 10-time step, revealing superior performance for unfiltered TDCCS compared to its filtered counterparts and TDCNCS with or without filtering. 
\begin{figure}[htbp!]
\begin{minipage}[b]{0.24\linewidth}
    \centering
    \includegraphics[trim=0cm 0cm 0cm 0cm, clip=true,width=\linewidth]{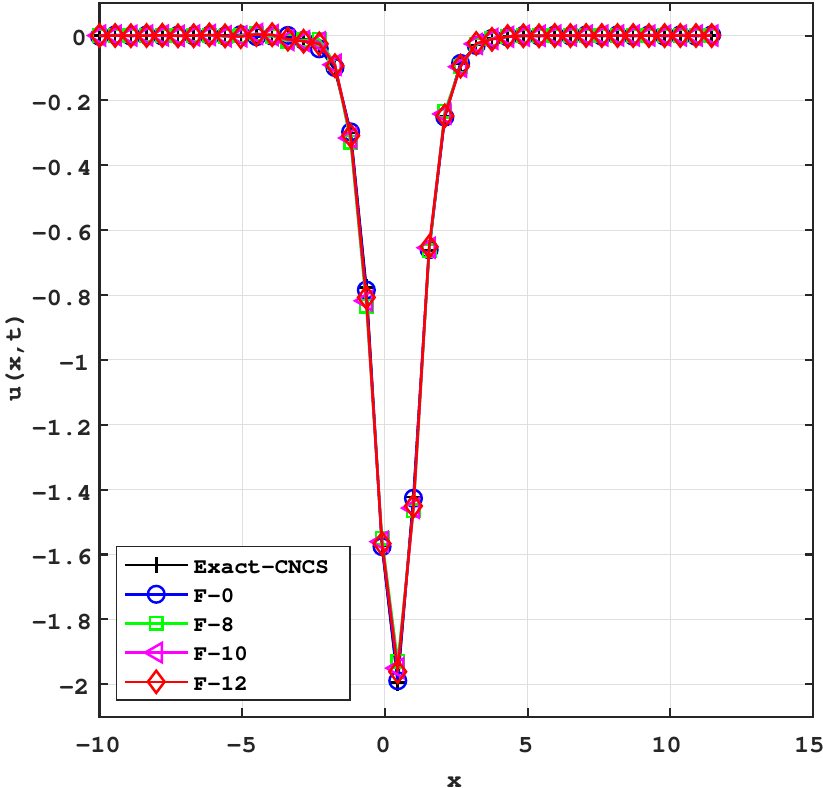}
    \subcaption*{\normalsize{\centering (a) TDCNCS, filtered every 2 time steps}}
  \end{minipage}\hfill
  \begin{minipage}[b]{0.24\linewidth}
    \centering
    \includegraphics[trim=0cm 0cm 0cm 0cm, clip=true,width=\linewidth]{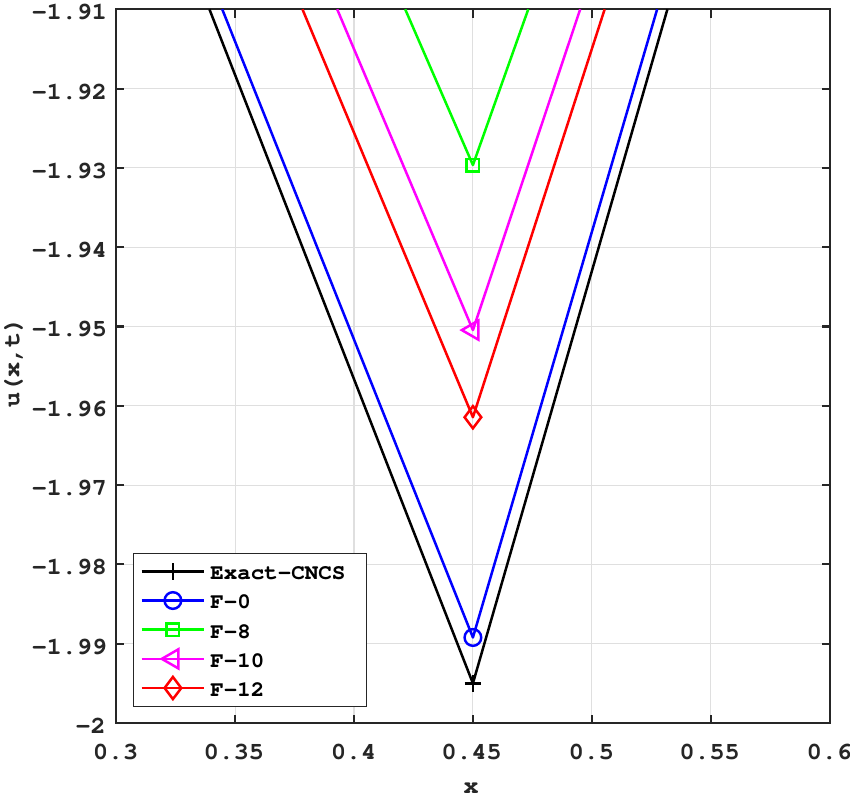}
    \subcaption*{\normalsize{\centering (b) Zoom in near smooth extrema}}
  \end{minipage}\hfill
  \begin{minipage}[b]{0.24\linewidth}
    \centering
    \includegraphics[trim=0cm 0cm 0cm 0cm, clip=true,width=\linewidth]{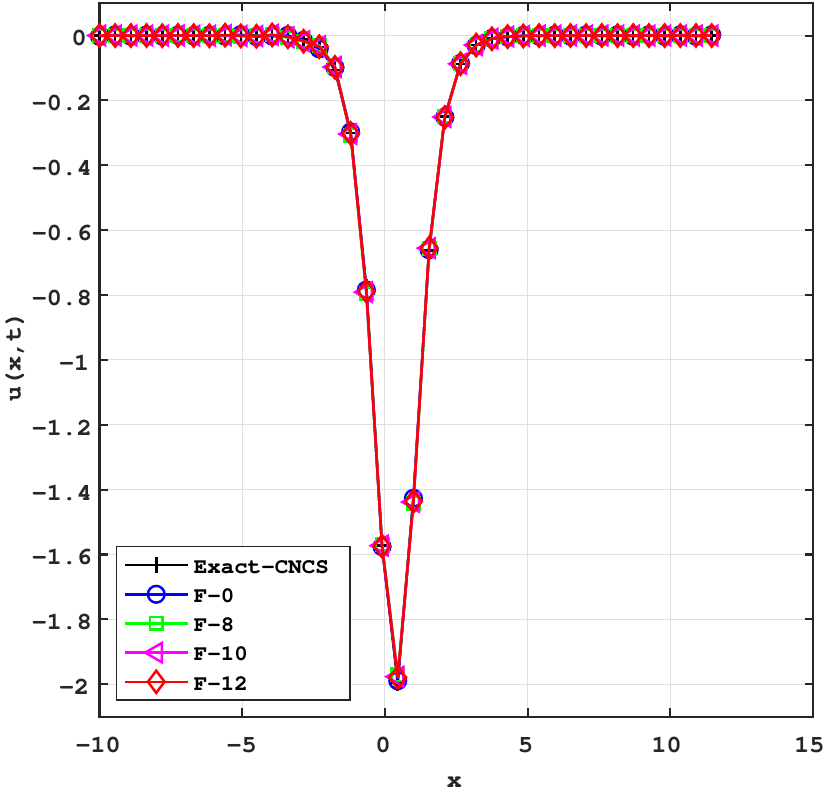}
    \subcaption*{\normalsize{\centering (c) TDCNCS, filtered every 10 time steps}}
  \end{minipage}
    \begin{minipage}[b]{0.24\linewidth}
    \centering
    \includegraphics[trim=0cm 0cm 0cm 0cm, clip=true,width=\linewidth]{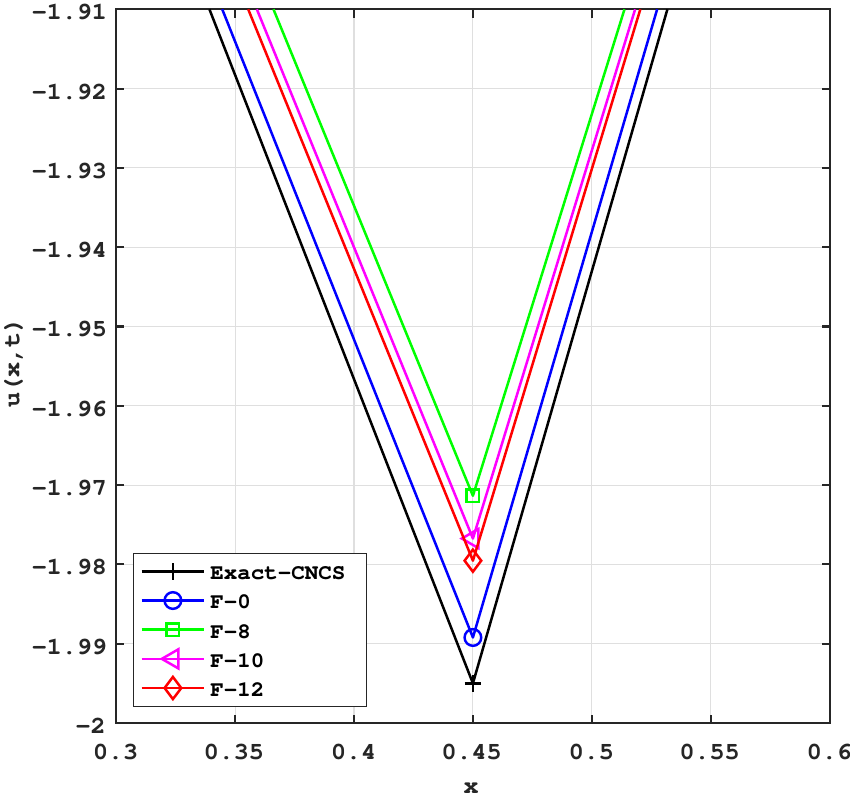}
    \subcaption*{\normalsize{\centering (d) Zoom in near smooth extrema}}
  \end{minipage}
  \begin{minipage}[b]{0.24\linewidth}
    \centering
    \includegraphics[trim=0cm 0cm 0cm 0cm, clip=true,width=\linewidth]{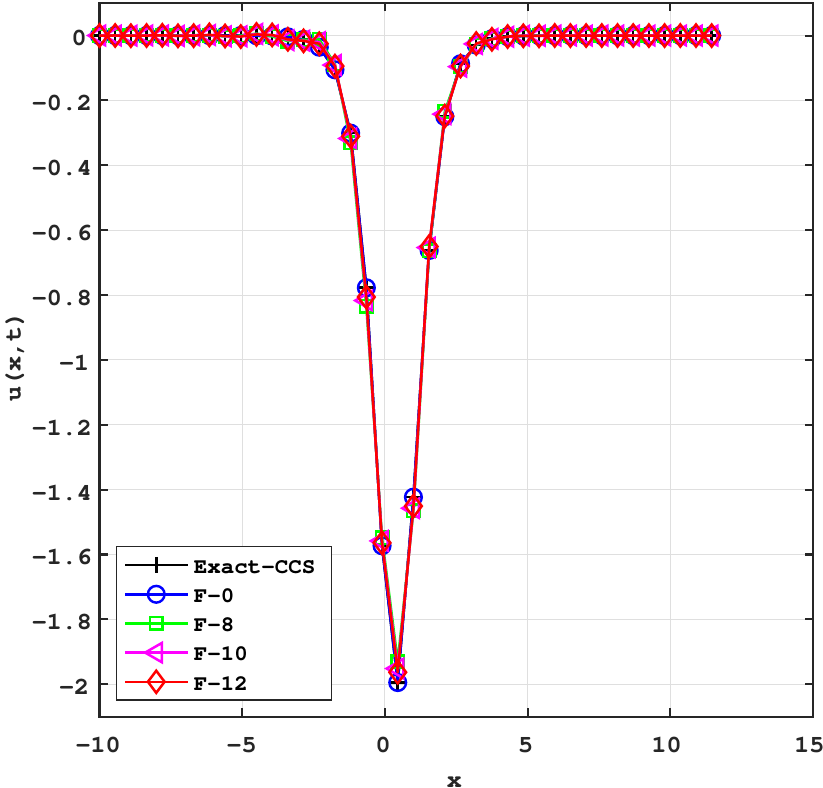}
    \subcaption*{\normalsize{\centering (e) TDCCS, filtered every 2 time steps}}
  \end{minipage}\hfill
  \begin{minipage}[b]{0.24\linewidth}
    \centering
    \includegraphics[trim=0cm 0cm 0cm 0cm, clip=true,width=\linewidth]{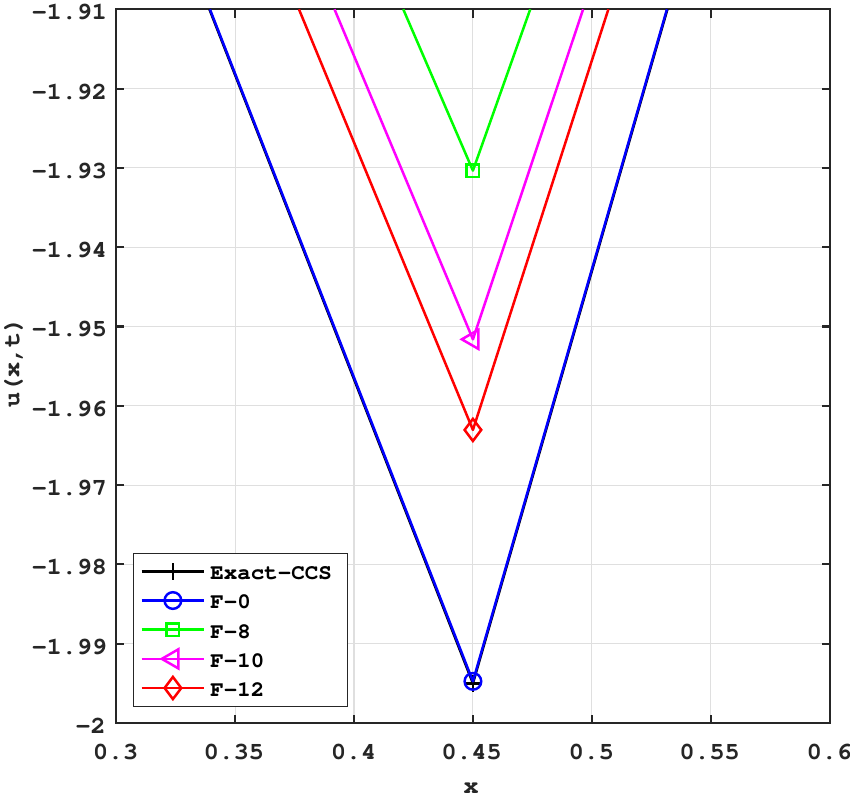}
    \subcaption*{\normalsize{\centering (f) Zoom in near smooth extrema}}
  \end{minipage}\hfill
      \begin{minipage}[b]{0.24\linewidth}
    \centering
    \includegraphics[trim=0cm 0cm 0cm 0cm, clip=true,width=\linewidth]{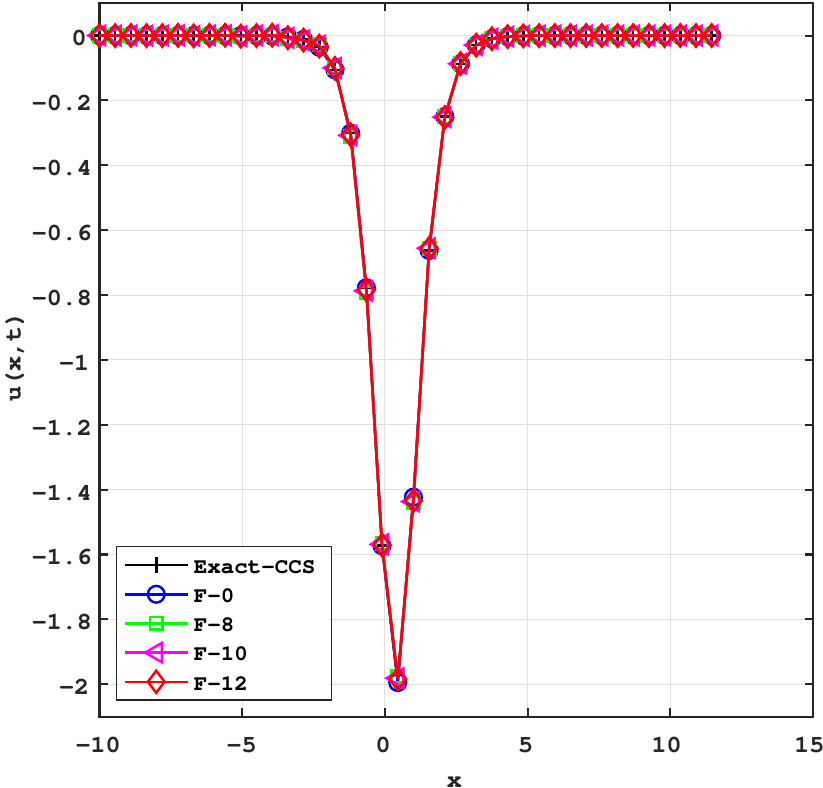}
    \subcaption*{\normalsize{\centering (g) TDCCS, filtered every 10 time steps}}
  \end{minipage}
    \begin{minipage}[b]{0.24\linewidth}
    \centering
    \includegraphics[trim=0cm 0cm 0cm 0cm, clip=true,width=\linewidth]{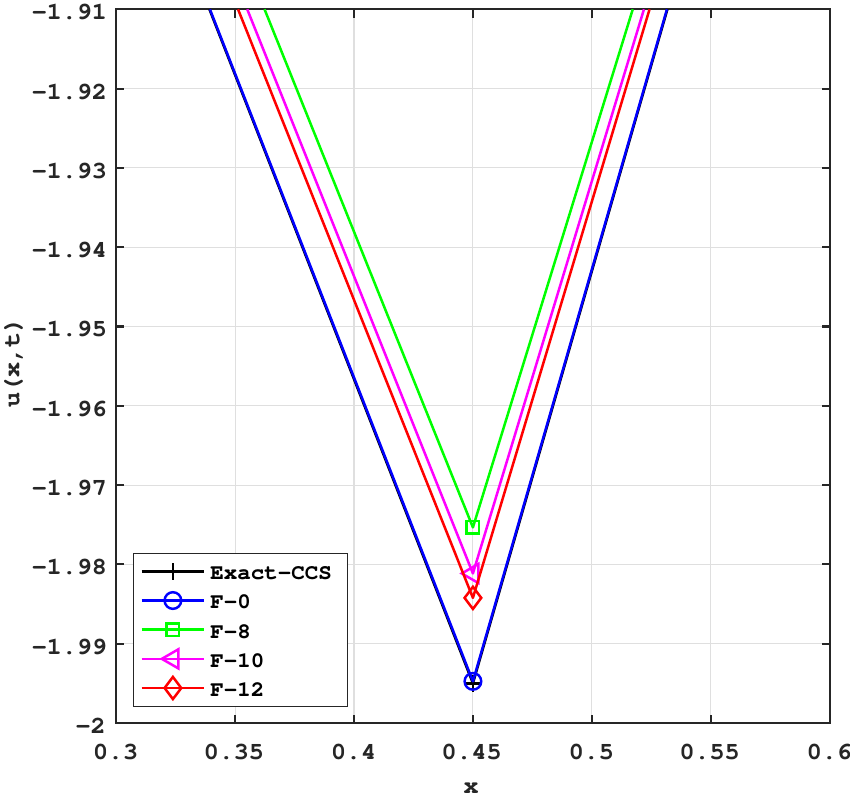}
    \subcaption*{\normalsize{\centering (h) Zoom in near smooth extrema}}
  \end{minipage}
  \caption{Analysis of filter at $N=40$, $t=0.1$ for Example \ref{example:2}. }\label{Figure:E2b}
\end{figure}
\end{example}
\begin{example}\label{example:3}
\normalfont We solve the following nonlinear KdV equations. To observe the effectiveness of our method in handling a nonlinear problem with a low coefficient for the third derivative term, we calculate the classical soliton solutions of the KdV equation.
\begin{equation}
u_t + \left(\frac{u}{2}^2\right)_x+ \epsilon u_{xxx} =0.
 \end{equation}
 \subsection*{Single soliton propagation} In the case of a single soliton, the initial condition is as follows,
 \begin{equation}\label{IC:3a}
u_0(x) = 3c\sech^2(k(x-x_0)), \quad x \in[0,2],
 \end{equation}
 with $k=0.5\sqrt{\frac{c}{\epsilon}}$, $c=0.3$, $x_0 = 0.5$ and $\epsilon=5 \times 10^{-4}$. The solution to this problem is a solitary wave, referred to as a \textit{soliton}, moving to the right with a speed $c$ given by $u(x,t)=3c\sech^2(k[(x-x_0)-ct])$. The solution is computed with periodic boundary conditions for $x \in [0, 2]$, $N=80$, $t \in [0, 3]$ using eighth-order TDCNCS and TDCCS schemes with step size $\Delta t = 0.01 \Delta x^3$. The solution and error at time $t =0, 1, 2$ and $3$ is shown in Figure (\ref{Figure:E3a}). The TDCCS scheme exhibits significantly lower absolute error compared to TDCNCS, resulting in at least an order of magnitude improvement in accuracy.
\begin{figure}[htbp!]
   \begin{minipage}[b]{0.45\linewidth}
    \centering
    \includegraphics[trim=0cm 0cm 0cm 0cm, clip=true,width=\linewidth]{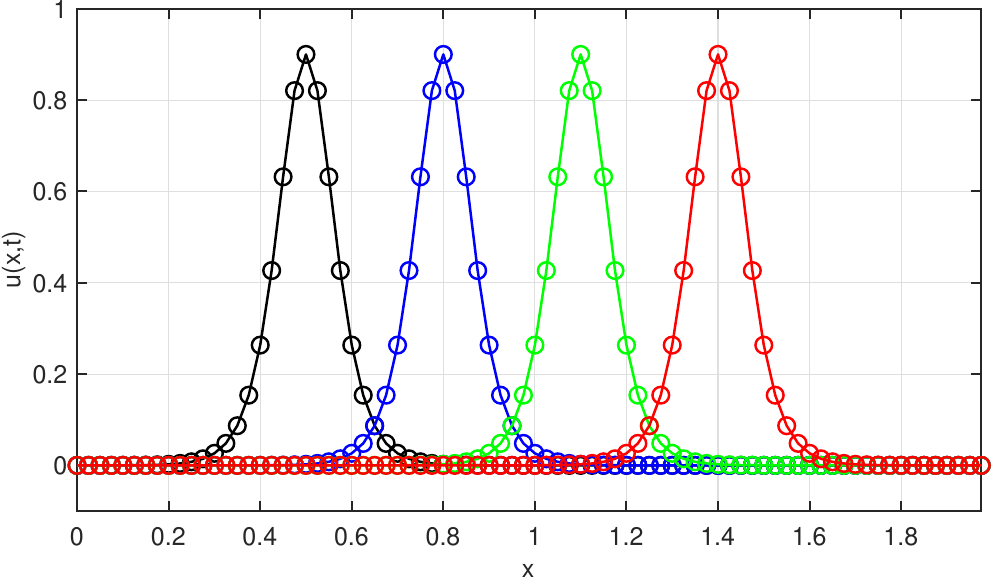}
  \end{minipage}\hfill
    \begin{minipage}[b]{0.45\linewidth}
    \centering
    \includegraphics[trim=0cm 0cm 0cm 0cm, clip=true,width=\linewidth]{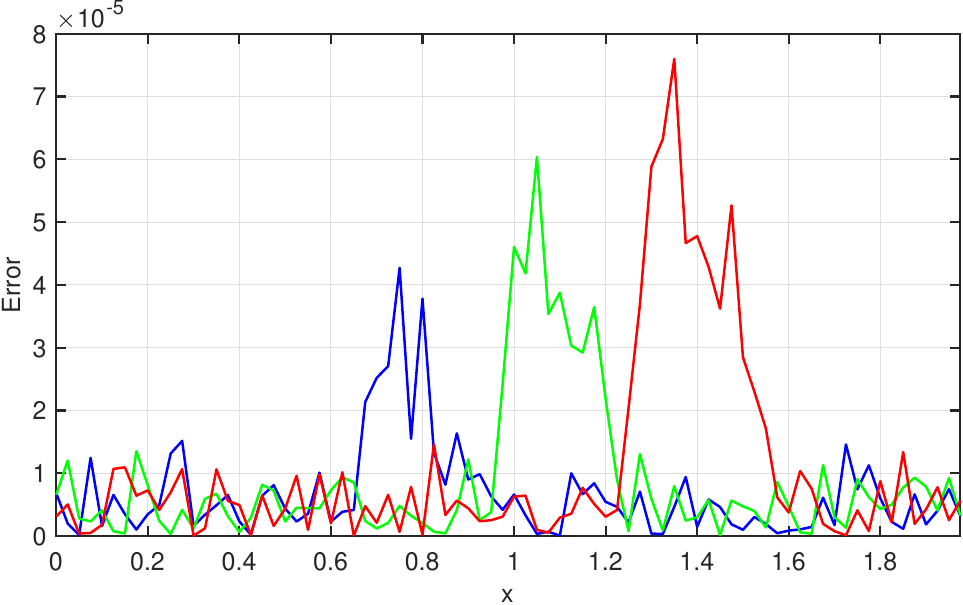}
  \end{minipage}
    \begin{minipage}[b]{0.45\linewidth}
    \centering
    \includegraphics[trim=0cm 0cm 0cm 0cm, clip=true,width=\linewidth]{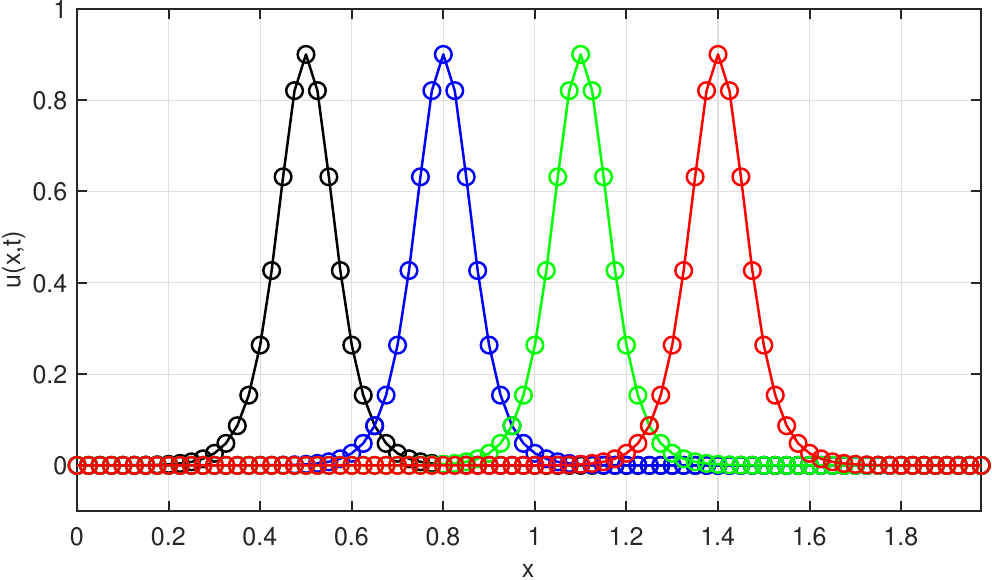}
  \end{minipage}\hfill
    \begin{minipage}[b]{0.45\linewidth}
    \centering
    \includegraphics[trim=0cm 0cm 0cm 0cm, clip=true,width=\linewidth]{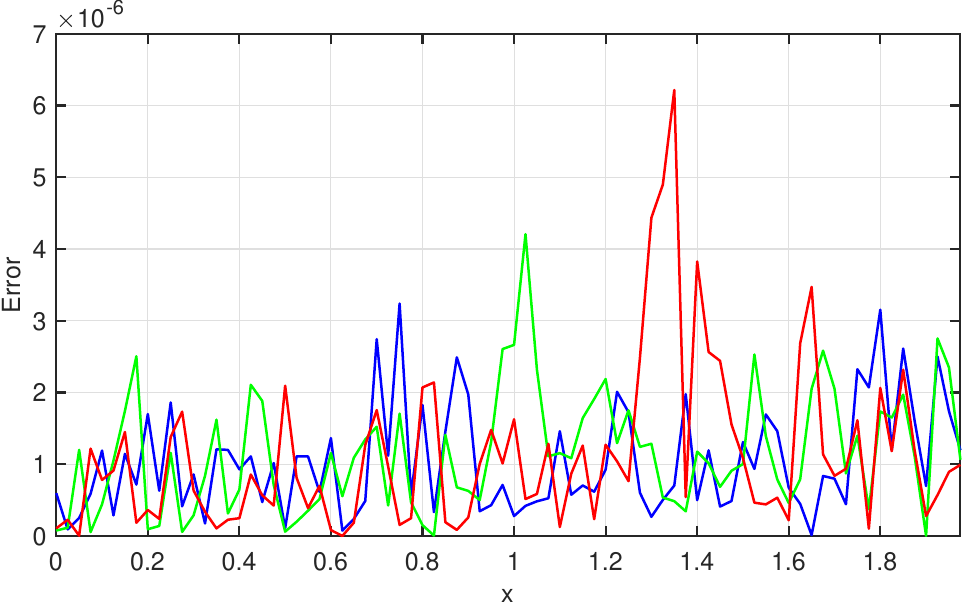}
  \end{minipage}
  \caption{ Solutions and errors for
initial condition \ref{IC:3a} of Example \ref{example:3}. The first and second rows are obtained by TDCNCS-T8 and TDCCS-T8, respectively. The left column shows the exact (solid line) and numerical values ($\bigcirc$) at $t =$ 0 (black), $t =$ 1 (blue), $t =$ 2 (green), and $t =$ 3 (red). The right column shows the corresponding errors varying with time.}\label{Figure:E3a}
\end{figure}
\subsection*{Double soliton collision} In the case of a double soliton collision, the initial condition is as follows,
 \begin{equation}\label{IC:3b}
u_0(x) = 3c_1\sech^2(k_1(x-x_1))+3c_2\sech^2(k_2(x-x_2)), \quad x \in[0,2],
 \end{equation}
 with $k_j=0.5\sqrt{\frac{c_j}{\epsilon}}$ for $j=1,2$, $c_1=0.3$, $c_2=0.1$, $x_1 = 0.4$, $x_2 = 0.8$ and $\epsilon=4.84 \times 10^{-4}$. The solution is computed for $t \in [0, 4]$, $N=100$ with periodic boundary conditions using eighth-order TDCNCS and TDCCS schemes with step size $\Delta t = 0.01 \Delta x^3$. The solution at time $t =0,1,2$ and solution contour up to $t=4$ is shown in Figure (\ref{Figure:E3b}).
 \begin{figure}[htbp!]
   \begin{minipage}[b]{0.24\linewidth}
    \centering
    \includegraphics[trim=0cm 0cm 0cm 0cm, clip=true,width=\linewidth]{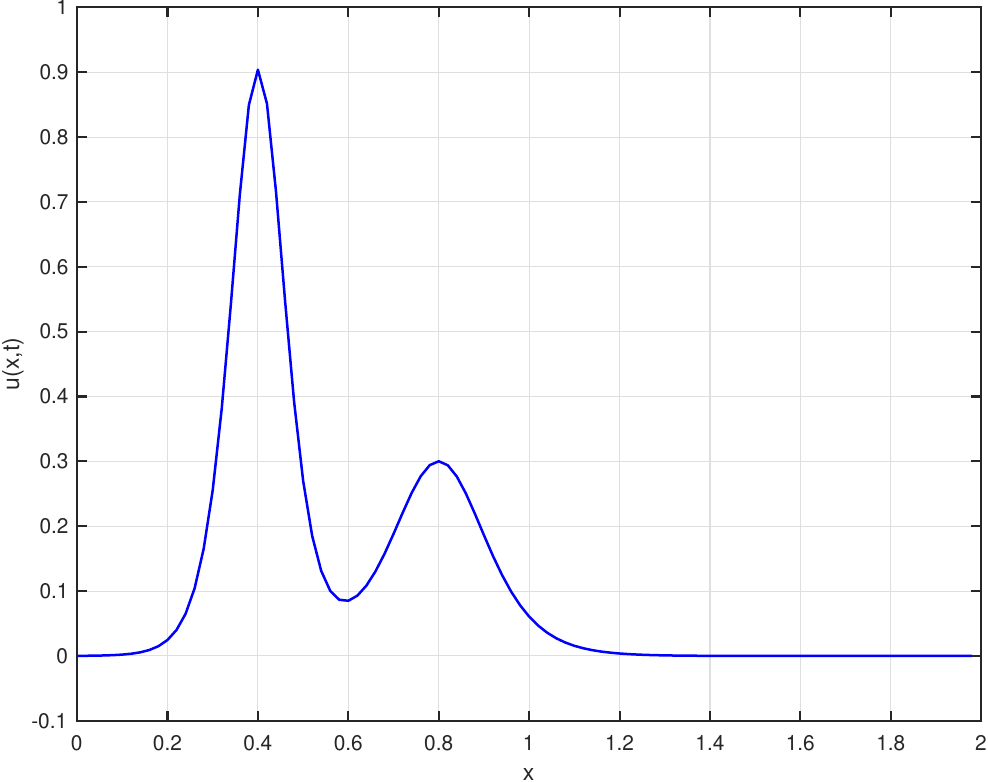}
    \subcaption*{\normalsize{\centering (a) t=0}}
  \end{minipage}\hfill
    \begin{minipage}[b]{0.24\linewidth}
    \centering
    \includegraphics[trim=0cm 0cm 0cm 0cm, clip=true,width=\linewidth]{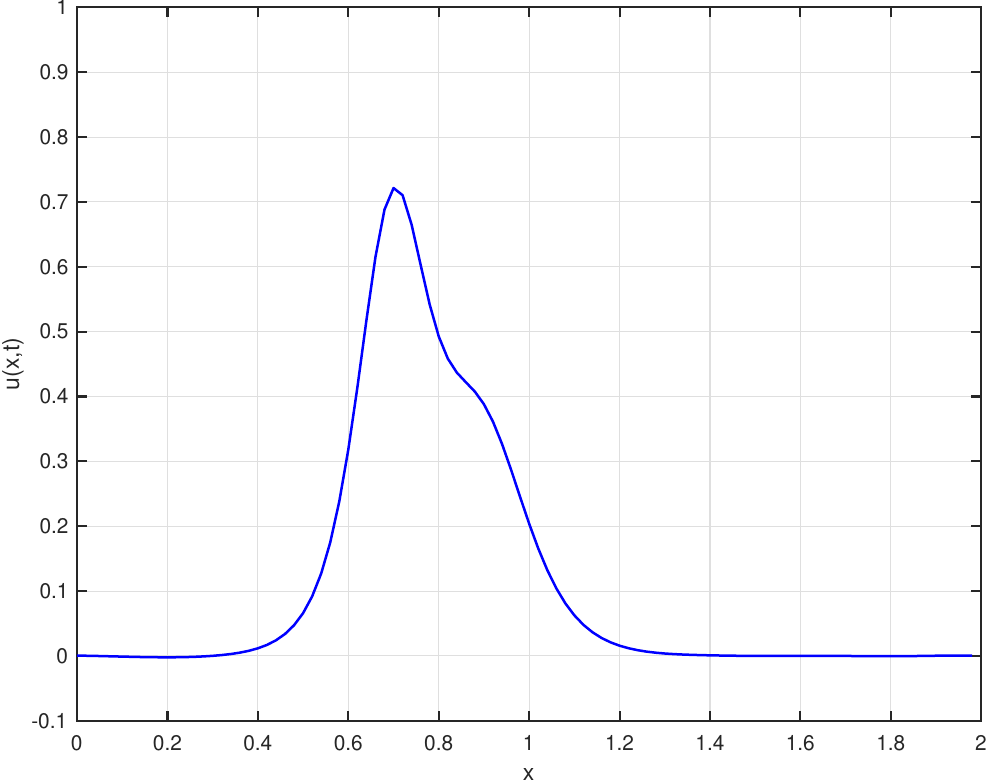}
    \subcaption*{\normalsize{\centering (b) t=1}}
  \end{minipage}\hfill
  \begin{minipage}[b]{0.24\linewidth}
    \centering
    \includegraphics[trim=0cm 0cm 0cm 0cm, clip=true,width=\linewidth]{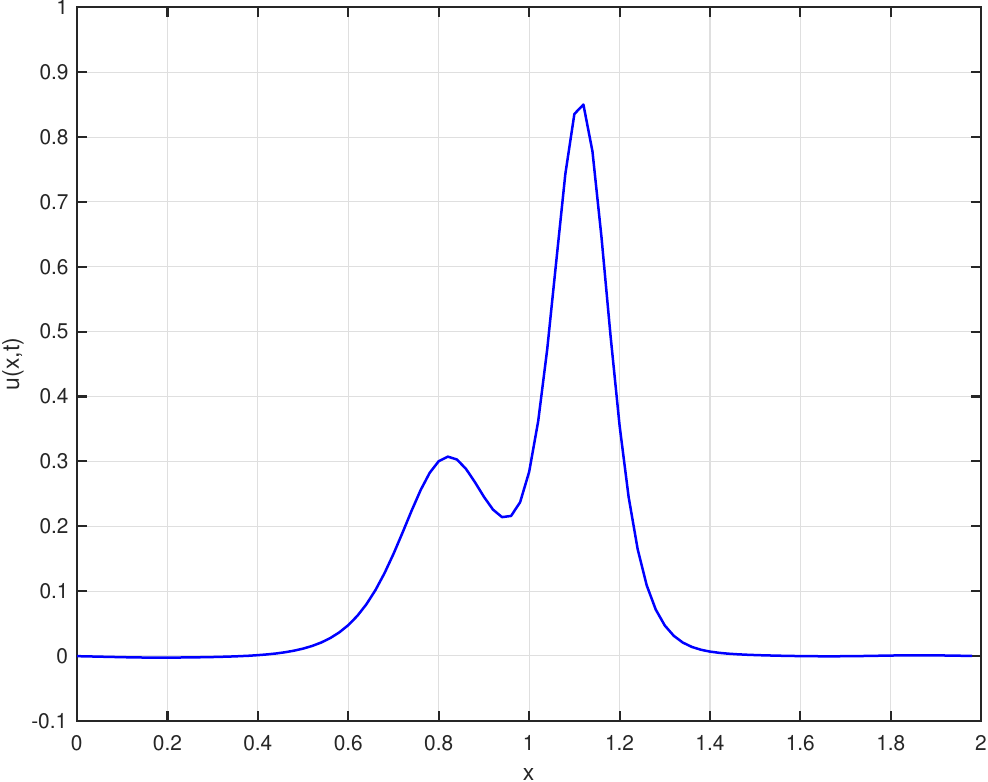}
    \subcaption*{\normalsize{\centering (c) t=2}}
  \end{minipage}\hfill
  \begin{minipage}[b]{0.24\linewidth}
    \centering
    \includegraphics[trim=0cm 0cm 0cm 0cm, clip=true,width=\linewidth]{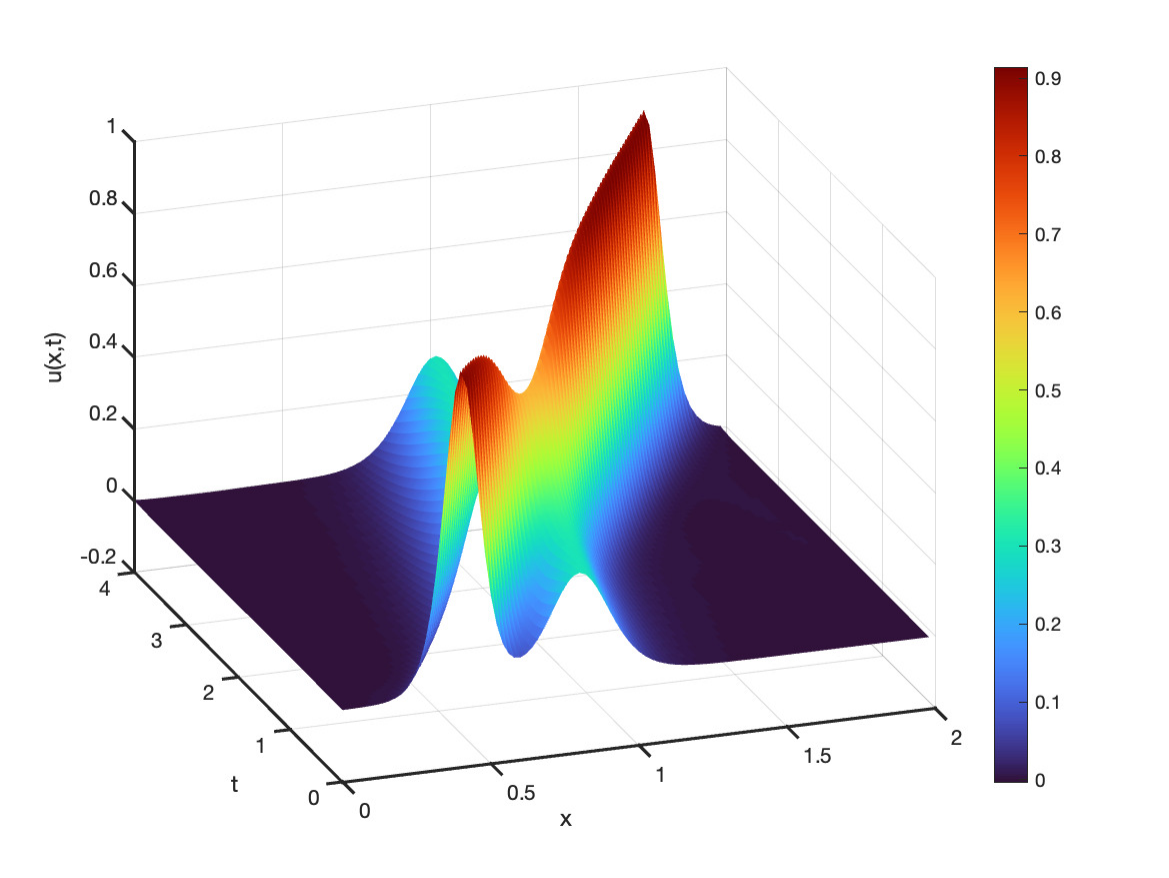}
    \subcaption*{\normalsize{\centering (d) t=4}}
  \end{minipage}\hfill
  \begin{minipage}[b]{0.24\linewidth}
    \centering
    \includegraphics[trim=0cm 0cm 0cm 0cm, clip=true,width=\linewidth]{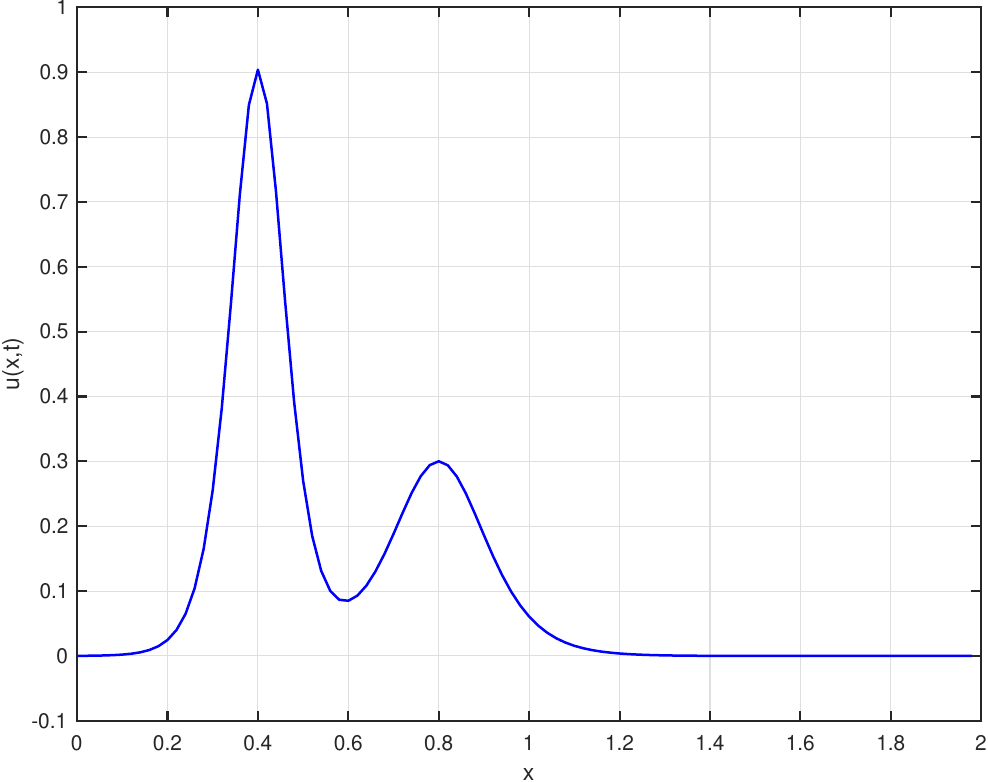}
    \subcaption*{\normalsize{\centering (e) t=0}}
  \end{minipage}\hfill
    \begin{minipage}[b]{0.24\linewidth}
    \centering
    \includegraphics[trim=0cm 0cm 0cm 0cm, clip=true,width=\linewidth]{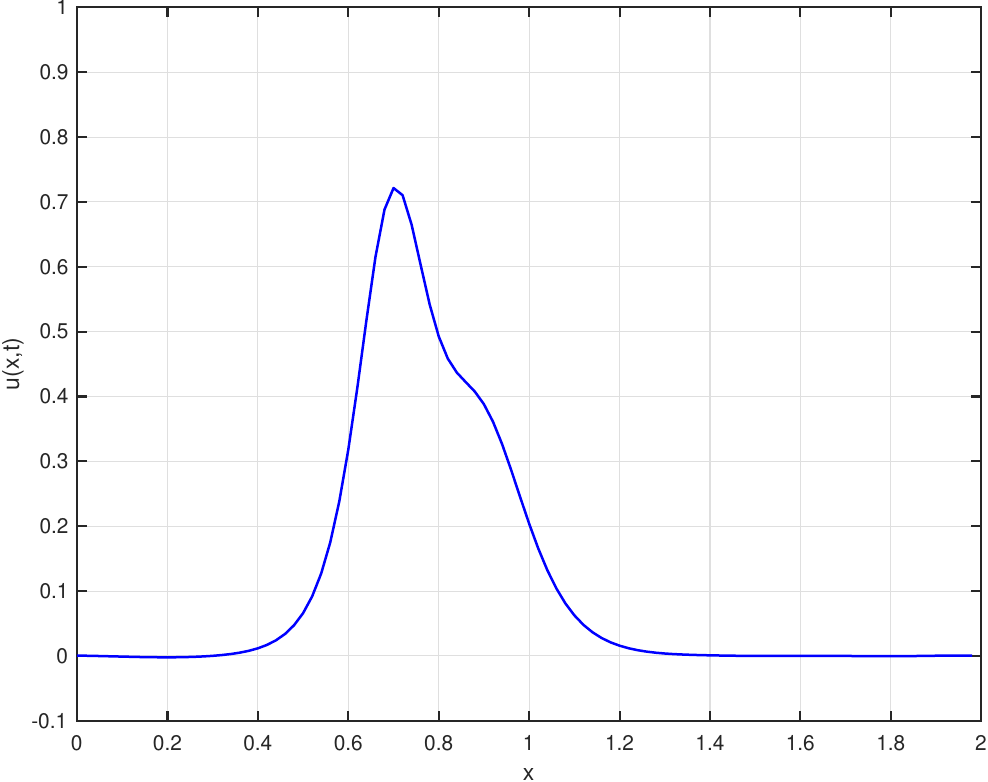}
    \subcaption*{\normalsize{\centering (f) t=1}}
  \end{minipage}\hfill
  \begin{minipage}[b]{0.24\linewidth}
    \centering
    \includegraphics[trim=0cm 0cm 0cm 0cm, clip=true,width=\linewidth]{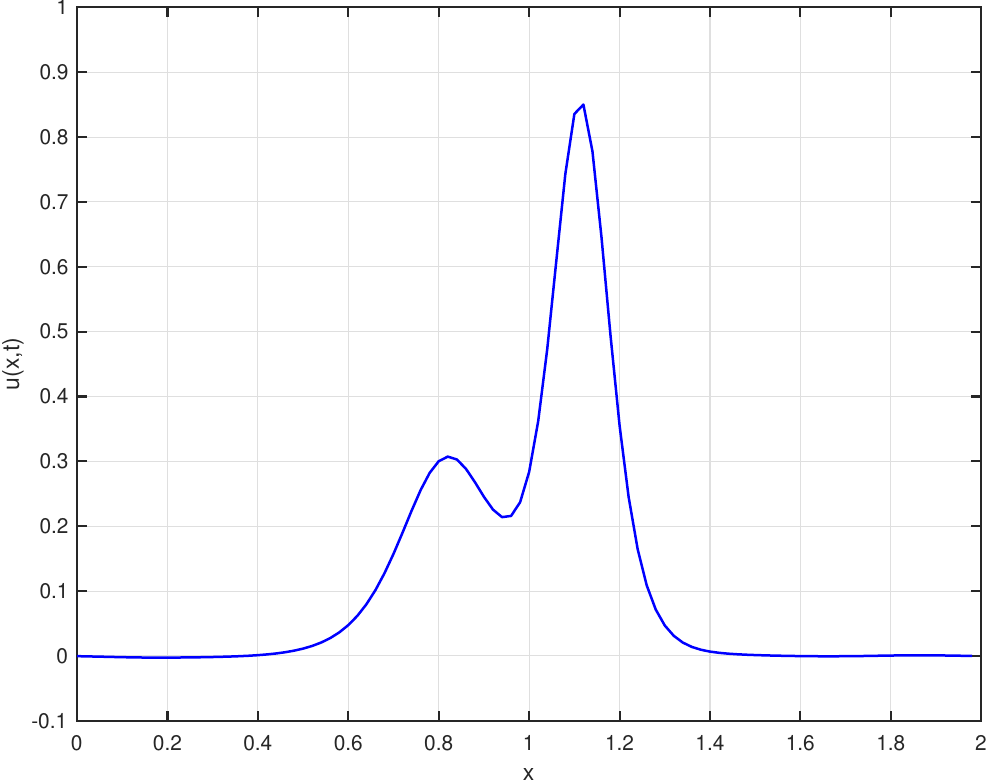}
    \subcaption*{\normalsize{\centering (g) t=2}}
  \end{minipage}\hfill
  \begin{minipage}[b]{0.24\linewidth}
    \centering
    \includegraphics[trim=0cm 0cm 0cm 0cm, clip=true,width=\linewidth]{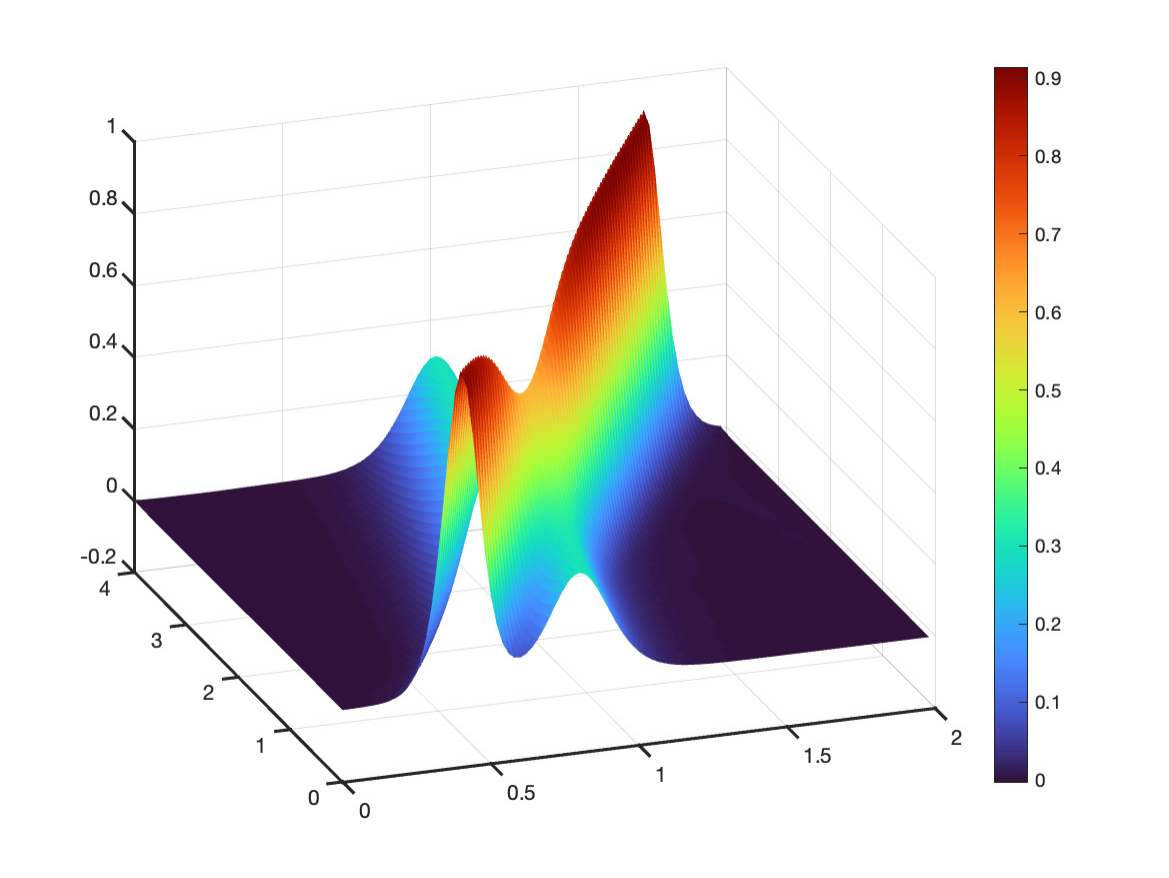}
    \subcaption*{\normalsize{\centering (h) t=4}}
  \end{minipage}\hfill
  \caption{Numerical solutions for
initial condition \ref{IC:3b} of Example \ref{example:3}.  The first row of four graphs is computed using TDCNCS and the second row of four graphs is computed using TDCCS. }\label{Figure:E3b}
\end{figure}
\subsection*{Triple soliton splitting} In the case of a triple soliton splitting, the initial condition is given by
 \begin{equation}\label{IC:3c}
u_0(x) = \frac{2}{3}\sech^2\left( \frac{x-1}{\sqrt{108 \epsilon}}\right), \quad x \in[0,3],
 \end{equation}
where $\epsilon=10^{-4}$. The solution is simulated over the interval $t \in [0, 4]$ using periodic boundary conditions, employing eighth-order TDCNCS and TDCCS schemes with a grid size of $N=150$ and time step $\Delta t = 0.5 \Delta x^2$ as in \cite{JV}. Figure (\ref{Figure:E3c}) illustrates the solution at time instances $t = 0, 1, 2,$ and  its contour plot up to $t=4$. Notably, TDCNCS exhibits small-scale numerical oscillations, whereas TDCCS produces a smooth numerical solution. To address the presence of high-frequency oscillations, we apply a 12th-order periodic filter with $\alpha = 0.4$ every 20 steps for TDCNCS and every 50 steps for TDCCS (time steps = 20000). Our experiments indicate that filtering every 50 steps does not sufficiently eliminate the small oscillations in TDCNCS.
 \begin{figure}[htbp!]
   \begin{minipage}[b]{0.24\linewidth}
    \centering
    \includegraphics[trim=0cm 0cm 0cm 0cm, clip=true,width=\linewidth]{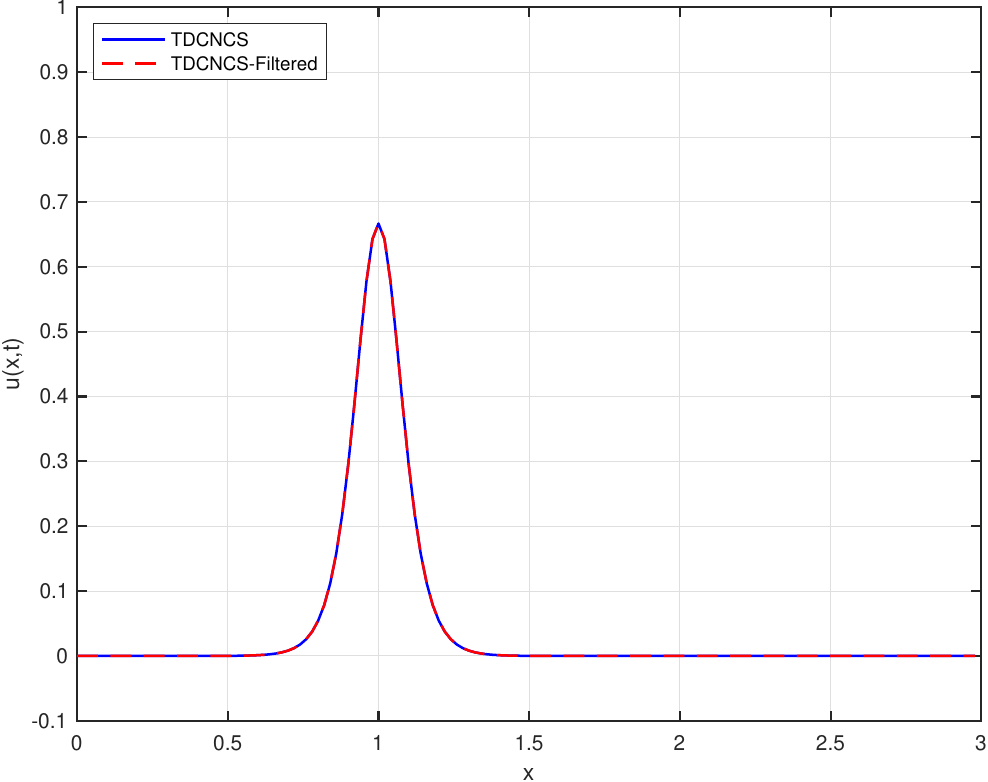}
    \subcaption*{\normalsize{\centering (a) t=0}}
  \end{minipage}\hfill
    \begin{minipage}[b]{0.24\linewidth}
    \centering
    \includegraphics[trim=0cm 0cm 0cm 0cm, clip=true,width=\linewidth]{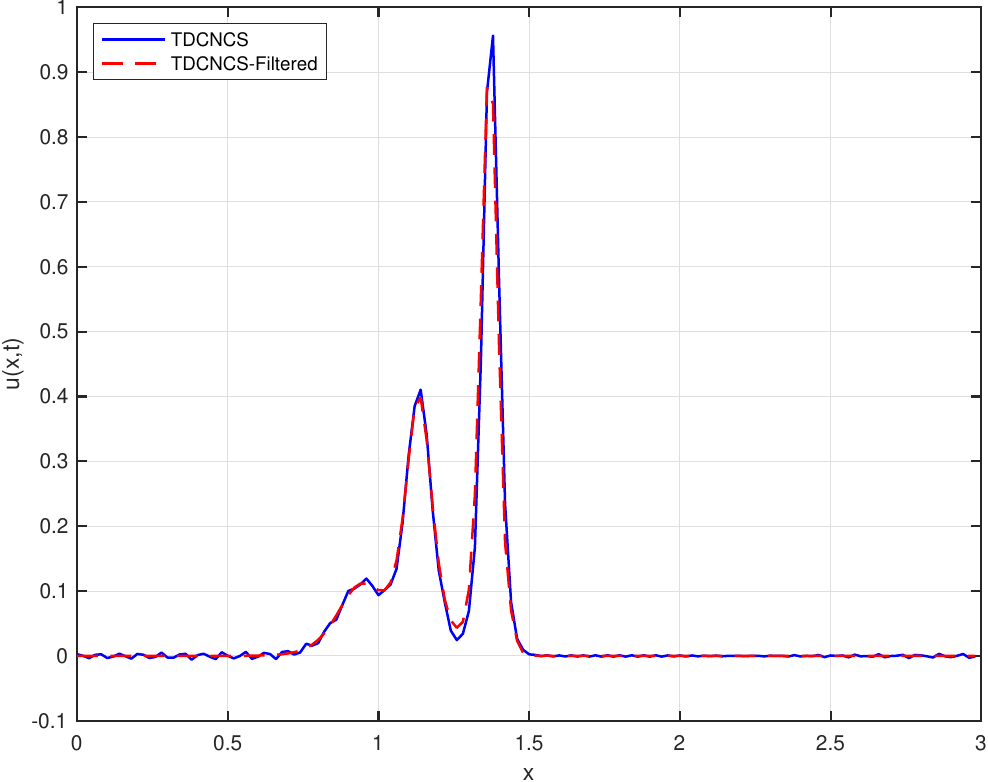}
    \subcaption*{\normalsize{\centering (b) t=1}}
  \end{minipage}\hfill
  \begin{minipage}[b]{0.24\linewidth}
    \centering
    \includegraphics[trim=0cm 0cm 0cm 0cm, clip=true,width=\linewidth]{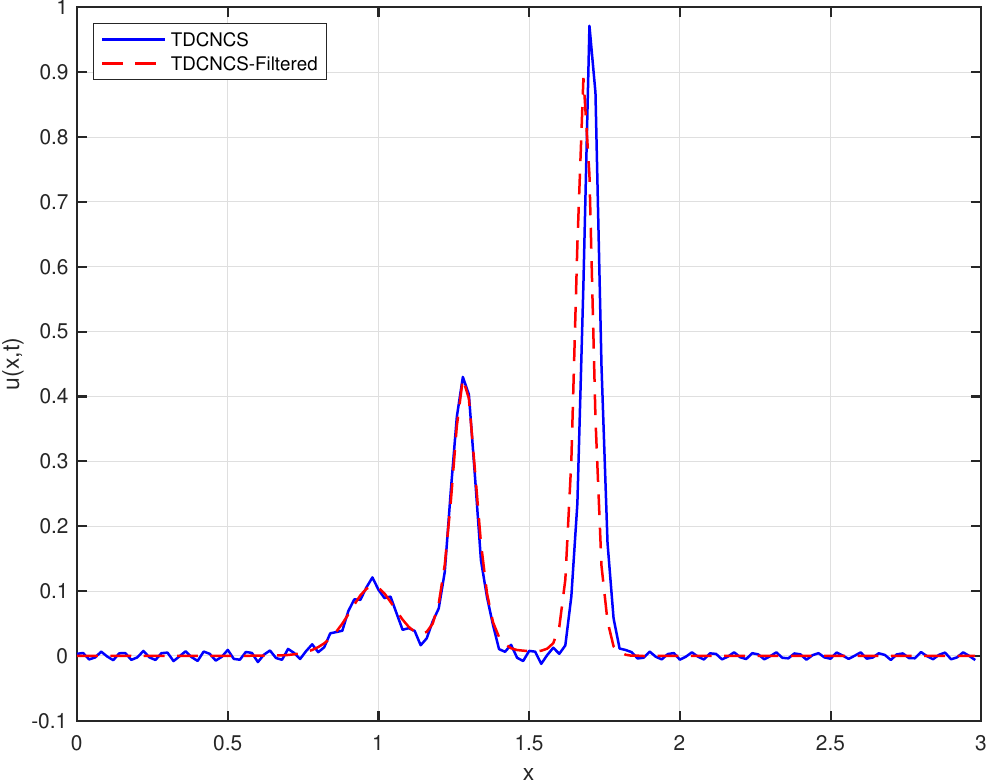}
    \subcaption*{\normalsize{\centering (c) t=2}}
  \end{minipage}\hfill
  \begin{minipage}[b]{0.24\linewidth}
    \centering
    \includegraphics[trim=0cm 0cm 0cm 0cm, clip=true,width=\linewidth]{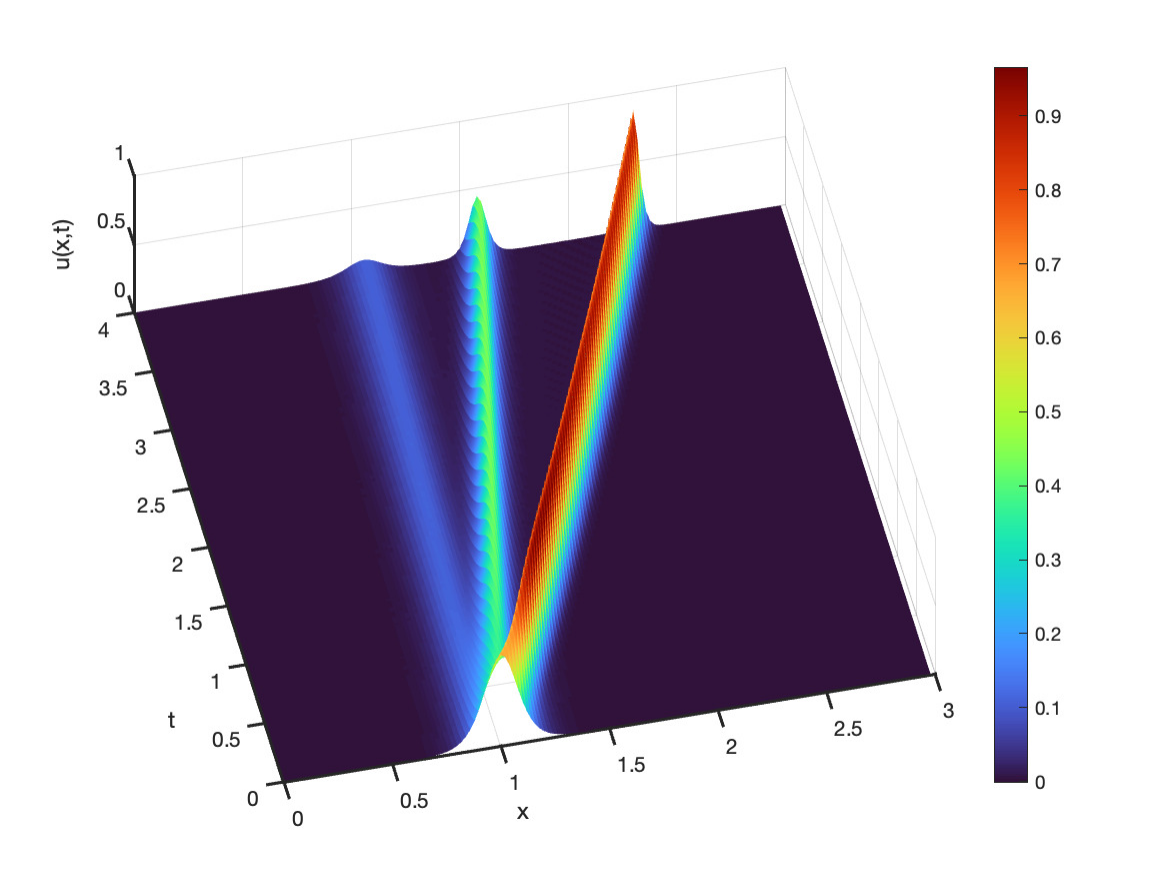}
    \subcaption*{\normalsize{\centering (d) t=4}}
  \end{minipage}\hfill
  \begin{minipage}[b]{0.24\linewidth}
    \centering
    \includegraphics[trim=0cm 0cm 0cm 0cm, clip=true,width=\linewidth]{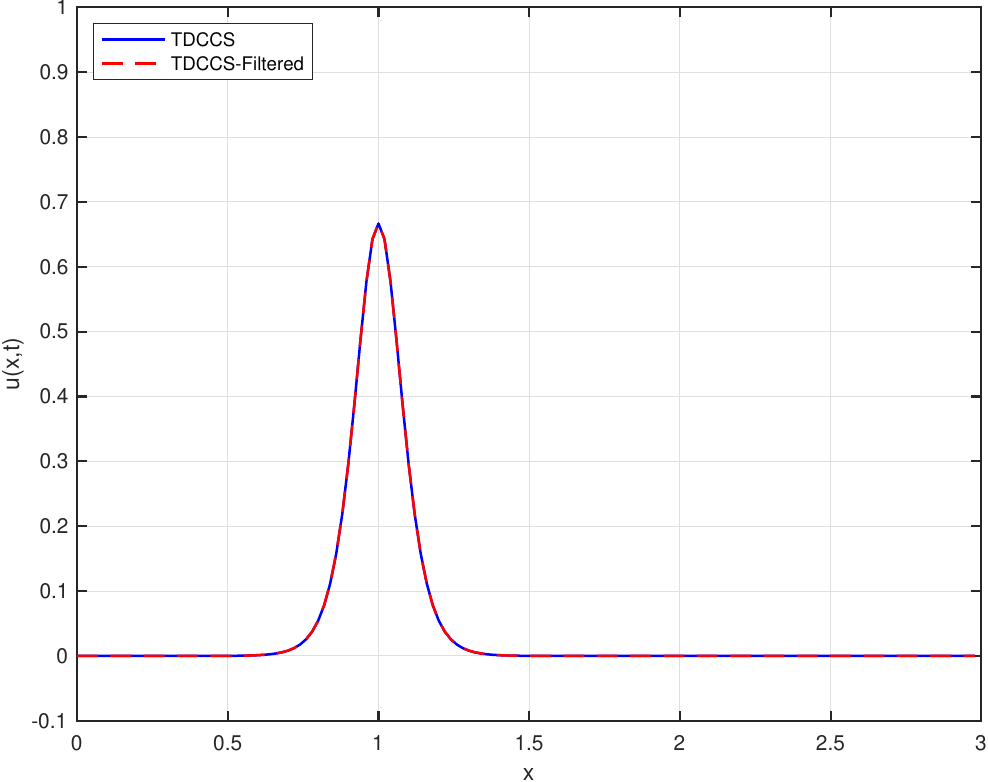}
    \subcaption*{\normalsize{\centering (e) t=0}}
  \end{minipage}\hfill
    \begin{minipage}[b]{0.24\linewidth}
    \centering
    \includegraphics[trim=0cm 0cm 0cm 0cm, clip=true,width=\linewidth]{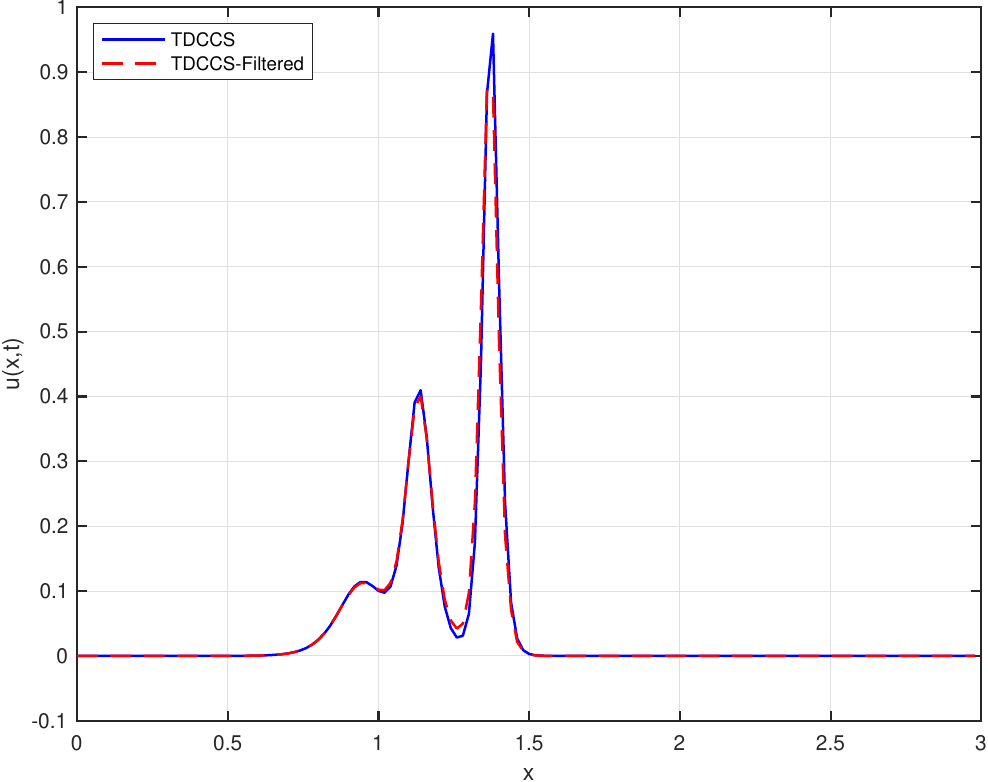}
    \subcaption*{\normalsize{\centering (f) t=1}}
  \end{minipage}\hfill
  \begin{minipage}[b]{0.24\linewidth}
    \centering
    \includegraphics[trim=0cm 0cm 0cm 0cm, clip=true,width=\linewidth]{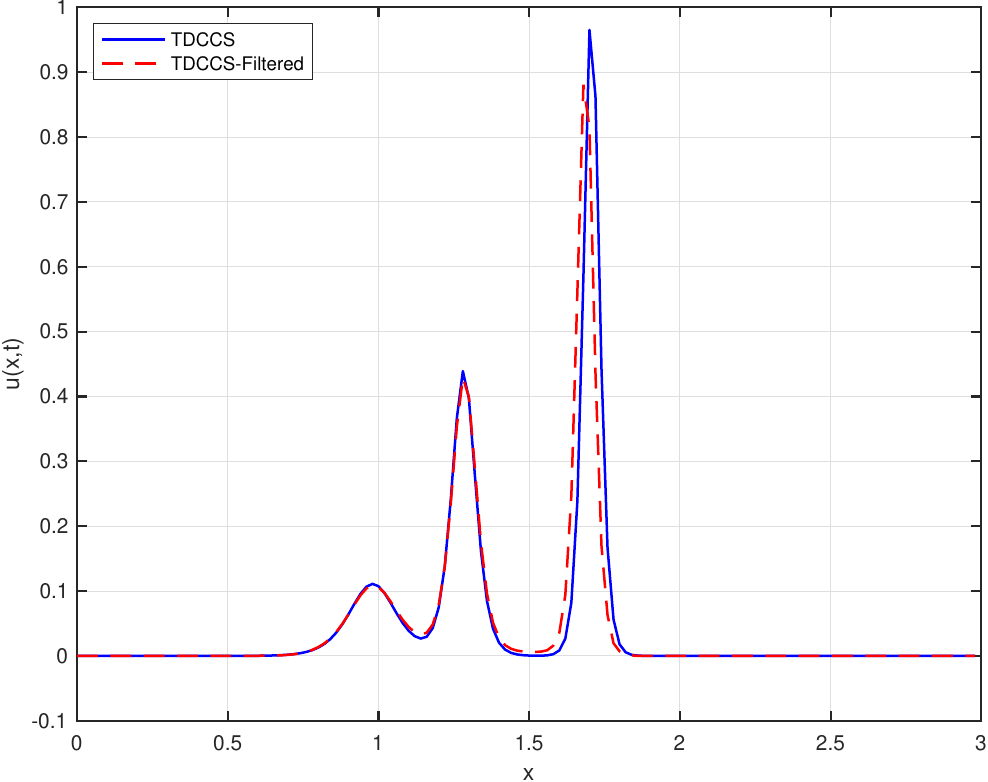}
    \subcaption*{\normalsize{\centering (g) t=2}}
  \end{minipage}\hfill
  \begin{minipage}[b]{0.24\linewidth}
    \centering
    \includegraphics[trim=0cm 0cm 0cm 0cm, clip=true,width=\linewidth]{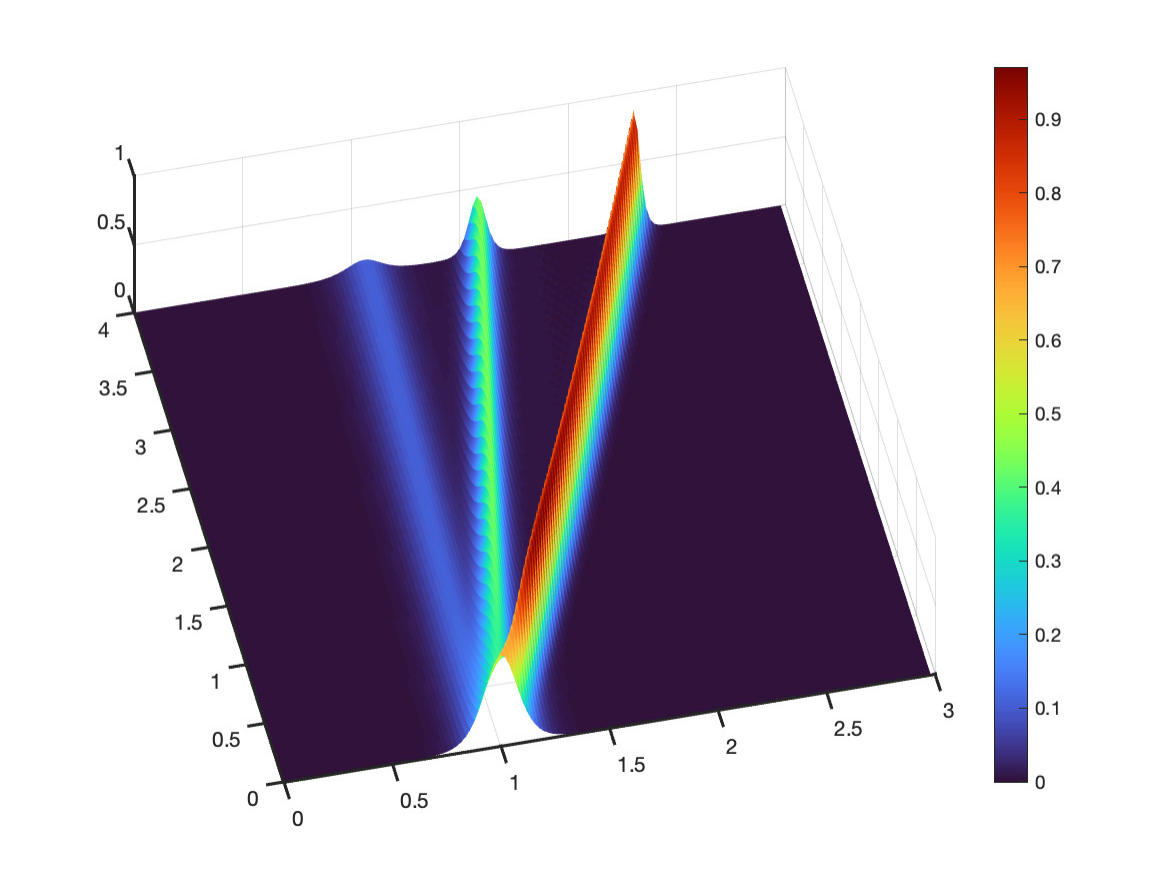}
    \subcaption*{\normalsize{\centering (h) t=4}}
  \end{minipage}\hfill
  \caption{Numerical solutions for
initial condition \ref{IC:3c} of Example \ref{example:3}. The first row of four graphs is computed using TDCNCS with (red) and without (blue) filter and the second row of four graphs is computed using TDCCS with (red) and without (blue) filter.}\label{Figure:E3c}
\end{figure}
\end{example}
\begin{example}\label{example:4}
\normalfont In this example, we investigate the zero dispersion limit of conservation laws, specifically the KdV equation (\ref{example:3}) with the continuous initial condition
 \begin{equation}\label{IC:4}
u_0(x) = 2 + \frac{1}{2}\sin(2 \pi x), \quad x \in[0,1],
 \end{equation}
with periodic boundary conditions. Theoretical and numerical analysis concerning the limit as $\epsilon \to 0^+$ are available in \cite{venakides1987zero} and \cite{Lax}. The primary goal is assessing our numerical method capability to resolve small-scale solution structures in this limit when $\epsilon$ is small. To address this, we compute solutions for $t = 0.5$ with $\epsilon = 10^{-4}, 10^{-5}, 10^{-6}$, and $10^{-7}$ using eighth-order TDCNCS and TDCCS, with a step size $\Delta t = \Delta x^2$ as in \cite{JV}. Since the exact solution is not known for this problem, we utilize a reference solution acquired through a significantly higher number of grid points ($N=1000$) using TDCNCS.  
For $\epsilon = 10^{-4}$ (with $N=100$) and $\epsilon = 10^{-5}$ (with $N=200$), Figures (\ref{Figure:E4a})(a) and (\ref{Figure:E4a})(b) depict the solutions obtained by TDCNCS (in red), TDCCS (in blue), in comparison with the reference solution (in black). Similarly, for $\epsilon = 10^{-6}$ (with $N=800$) and $\epsilon = 10^{-7}$ (with $N=1600$) using both TDCNCS and TDCCS, Figures (\ref{Figure:E4a})(c)-(f) showcase these solutions, confirmed as \enquote{converged}, exhibiting physical oscillations typical in dispersive limits \cite{Lax}. Moreover, inadequate mesh refinement results in failure to achieve a converged solution. For instance, when $\epsilon = 10^{-6}$, we observed that numerical solutions obtained with 200, 300, 400, and even 600 uniform cells did not converge to the solution obtained with 800 cells.  As $\epsilon \to 0^+$, the problem becomes more demanding, necessitating an increasingly fine mesh for convergence to the true solution.  Our method demonstrates high effectiveness in computing such solutions.
 It is observed that using a low-pass filter has no impact on the solution, as the rapidly oscillating physical solution does not constitute high-frequency noise, especially with very fine meshes.
\begin{figure}[htbp!]
  \begin{minipage}[b]{0.3\linewidth}
    \centering
    \includegraphics[trim=0cm 0cm 0cm 0cm, clip=true,width=\linewidth]{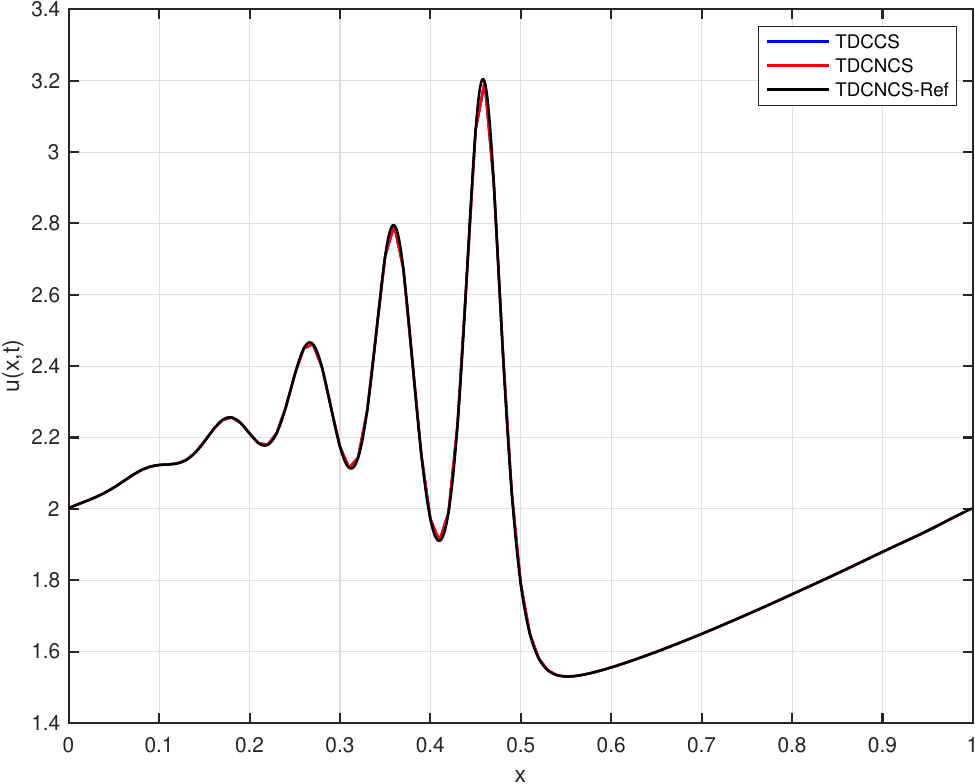}
    \subcaption*{\normalsize{\centering (a) $\epsilon = 10^{-4}$, $N =100$}}
  \end{minipage}\hfill
    \begin{minipage}[b]{0.3\linewidth}
    \centering
    \includegraphics[trim=0cm 0cm 0cm 0cm, clip=true,width=\linewidth]{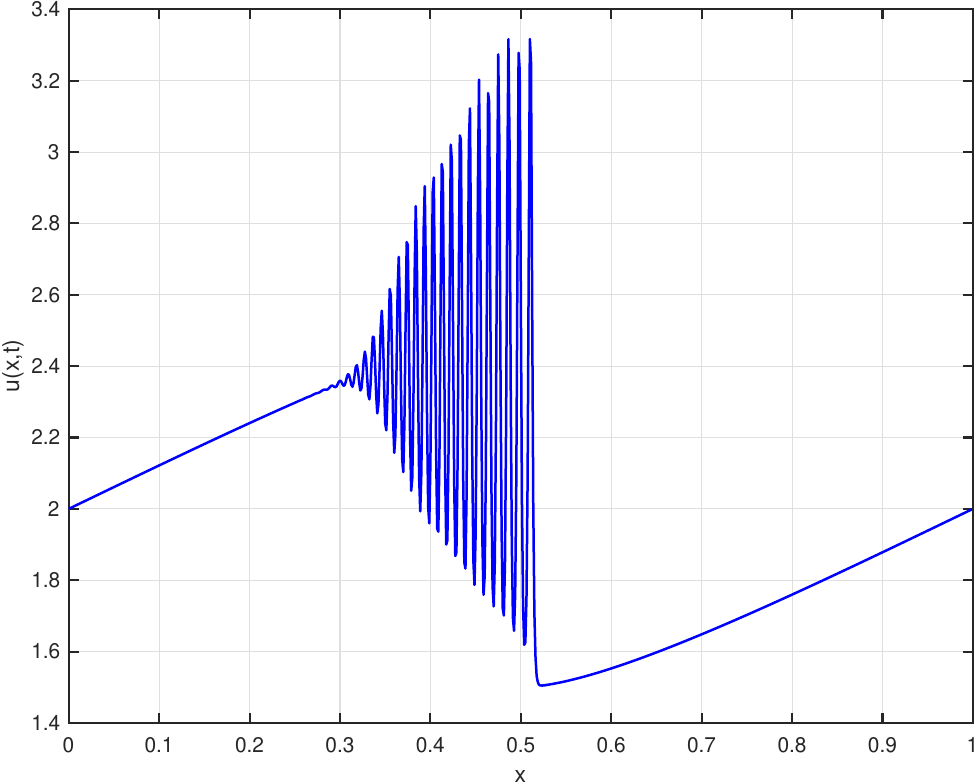}
    \subcaption*{\normalsize{\centering (c) $\epsilon = 10^{-6}$,$N =800$,TDCNCS}}
  \end{minipage}\hfill
  \begin{minipage}[b]{0.3\linewidth}
    \centering
    \includegraphics[trim=0cm 0cm 0cm 0cm, clip=true,width=\linewidth]{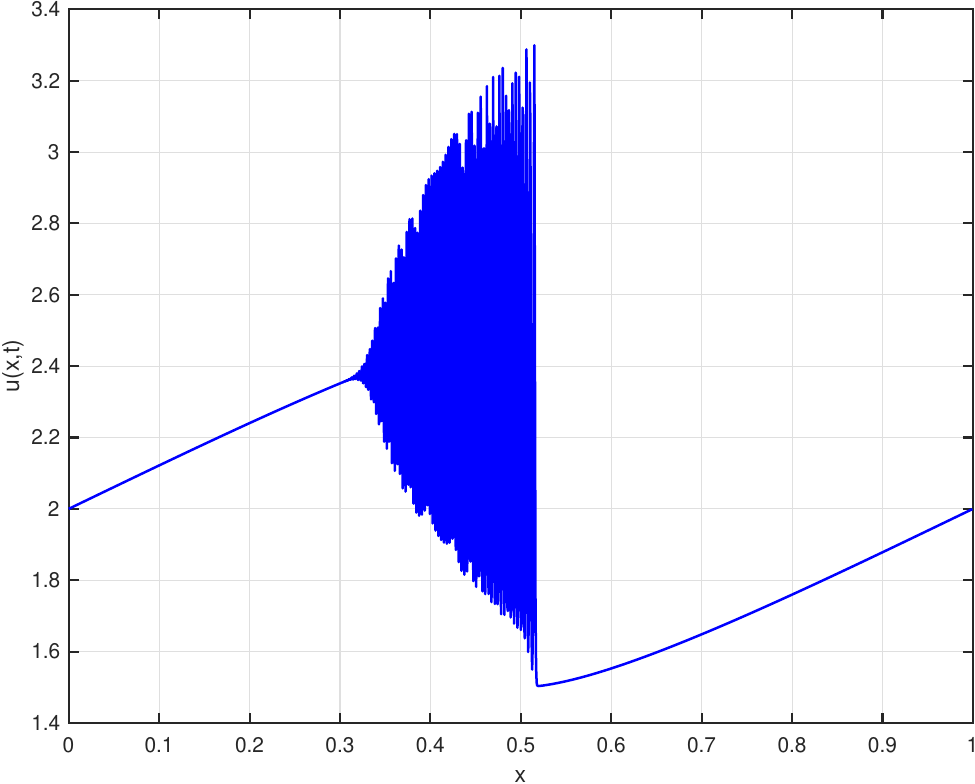}
    \subcaption*{\normalsize{\centering (e)$\epsilon=10^{-7}$,$N=1600$,TDCNCS}}
  \end{minipage}\hfill
    \begin{minipage}[b]{0.3\linewidth}
    \centering
    \includegraphics[trim=0cm 0cm 0cm 0cm, clip=true,width=\linewidth]{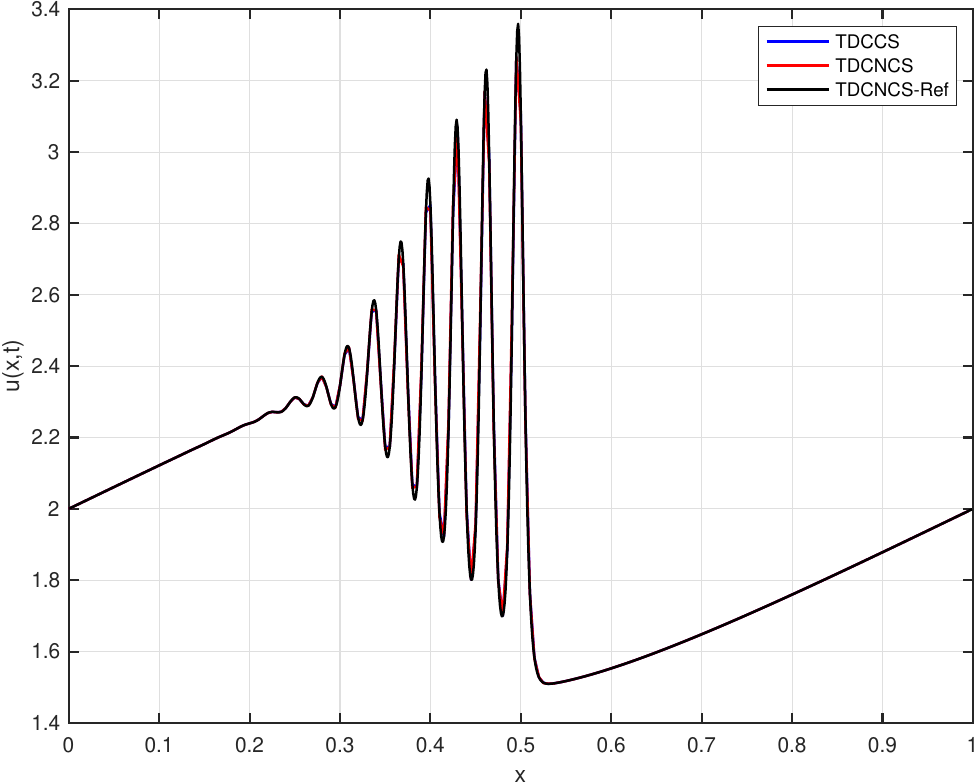}
    \subcaption*{\normalsize{\centering (b)  $\epsilon = 10^{-5}$, $N =200$}}
  \end{minipage}\hfill
    \begin{minipage}[b]{0.3\linewidth}
    \centering
    \includegraphics[trim=0cm 0cm 0cm 0cm, clip=true,width=\linewidth]{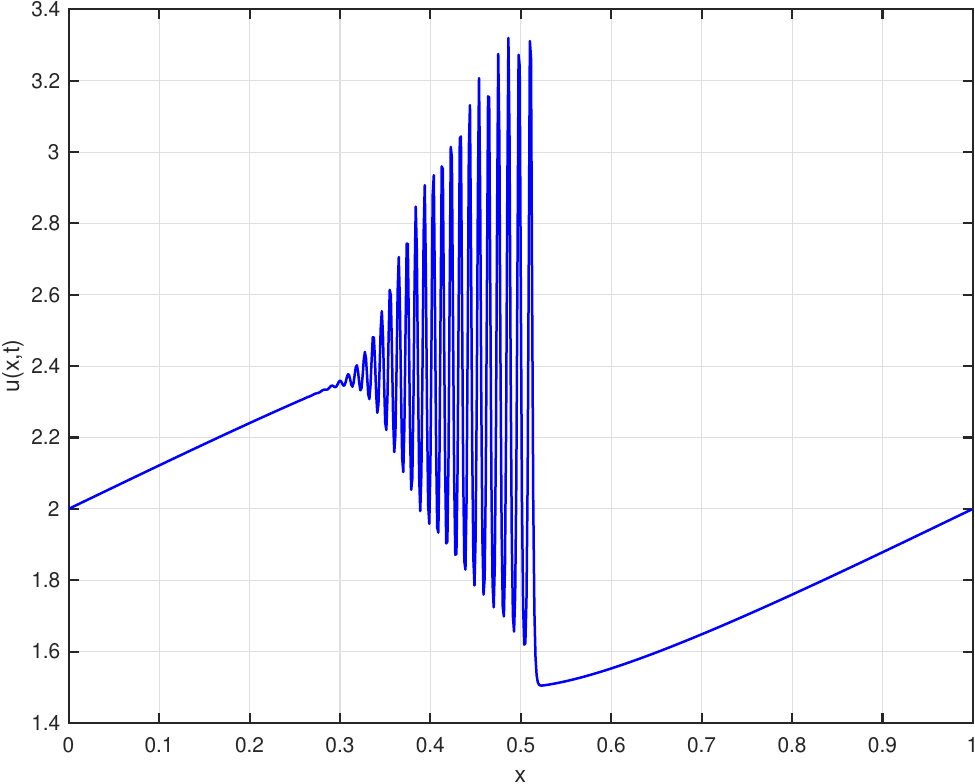}
    \subcaption*{\normalsize{\centering (d) $\epsilon = 10^{-6}$,$N =800$,TDCCS}}
  \end{minipage}\hfill
    \begin{minipage}[b]{0.3\linewidth}
    \centering
    \includegraphics[trim=0cm 0cm 0cm 0cm, clip=true,width=\linewidth]{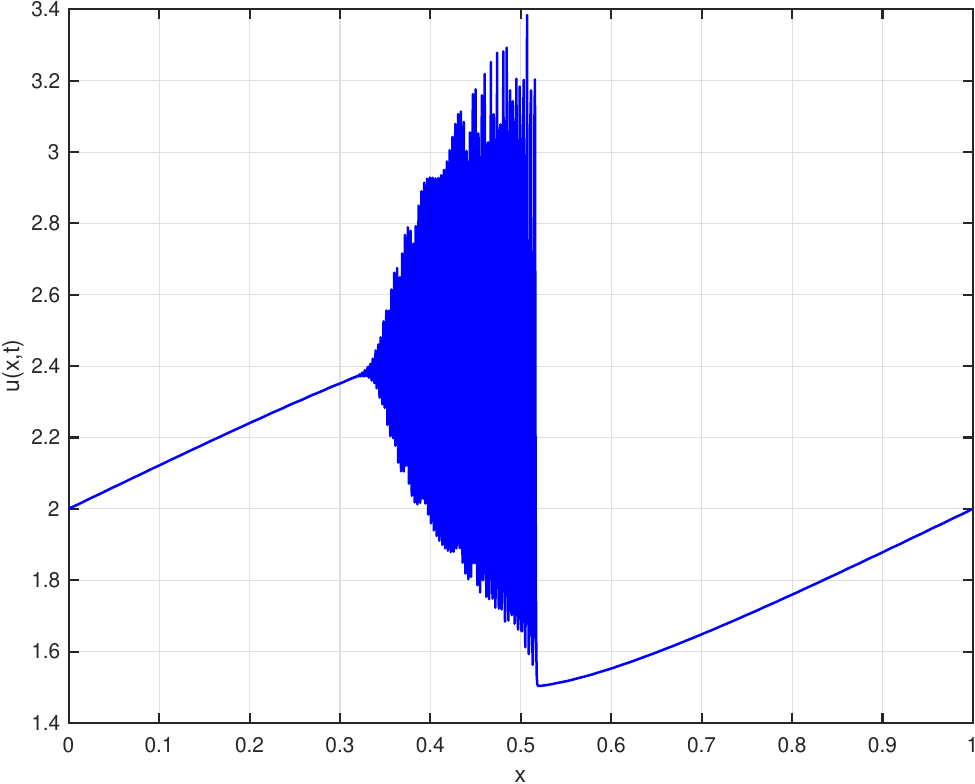}
    \subcaption*{\normalsize{\centering (f)  $\epsilon = 10^{-7}$,$N =1600$,TDCCS}}
  \end{minipage}\hfill
  \caption{Numerical solutions for
initial condition \ref{IC:4} of Example \ref{example:3}.}\label{Figure:E4a}
\end{figure}

\par
Finally, we investigate the behavior of a top hat discontinuous initial condition represented by 
  \begin{equation}
     u(x,0)=\begin{cases}
			1, & \text{if} \quad 0.25<x<4, \\
			0, & \text{else}.
		\end{cases}
  \label{IC:4b}
 \end{equation}
with periodic boundary conditions. We analyze the temporal evolution of a series of left-propagating waves emanating from the initial discontinuity point over multiple time steps. The solution is simulated at $N = 1000$ with a step size $\Delta t = 0.5 \Delta x^2$, $\epsilon = 10^{-4}$ for $t=0.01, 0.05$ respectively, and plotted in Figure (\ref{Figure:E4c}). Within the framework of the TDCNCS scheme, we identify the development of small-amplitude irregularities at locations besides the left side of the discontinuity. However, these irregularities are transient and vanish after the application of a 12th-order filter (every 10 steps) to the aforementioned results (TDCNCS-F). This filtering process leads to the emergence of a well-defined left-propagating wave devoid of such spurious features. Similarly, the TDCCS scheme exhibits less pronounced variations in the solution compared to TDCNCS. These variations are further significantly mitigated by the application of the same filter. The filtered solutions obtained from both TDCCS and TDCNCS are in good agreement with previous studies~\cite{mazaheri2016first, LS}. These solutions depict the transformation of fine-scale continuous waves, as time progresses, into solitary waves within both numerical schemes.

\begin{figure}[htbp!]
  \begin{minipage}[b]{0.25\linewidth}
    \centering
    \includegraphics[trim=0cm 0cm 0cm 0cm, clip=true,width=\linewidth]{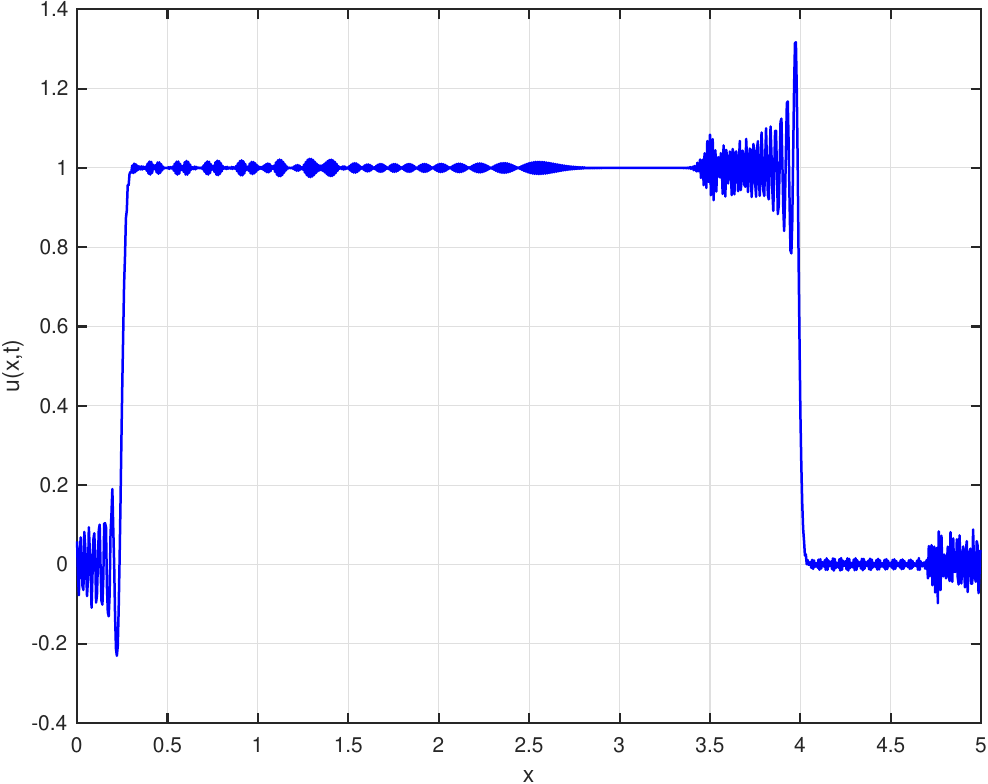}
    \subcaption*{\normalsize{\centering (a)TDCNCS, $t=0.01$}}
  \end{minipage}\hfill
    \begin{minipage}[b]{0.25\linewidth}
    \centering
    \includegraphics[trim=0cm 0cm 0cm 0cm, clip=true,width=\linewidth]{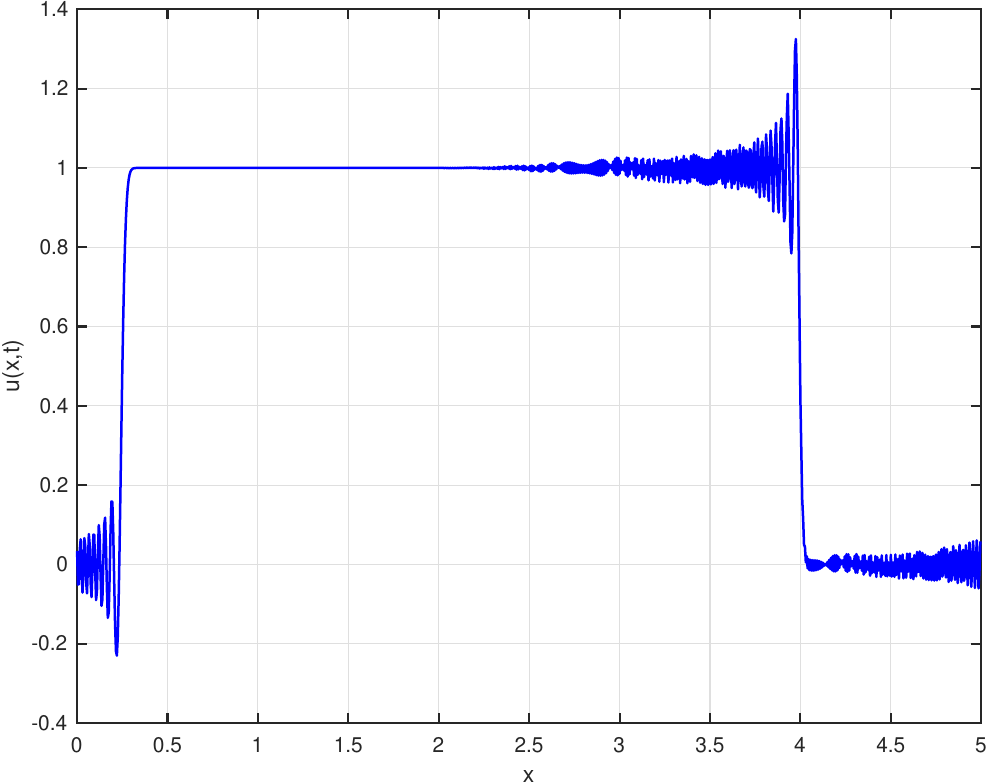}
    \subcaption*{\normalsize{\centering (b) TDCCS, $t=0.01$}}
  \end{minipage}\hfill
    \begin{minipage}[b]{0.25\linewidth}
    \centering
    \includegraphics[trim=0cm 0cm 0cm 0cm, clip=true,width=\linewidth]{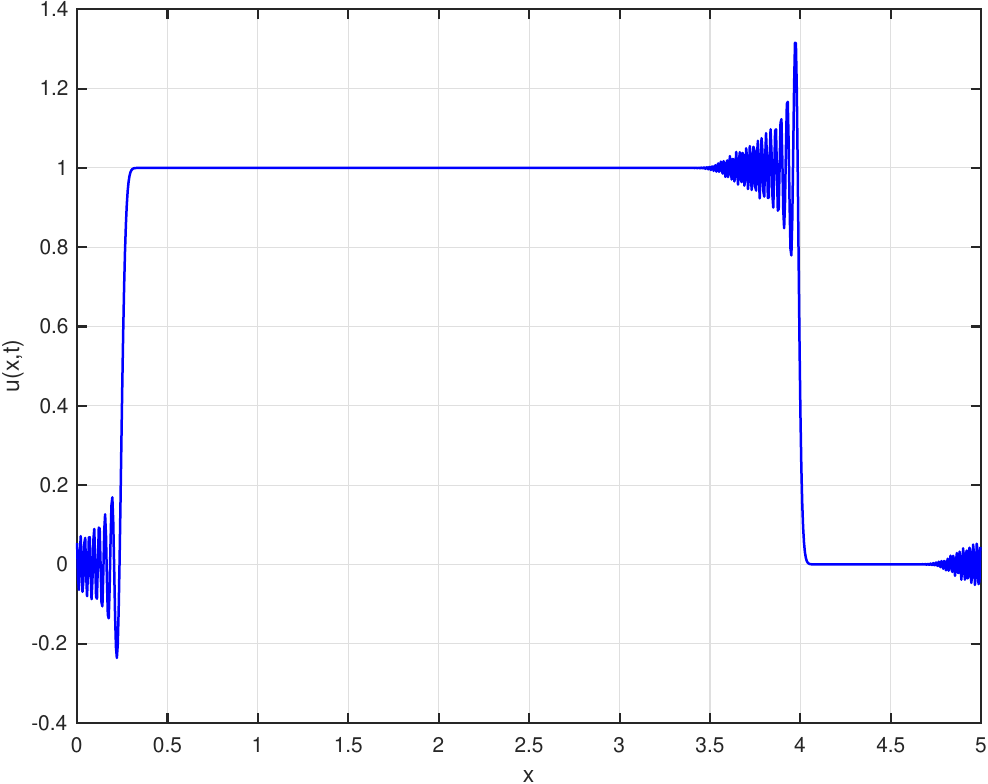}
    \subcaption*{\normalsize{\centering (c)TDCNCS-F, $t=0.01$}}
  \end{minipage}\hfill
  \begin{minipage}[b]{0.25\linewidth}
    \centering
    \includegraphics[trim=0cm 0cm 0cm 0cm, clip=true,width=\linewidth]{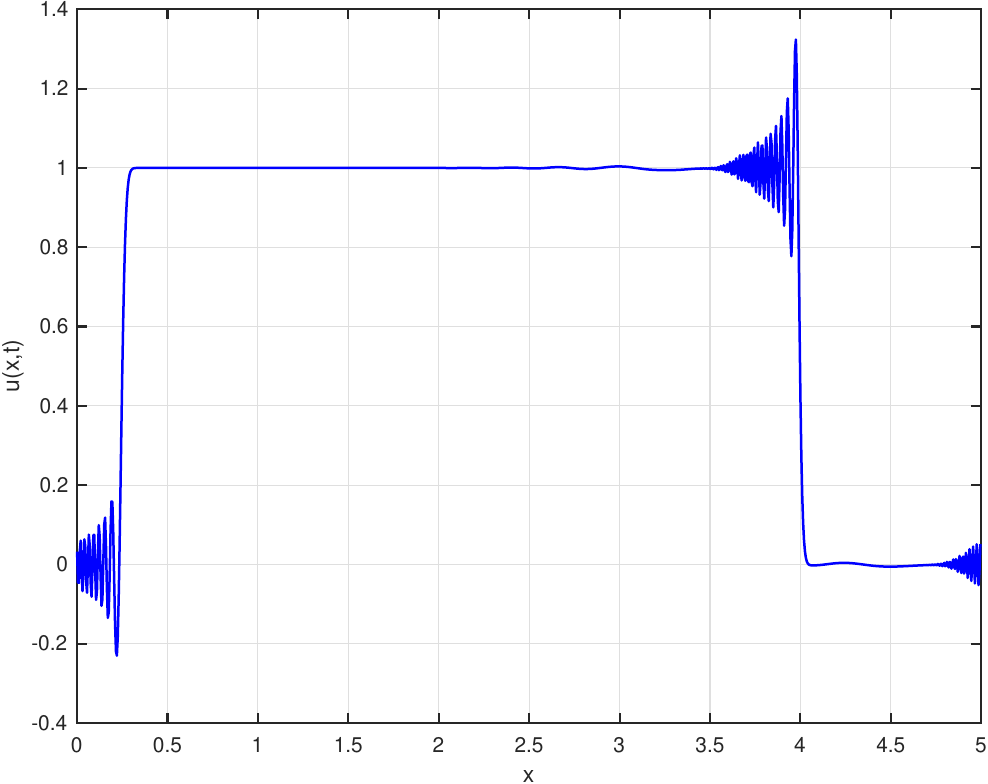}
    \subcaption*{\normalsize{\centering (d) TDCCS-F, $t=0.01$}}
  \end{minipage}\hfill
  \begin{minipage}[b]{0.25\linewidth}
    \centering
    \includegraphics[trim=0cm 0cm 0cm 0cm, clip=true,width=\linewidth]{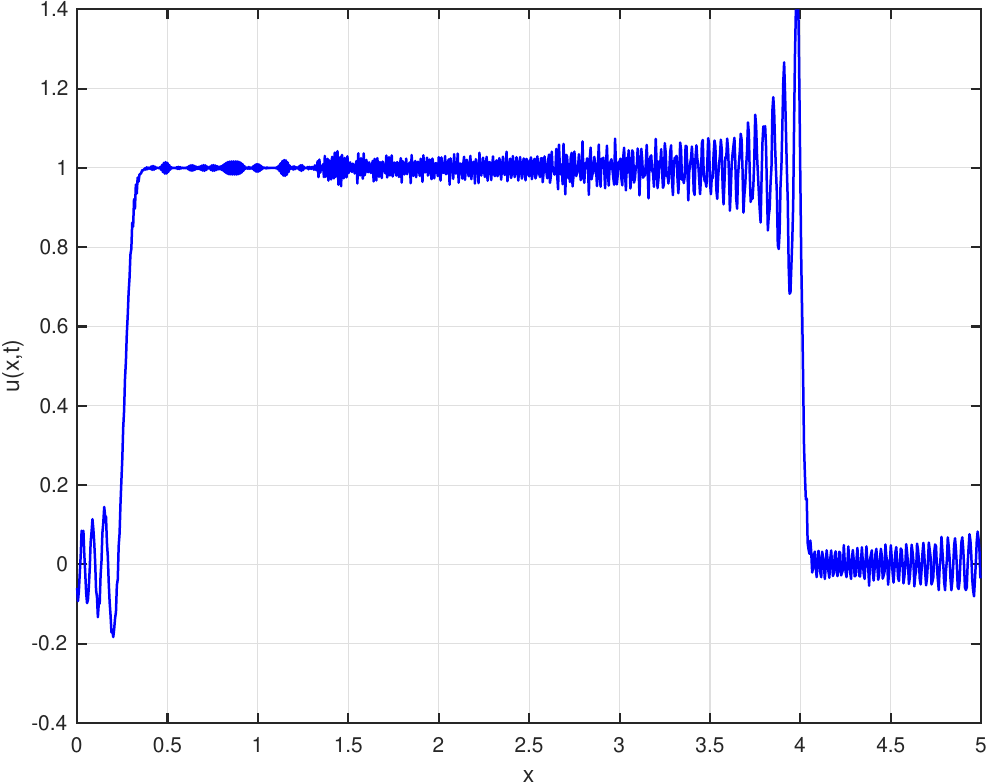}
    \subcaption*{\normalsize{\centering (e)TDCNCS, $t=0.05$}}
  \end{minipage}\hfill
    \begin{minipage}[b]{0.25\linewidth}
    \centering
    \includegraphics[trim=0cm 0cm 0cm 0cm, clip=true,width=\linewidth]{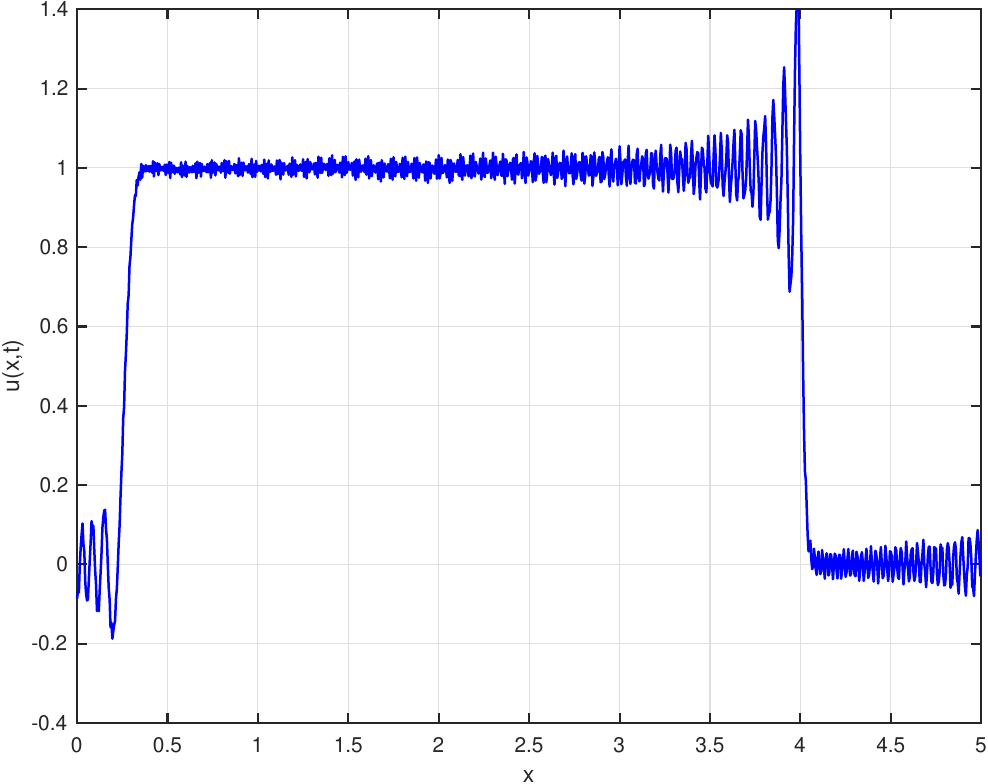}
    \subcaption*{\normalsize{\centering (f) TDCCS, $t=0.05$}}
  \end{minipage}\hfill
    \begin{minipage}[b]{0.25\linewidth}
    \centering
    \includegraphics[trim=0cm 0cm 0cm 0cm, clip=true,width=\linewidth]{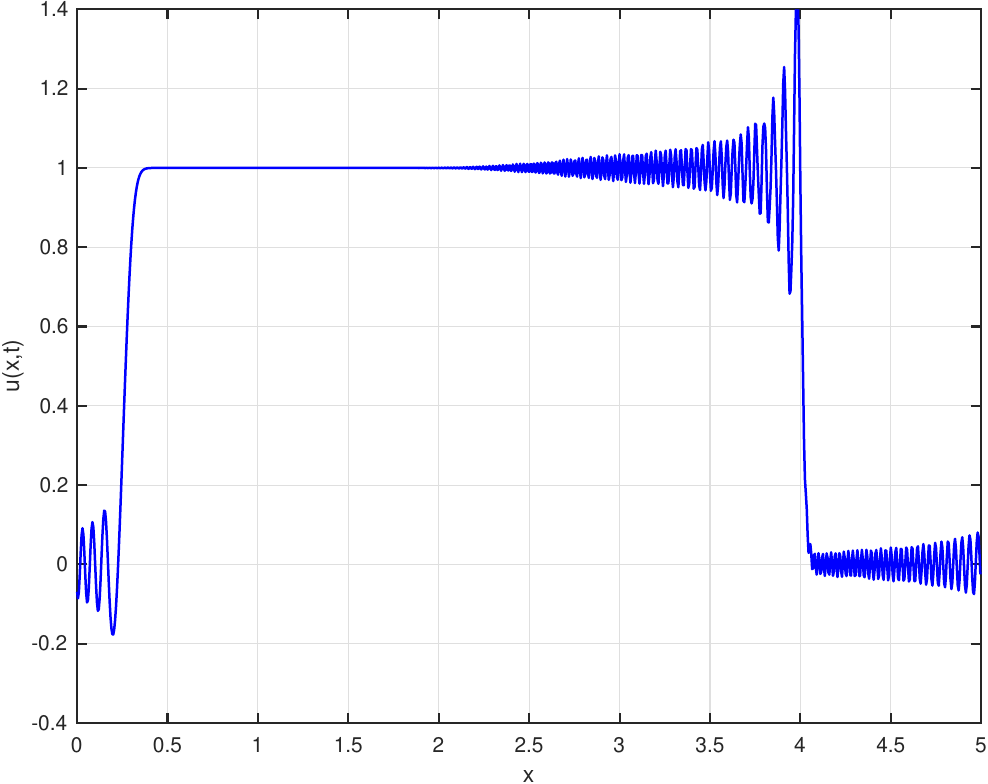}
    \subcaption*{\normalsize{\centering (g)TDCNCS-F, $t=0.05$}}
  \end{minipage}\hfill
  \begin{minipage}[b]{0.25\linewidth}
    \centering
    \includegraphics[trim=0cm 0cm 0cm 0cm, clip=true,width=\linewidth]{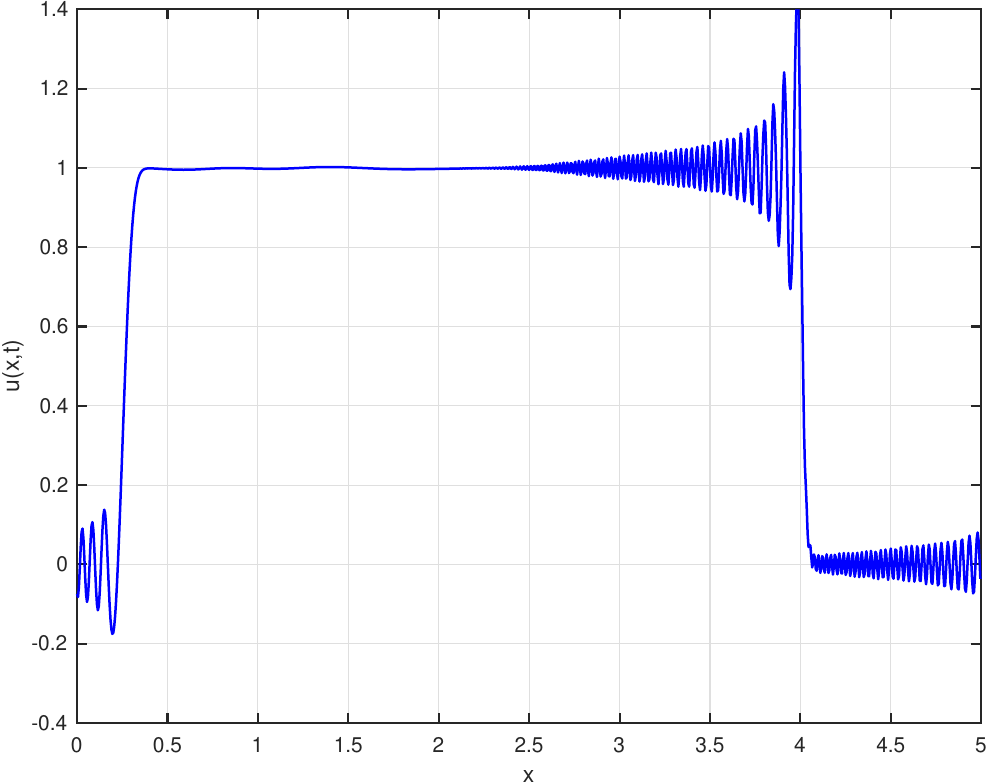}
    \subcaption*{\normalsize{\centering (h) TDCCS-F, $t=0.05$}}
  \end{minipage}\hfill
  \caption{Numerical solutions for
initial condition \ref{IC:4b} of Example \ref{example:3}.}\label{Figure:E4c}
\end{figure}
\end{example}
\section{Conclusion}
\label{sec:con}
In this paper, we introduce a novel family of compact schemes, denoted as third derivative central compact schemes (TDCCS) and third derivative cell-node compact schemes (TDCNCS) extended up to tenth-order accuracy, specifically designed for handling spatial derivatives in dispersion equations. This design is based on the cell-node compact scheme. The TDCCS computes third-order derivatives at the cell-nodes, incorporating values from both the cell-nodes and cell-centers. The values on the cell-centers are determined using the same scheme designed for cell-nodes and are treated as independent variables in the modeling process. This approach has higher memory requirements. Various tests are conducted, encompassing scenarios such as a one-dimensional linear KdV equation, a non-linear KdV equation, and a convection-dominated problem where the coefficients of the third derivative terms are small. A comprehensive comparison is made with third derivative cell-node compact schemes (TDCNCS). This comparative analysis demonstrates the superior accuracy and resolution of the TDCCS over the TDCNCS. In future work, we will study to enhance the proposed method by developing stability-optimized Runge-Kutta techniques to mitigate the constraints imposed by a low CFL number.
\section*{Acknowledgements}
The author Samala Rathan is supported by  NBHM, DAE, India (Ref. No. 02011/46/2021 NBHM(R.P.)/R \& D II/14874). Debojyoti Ghosh contributed to this article under the auspices of the U.S. Department of Energy by Lawrence Livermore National Laboratory under Contract No. DE-AC52-07NA27344.
\section*{Conflict of interest}
The authors declare no potential conflict of interest.
\section*{Data availability}
The data that support the findings of this study are available upon reasonable request.
\appendix
\section{Other possible combinations}
\label{Appendix:A}

\begin{enumerate}
\item \textbf{TDCCS-1:} 
\begin{equation*}
\begin{split}
\beta f''_{j-2}+\alpha f'''_{j-1}+f'''_{j}+\alpha f'''_{j+1}+\beta f''_{j+2} &=  a \frac{4f_{j+1}-8f_{j+\frac{1}{2}}+8f_{j-\frac{1}{2}} -4f_{j-1}}{h^3}+b  \frac{8f_{j+\frac{3}{2}}-12f_{j+1}+12f_{j-1} -8f_{j-\frac{3}{2}}}{5h^3}\\ &+ c  \frac{8f_{j+\frac{5}{2}}-10f_{j+2}+10f_{j-2} -8f_{j-\frac{5}{2}}}{15h^3}\\
\end{split}
\end{equation*}
The modified wavenumber $\omega'''$ with TDCCS-1 is
\begin{equation*}
    \omega'''_{\text{TDCCS-1}} = \frac{2a[8 \sin(\frac{\omega}{2}) - 4 \sin(\omega)]+\frac{2b}{5}[12 \sin(\omega) - 8 \sin(\frac{3\omega}{2})]+\frac{2c}{15}[10 \sin(2\omega) - 8 \sin(\frac{5\omega}{2})]}{1+2\alpha \cos(\omega)+2\beta \cos(2\omega)}
\end{equation*}
Figures \ref{Fig:F_5}(a) and \ref{Fig:F_5}(b) show a plot of modified wavenumber versus wavenumber for third derivative approximations for the Taylor expansion-based method and least square optimization-based method, respectively. Figures \ref{Fig:F_5}(c) and \ref{Fig:F_5}(d) shows the variations of the relative modified wavenumber factor $R_{\omega'''}$ versus wavenumbers for TDCCS-TE-1 and TDCCS-LS-1, respectively.
\begin{figure}[htbp!]
  \begin{minipage}[b]{0.24\linewidth}
    \centering
    \includegraphics[trim=0cm 0cm 0cm 0cm, clip=true,width=\linewidth]{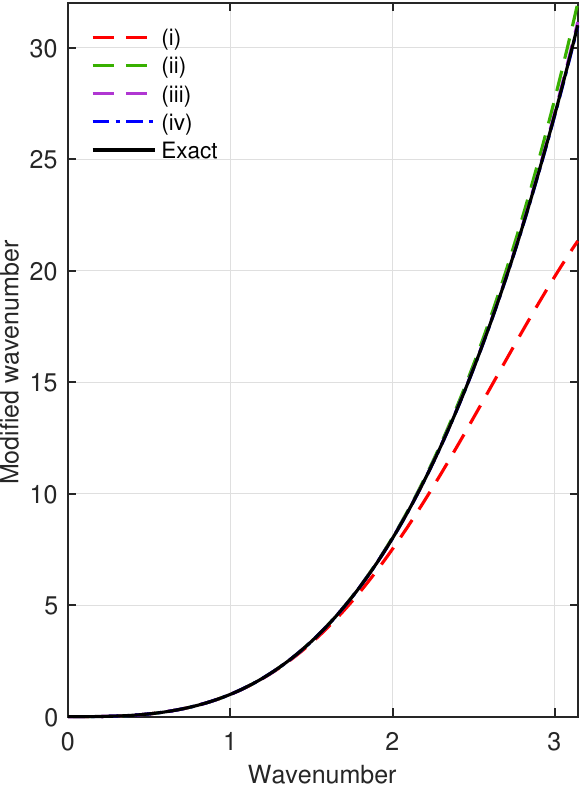}
    \subcaption*{\normalsize{\centering (a) TDCCS-TE-1}}
  \end{minipage}\hfill
    \begin{minipage}[b]{0.24\linewidth}
    \centering
    \includegraphics[trim=0cm 0cm 0cm 0cm, clip=true,width=\linewidth]{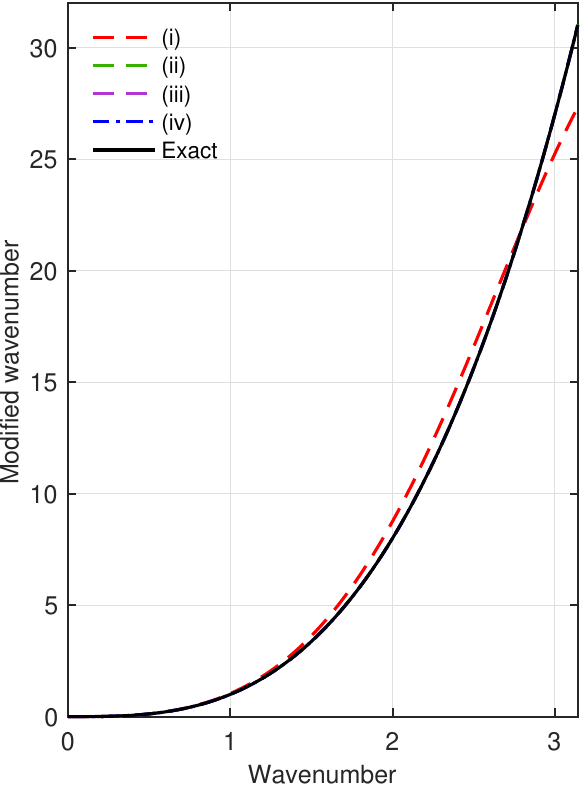}
    \subcaption*{\normalsize{\centering (b) TDCCS-LS-1}}
  \end{minipage}
  \begin{minipage}[b]{0.24\linewidth}
    \centering
    \includegraphics[trim=0.6cm 0cm 0cm 0cm, clip=true,width=\linewidth]{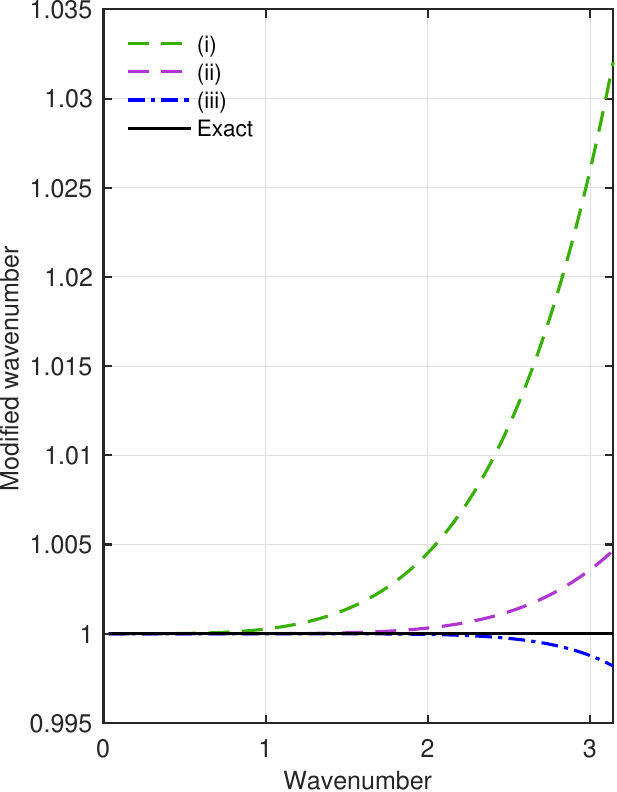}
    \subcaption*{\normalsize{\centering (c) TDCCS-TE-1}}
  \end{minipage}\hfill
    \begin{minipage}[b]{0.24\linewidth}
    \centering
    \includegraphics[trim=0.8cm 0cm 0cm 0cm, clip=true,width=\linewidth]{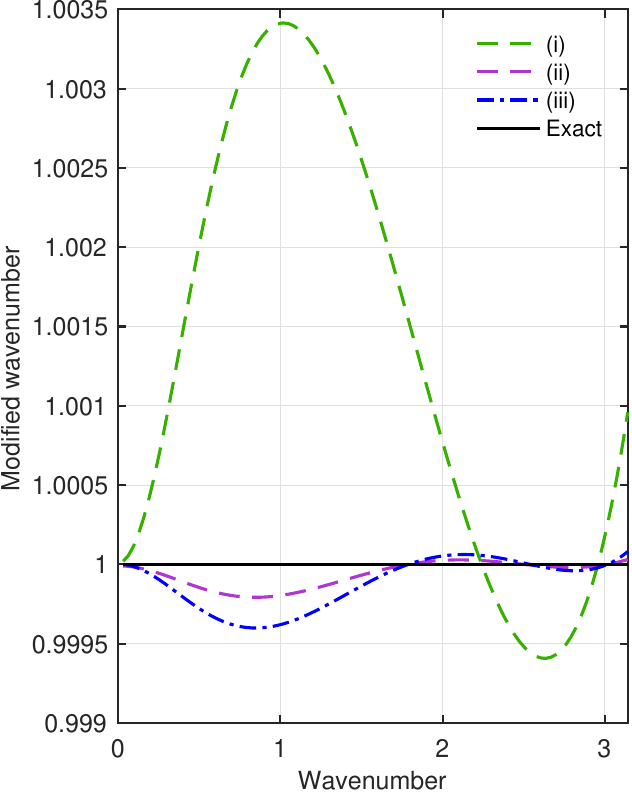}
    \subcaption*{\normalsize{\centering (d) TDCCS-LS-1}}
  \end{minipage}
  \begin{minipage}[b]{0.24\linewidth}
    \centering
    \includegraphics[trim=0cm 0cm 0cm 0cm, clip=true,width=\linewidth]{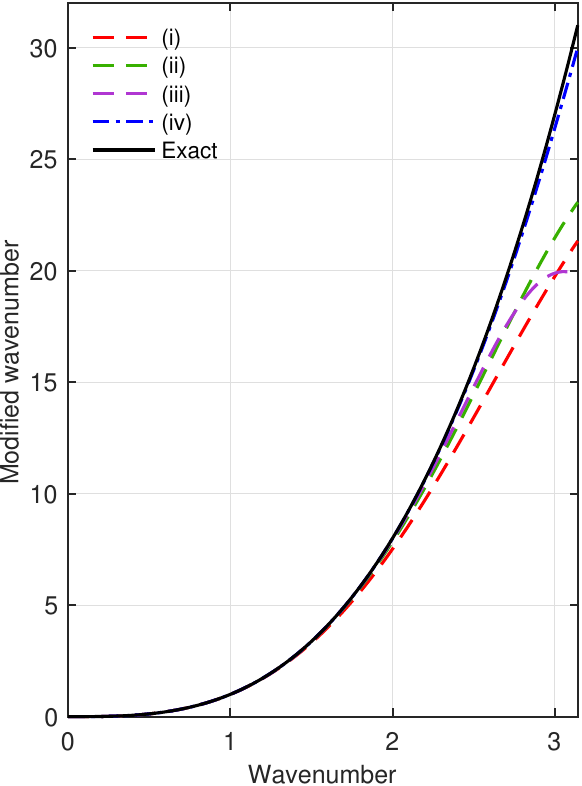}
    \subcaption*{\normalsize{\centering (e) TDCCS-TE-2}}
  \end{minipage}\hfill
    \begin{minipage}[b]{0.24\linewidth}
    \centering
    \includegraphics[trim=0cm 0cm 0cm 0cm, clip=true,width=\linewidth]{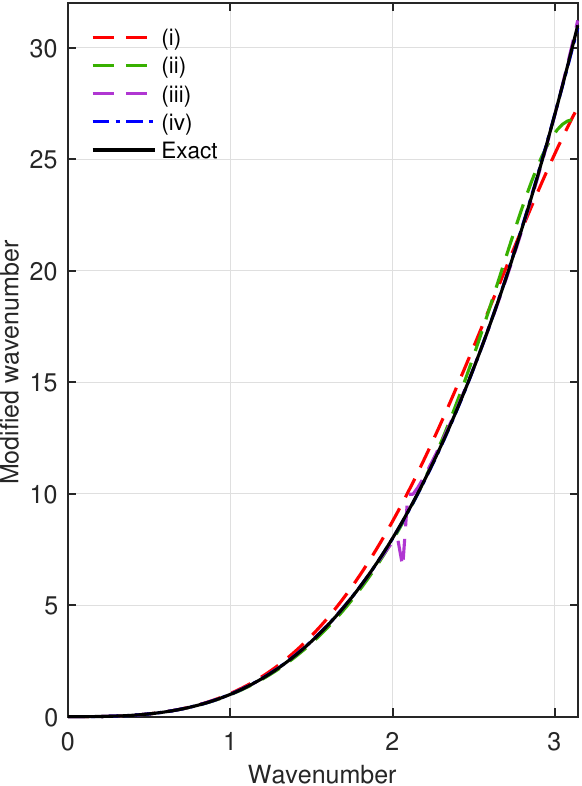}
    \subcaption*{\normalsize{\centering (f) TDCCS-LS-2}}
  \end{minipage}
  \begin{minipage}[b]{0.24\linewidth}
    \centering
    \includegraphics[trim=0.3cm 0cm 0cm 0cm, clip=true,width=\linewidth]{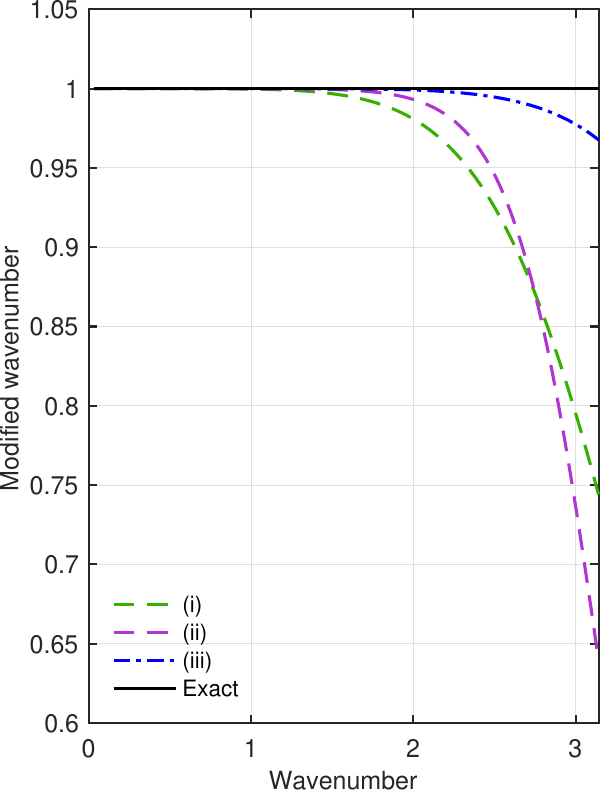}
    \subcaption*{\normalsize{\centering (g) TDCCS-TE-2}}
  \end{minipage}\hfill
    \begin{minipage}[b]{0.24\linewidth}
    \centering
    \includegraphics[trim=0.2cm 0cm 0cm 0cm, clip=true,width=\linewidth]{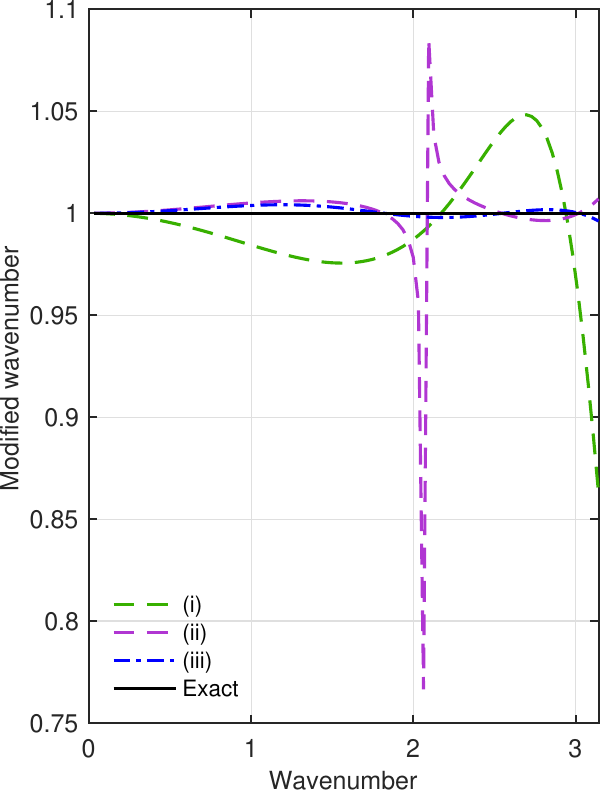}
    \subcaption*{\normalsize{\centering (h) TDCCS-LS-2}}
  \end{minipage}
  \begin{minipage}[b]{0.24\linewidth}
    \centering
    \includegraphics[trim=0cm 0cm 0cm 0cm, clip=true,width=\linewidth]{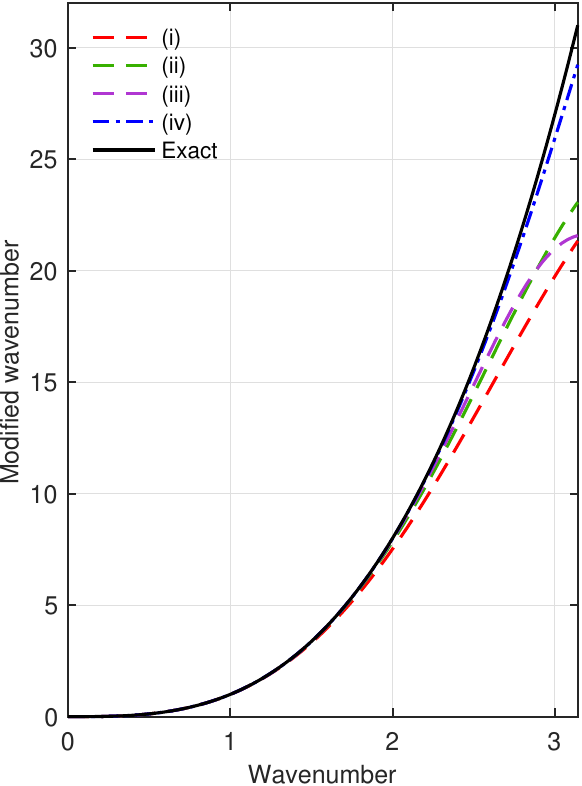}
    \subcaption*{\normalsize{\centering (j) TDCCS-TE-3}}
  \end{minipage}\hfill
    \begin{minipage}[b]{0.24\linewidth}
    \centering
    \includegraphics[trim=0cm 0cm 0cm 0cm, clip=true,width=\linewidth]{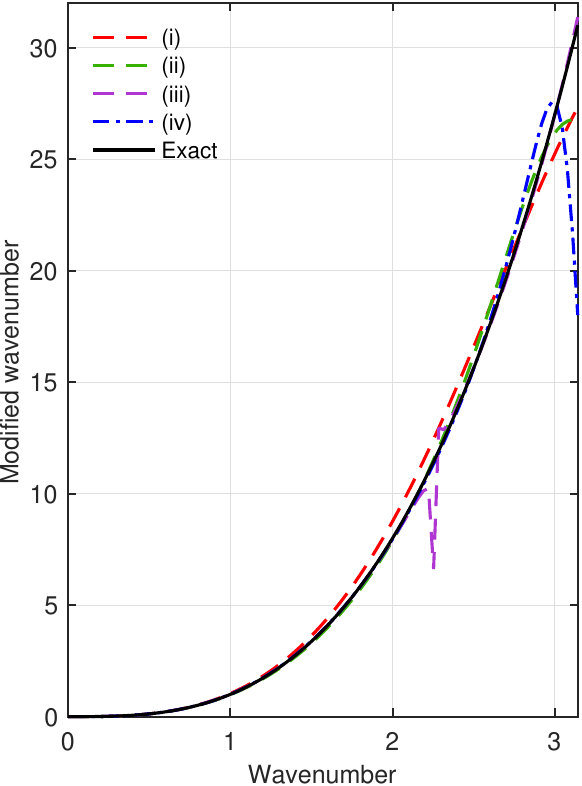}
    \subcaption*{\normalsize{\centering (k) TDCCS-LS-3}}
  \end{minipage}
  \begin{minipage}[b]{0.24\linewidth}
    \centering
    \includegraphics[trim=0.3cm 0cm 0cm 0cm, clip=true,width=\linewidth]{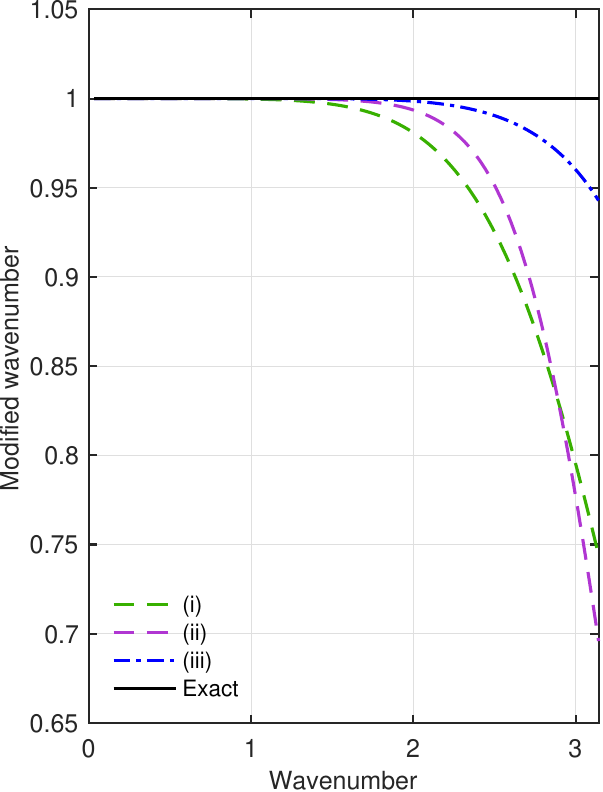}
    \subcaption*{\normalsize{\centering (l) TDCCS-TE-3}}
  \end{minipage}\hfill
    \begin{minipage}[b]{0.24\linewidth}
    \centering
    \includegraphics[trim=0.2cm 0cm 0cm 0cm, clip=true,width=\linewidth]{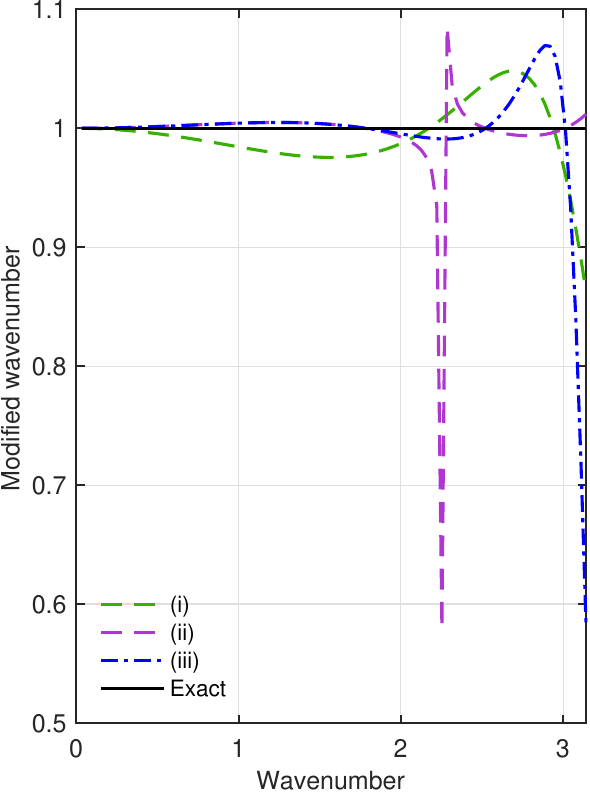}
    \subcaption*{\normalsize{\centering (m) TDCCS-LS-3}}
  \end{minipage}
 \caption{ Plot of modified wavenumber versus wavenumber for (a) TDCCS-TE-1, (b) TDCCS-LS-1, (e) TDCCS-TE-2, (f) TDCCS-LS-2, (j) TDCCS-TE-3 and (k) TDCCS-LS-3: (i) T4; (ii) T6; (iii) T8; (iv) P10. Variations of the relative modified wavenumber factor $R_{\omega'''}$ with wavenumber for (c) TDCCS-TE-1, (d) TDCCS-LS-1, (g) TDCCS-TE-2, (h) TDCCS-LS-2, (l) TDCCS-TE-3 and (m) TDCCS-LS-3:(i) T6; (ii) T8; (iii) P10.}\label{Fig:F_5}
\end{figure}
\item \textbf{TDCCS-2:} 
\begin{equation*}
\begin{split}
\beta f''_{j-2}+\alpha f'''_{j-1}+f'''_{j}+\alpha f'''_{j+1}+\beta f''_{j+2} &=  a \frac{4f_{j+1}-8f_{j+\frac{1}{2}}+8f_{j-\frac{1}{2}} -4f_{j-1}}{h^3}+b  \frac{6f_{j+2}-8f_{j+\frac{3}{2}}+8f_{j-\frac{3}{2}} -6f_{j-2}}{7h^3}\\ &+ c  \frac{8f_{j+\frac{5}{2}}-10f_{j+2}+10f_{j-2} -8f_{j-\frac{5}{2}}}{15h^3}\\
\end{split}
\end{equation*}
The modified wavenumber $\omega'''$ with TDCCS-2 is
\begin{equation*}
    \omega'''_{\text{TDCCS-2}} = \frac{2a[8 \sin(\frac{\omega}{2}) - 4 \sin(\omega)]+\frac{2b}{7}[8 \sin(\frac{3\omega}{2}) - 6 \sin(2\omega)]+\frac{2c}{15}[10 \sin(2\omega) - 8 \sin(\frac{5\omega}{2})]}{1+2\alpha \cos(\omega)+2\beta \cos(2\omega)}
\end{equation*}
Figures \ref{Fig:F_5}(e) and \ref{Fig:F_5}(f) show a plot of modified wavenumber versus wavenumber for third derivative approximations for the Taylor expansion-based method and least square optimization-based method, respectively. Figures \ref{Fig:F_5}(g) and \ref{Fig:F_5}(h) shows the variations of the relative modified wavenumber factor $R_{\omega'''}$ versus wavenumbers for TDCCS-TE-2 and TDCCS-LS-2, respectively.
\item \textbf{TDCCS-3:} 
\begin{equation*}
\begin{split}
\beta f''_{j-2}+\alpha f'''_{j-1}+f'''_{j}+\alpha f'''_{j+1}+\beta f''_{j+2} &=  a \frac{4f_{j+1}-8f_{j+\frac{1}{2}}+8f_{j-\frac{1}{2}} -4f_{j-1}}{h^3}+b   \frac{6f_{j+2}-8f_{j+\frac{3}{2}}+8f_{j-\frac{3}{2}} -6f_{j-2}}{7h^3}\\ &+ c  \frac{8f_{j+\frac{5}{2}}-20f_{j+1}+20f_{j-1} -8f_{j-\frac{5}{2}}}{35h^3}\\
\end{split}
\end{equation*}
The modified wavenumber $\omega'''$ with TDCCS-3 is
\begin{equation*}
    \omega'''_{\text{TDCCS-3}} = \frac{2a[8 \sin(\frac{\omega}{2}) - 4 \sin(\omega)]+\frac{2b}{7}[8 \sin(\frac{3\omega}{2}) - 6 \sin(2\omega)]+\frac{2c}{35}[20 \sin(\omega) - 8 \sin(\frac{5\omega}{2})]}{1+2\alpha \cos(\omega)+2\beta \cos(2\omega)}
\end{equation*}
Figures \ref{Fig:F_5}(j) and \ref{Fig:F_5}(k) show a plot of modified wavenumber versus wavenumber for third derivative approximations for the Taylor expansion-based method and least square optimization-based method, respectively. Figures \ref{Fig:F_5}(l) and \ref{Fig:F_5}(m) shows the variations of the relative modified wavenumber factor $R_{\omega'''}$ versus wavenumbers for TDCCS-TE-3 and TDCCS-LS-3, respectively.
\end{enumerate}

\bibliographystyle{plain}

\end{document}